\documentclass[preprint,11pt]{elsarticle}
\usepackage{amsmath,amssymb,amsthm,mathtools}
\usepackage{graphicx}
\usepackage{booktabs}
\usepackage{algorithm,algpseudocode}
\usepackage{tikz}
\usetikzlibrary{arrows.meta,positioning,fit,backgrounds}
\usepackage[margin=1in]{geometry}
\usepackage[hyphens]{url}
\usepackage[colorlinks=true,linkcolor=blue,citecolor=blue,urlcolor=blue,pdfauthor={Minghui Li, Zixian Liu}]{hyperref}
\graphicspath{{figs/}}

\newtheorem{theorem}{Theorem}[section]
\newtheorem{proposition}[theorem]{Proposition}

\newtheorem{corollary}[theorem]{Corollary}
\theoremstyle{definition}
\newtheorem{definition}[theorem]{Definition}
\theoremstyle{remark}
\newtheorem{remark}[theorem]{Remark}

\newcommand{\fl}{\mathrm{fl}}
\newcommand{\R}{\mathbb{R}}
\newcommand{\Z}{\mathbb{Z}}
\newcommand{\eps}{\varepsilon}
\newcommand{\umach}{\mathbf{u}}
\newcommand{\Fix}{\mathrm{Fix}}
\newcommand{\bsig}{\boldsymbol{\sigma}}
\newcommand{\ee}{\mathrm{e}}

\makeatletter
\def\ps@pprintTitle{%
  \let\@oddhead\@empty
  \let\@evenhead\@empty
  \def\@oddfoot{\hfill\thepage\hfill}%
  \let\@evenfoot\@oddfoot}
\makeatother

\begin{document}
\begin{frontmatter}
\title{Exact long-time symmetry preservation in phase-field simulations:
isotypic spectral schemes that are equivariant in floating-point arithmetic}

\author[a]{Minghui Li\corref{cor}}
\ead{minghuili_xmu@163.com}
\author[a]{Zixian Liu}
\cortext[cor]{Corresponding author.}
\address[a]{School of Mathematical Sciences, Jiangxi Science and Technology Normal University, Nanchang, Jiangxi, China}

\begin{abstract}
Long-time spectral simulations of phase-field models often lose odd or
reflection symmetries and collapse to wrong steady states---even with
high-order, energy-stable integrators. We prove that for grid-equivariant
spectral discretizations this breakdown is not a discretization error, but is
triggered by non-equivariant floating-point rounding in the FFT (and unaligned
shifted sampling), leading to conditional breakdown on a logarithmic time
scale governed by the transverse instability rate and machine precision
whenever an unstable asymmetric mode is excited. Rather than applying
artificial filters or zero-mode resets, we decompose the discrete field into
orthogonal symmetry channels with an exact channel algebra. The same
half-domain mechanism extends directly beyond periodic domains to homogeneous
Neumann and Dirichlet boundary conditions (even without a constant zero mode)
and serves as a modular drop-in backend underneath standard spectral
integrators and scientific-computing workflows (including MATLAB
\texttt{ode15s}, Python/SciPy, and standalone Chebfun/Shenfun-style drivers).
Algebraically identical to standard Fourier, cosine, or sine schemes on the
same collocation grid, it preserves every target subgroup symmetry
\emph{exactly} (zero transverse channels from zero-transverse initialization,
with reduced or zero-skipping modes outperforming standard FFT in favourable
multi-dimensional configurations), resolves genuine sub-roundoff asymmetries
down to $10^{-60}$ in the polynomial channel algebra, and handles off-grid
spatial shifts with double-precision time stepping via translation-invariant
wavenumber parity and a one-time canonical-frame initialization. Benchmarks
across thirteen integrators and 1D/2D/3D Allen--Cahn, Flory--Huggins,
conservative Allen--Cahn, and Cahn--Hilliard models validate exact long-time
symmetry preservation.
\end{abstract}

\begin{keyword}
Allen--Cahn equation \sep symmetry breaking \sep floating-point arithmetic \sep
equivariant discretization \sep discrete sine/cosine transform \sep
exponential integrators \sep long-time simulation
\end{keyword}
\end{frontmatter}

\section{Introduction}\label{sec:intro}

Long-time spectral simulations of phase-field models such as the Allen--Cahn
equation
\begin{equation}\label{eq:AC}
  \partial_t u = \kappa^2\Delta u - \gamma f(u),\qquad f(u)=u^3-u,\qquad \mathbf x\in\Omega=[-\pi,\pi)^d,
\end{equation}
suffer from a striking numerical pathology: even with high-order, energy-stable
time integrators, solutions starting from an odd, shifted-odd, anti-periodic,
or multi-dimensional symmetric initial state routinely lose that symmetry and
collapse to a wrong steady state. Li, Quan, Tang and
Yang~\cite{LQTY2022,LQTY2021sharp} documented this failure across standard
IMEX, convex-splitting, BDF2, SAV and ETD schemes on $u_0=\sin x$ and cured it
with an artificial \emph{parity filter} zeroing forbidden Fourier modes after
every step---which fails on shifted symmetries and erases genuine small
asymmetries. Wang~\cite{Wang2023} later proposed the exponential spectral
process (ESP) method as a filter-free alternative, yet re-testing shows that
ESP breaks at large steps on unshifted 1D data, breaks at Wang's own parameters
on shifted 1D data, and collapses on Wang's 2D
benchmark~\cite[Ex.~4.1]{Wang2023} once continued past $t=20$
(Sections~\ref{sec:num1d}--\ref{sec:numshift}). Likewise, a~posteriori
zero-mode suppression fails in 2D and 3D whenever the unstable asymmetric mode
has zero mean.

\paragraph{The core idea in one line}
The root cause is \emph{not} time-discretization error---every standard
pseudo-spectral scheme is exactly equivariant in exact arithmetic---but tiny
$O(\umach)$ ($\umach\approx 1.1\times 10^{-16}$) asymmetric rounding errors
injected by the FFT, which grow exponentially at rate $\lambda_\perp>0$
whenever an unstable transverse mode is excited. The cure is equally direct.
In 1D, instead of evolving a single mixed array $u$, store its even and odd
components separately as $u=u_{\mathrm e}+u_{\mathrm o}$, apply a cosine
transform to $u_{\mathrm e}$ and a sine transform to $u_{\mathrm o}$, and
evaluate the cubic through the exact identity
\begin{equation}\label{eq:intro1d}
  u^3 = \underbrace{\big(u_{\mathrm e}^3 + 3u_{\mathrm e} u_{\mathrm o}^2\big)}_{\text{even channel}}
      + \underbrace{\big(3u_{\mathrm e}^2 u_{\mathrm o} + u_{\mathrm o}^3\big)}_{\text{odd channel}}.
\end{equation}
When $u$ is odd ($u_{\mathrm e}=0$), every term in the even channel carries an
explicit factor of $0$, so IEEE arithmetic ($0\cdot x=0$, $0+0=0$, ignoring
the sign bit of $\pm 0.0$ and assuming finite intermediates) keeps
$u_{\mathrm e}\equiv 0$ exactly forever without filters or projections; yet
when $u_{\mathrm e}\neq 0$, \eqref{eq:intro1d} is algebraically identical to
standard $(u_{\mathrm e}+u_{\mathrm o})^3$.

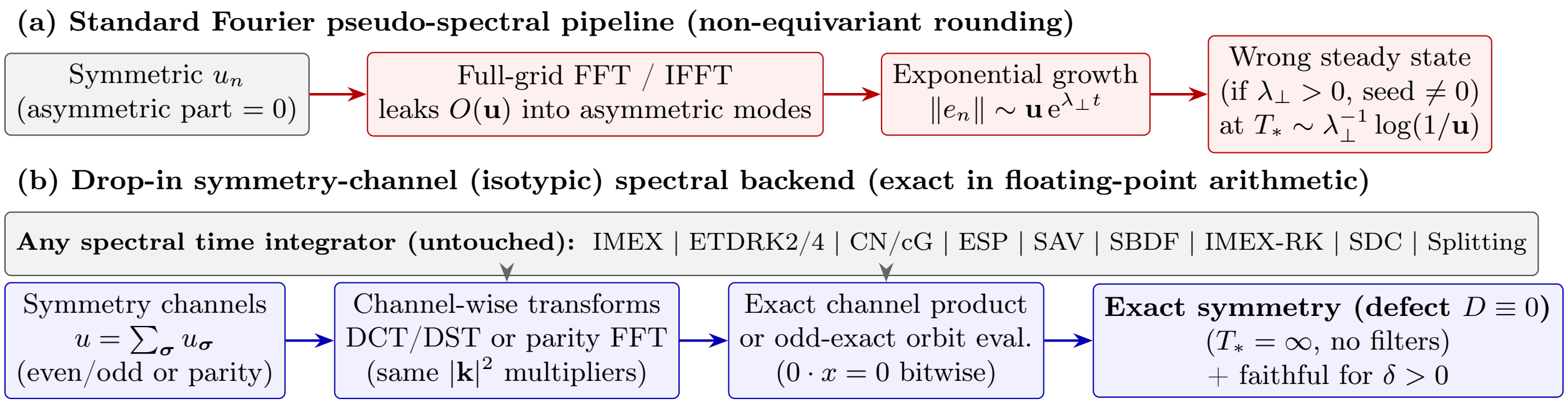
\begin{figure}[t]
\centering
\resizebox{\textwidth}{!}{%
\begin{tikzpicture}[
  font=\footnotesize,
  node distance=5mm and 5.5mm,
  box/.style={draw, rounded corners=2pt, align=center, minimum height=9mm, inner sep=3.5pt},
  bad/.style={box, draw=red!70!black, fill=red!6},
  good/.style={box, draw=blue!70!black, fill=blue!6},
  top/.style={box, draw=black!70, fill=gray!10},
  arr/.style={-{Stealth[length=2.2mm]}, thick}
]
% (a) Standard FFT pipeline
\node[font=\footnotesize\bfseries, anchor=west] (la) at (0, 2.35)
  {(a) Standard Fourier pseudo-spectral pipeline (non-equivariant rounding)};
\node[top, anchor=west] (a1) at (0, 1.55) {Symmetric $u_n$\\(asymmetric part $= 0$)};
\node[bad, right=6.5mm of a1] (a2) {Full-grid FFT / IFFT\\leaks $O(\umach)$ into asymmetric modes};
\node[bad, right=6.5mm of a2] (a3) {Exponential growth\\$\|e_n\|\sim \umach\,\ee^{\lambda_\perp t}$};
\node[bad, right=6.5mm of a3] (a4) {Wrong steady state\\(if $\lambda_\perp>0$, seed $\neq 0$)\\at $T_\ast\sim\lambda_\perp^{-1}\log(1/\umach)$};
\draw[arr, red!70!black] (a1) -- (a2);
\draw[arr, red!70!black] (a2) -- (a3);
\draw[arr, red!70!black] (a3) -- (a4);

% (b) Drop-in Isotypic Spectral Backend
\node[font=\footnotesize\bfseries, anchor=west] (lb) at (0, 0.55)
  {(b) Drop-in symmetry-channel (isotypic) spectral backend (exact in floating-point arithmetic)};
\node[good, anchor=west] (b1) at (0, -1.25) {Symmetry channels\\$u=\sum_{\bsig} u_{\bsig}$\\(even/odd or parity)};
\node[good, right=5.5mm of b1] (b2) {Channel-wise transforms\\DCT/DST or parity FFT\\(same $|\mathbf k|^2$ multipliers)};
\node[good, right=5.5mm of b2] (b3) {Exact channel product\\or odd-exact orbit eval.\\($0\cdot x=0$ bitwise)};
\node[good, right=5.5mm of b3] (b4) {\textbf{Exact symmetry (defect $D\equiv 0$)}\\($T_\ast=\infty$, no filters)\\+ faithful for $\delta>0$};
\node[top, anchor=west, minimum width=15.05cm, minimum height=7.2mm, font=\scriptsize] (btop) at (0, -0.15)
  {\textbf{Any spectral time integrator (untouched):}\ \ IMEX \textbar\ ETDRK2/4 \textbar\ CN/cG \textbar\ ESP \textbar\ SAV \textbar\ SBDF \textbar\ IMEX-RK \textbar\ SDC \textbar\ Splitting};
\draw[arr, blue!70!black] (b1) -- (b2);
\draw[arr, blue!70!black] (b2) -- (b3);
\draw[arr, blue!70!black] (b3) -- (b4);
\draw[arr, black!60, dashed] (btop.south -| b2.north) -- (b2.north);
\draw[arr, black!60, dashed] (btop.south -| b3.north) -- (b3.north);
\end{tikzpicture}}
\caption{Conceptual comparison between (a)~the standard Fourier pseudo-spectral
pipeline, where non-equivariant FFT rounding seeds unstable asymmetric modes
regardless of time-integrator order, and (b)~the proposed drop-in
symmetry-channel (isotypic) backend, which isolates symmetry components
algebraically underneath any spectral time integrator with zero artificial
intervention.}\label{fig:overview}
\end{figure}

\paragraph{Main contributions and broader impact}
Extending~\eqref{eq:intro1d} into a unified $d$-dimensional \emph{isotypic
spectral framework} (Figure~\ref{fig:overview}) yields four key advances:
\begin{enumerate}
\item \textbf{Floating-point diagnosis and conditional breakdown law:}
  Standard Fourier schemes are exactly equivariant in exact arithmetic
  (Proposition~\ref{prop:equiv}). When rounding excites an unstable transverse
  mode under the stated hypotheses, breakdown occurs at
  $T_\ast=\lambda_\perp^{-1}\log(\eta/\umach)+O(1)$ (threshold $\eta$, e.g.\
  $10^{-2}$; Theorem~\ref{thm:breakdown})---compounded for unaligned shifted
  data by an $O(\umach)$ initial-sampling obstruction
  (Proposition~\ref{prop:obstruction}). Higher temporal order or smaller step
  sizes cannot remove this mechanism.
\item \textbf{Exact zero-channel invariance, polynomial relative faithfulness,
  and non-periodic boundary conditions:} Splitting the field into $2^d$
  orthogonal reflection ($A_{\mathrm R}$, via DCT/DST) or wavenumber-parity
  ($\mathcal T_2$) channels is algebraically identical to the standard scheme
  on the same collocation grid (Theorem~\ref{thm:identical}), keeps every
  subgroup subspace of $A_{\mathrm R}$ or $\mathcal T_2$ exactly zero-invariant
  from zero-transverse initial data ($D\equiv 0$, Theorem~\ref{thm:exact},
  including general odd potentials via odd-exact orbit evaluation,
  Proposition~\ref{prop:orbit}), and for the polynomial channel
  algebra~\eqref{eq:chalg} propagates genuine small asymmetries with full
  relative accuracy down to $\delta=10^{-60}$ against a 120-digit reference
  (Theorem~\ref{thm:relative}, Corollary~\ref{cor:faithful}). Crucially, the
  exact same half-domain channel splitting on $[0,\pi/2]^d$ applies
  \emph{beyond periodic domains} to \textbf{homogeneous Neumann and Dirichlet
  boundary conditions} on $[0,\pi]^d$, where standard DCT/DST solvers likewise
  suffer midpoint symmetry collapse---even in 1D Dirichlet problems with no
  $k=0$ constant mode (Remark~\ref{rem:ext}).
\item \textbf{Arbitrary off-grid shifts in double-precision time stepping:}
  Half-period parity ($\mathcal T_2$) commutes with every spatial shift and
  preserves anti-periodic (odd-harmonic) symmetries exactly in pure IEEE double
  precision from parity-consistent initialization (Theorem~\ref{thm:t2}). For
  mixed even-and-odd harmonics under unknown off-grid shifts, our
  \emph{canonical-frame isotypic spectral method} (CF-ISM) performs a one-time
  $d_0$-digit gauge initialization at $t=0$ and executes 100\% of the $t>0$
  time stepping in hardware double precision, extending the symmetry horizon
  as $\lambda_\perp^{-1}d_0\ln 10+O(1)$ and attaining exact zero transverse
  channels once $C_0 10^{-d_0}<2^{-1075}$ (Theorem~\ref{thm:cfism}).
\item \textbf{A modular drop-in backend for solvers and scientific libraries:}
  Because spectral PDE drivers decouple diagonal linear multipliers
  $m(|\mathbf k|^2)$ from transforms and pointwise nonlinearities, replacing
  those low-level primitives (plus orthogonal channel inner products for SAV
  and symmetry-preserving stopping norms for iterative/adaptive solvers)
  upgrades existing time steppers---validated across 13 schemes (IMEX1,
  ETDRK2/4, CN, ESP, SAV1, BDF2-SAV, ARS(4,4,3), SBDF3/4, IFRK4, SDC4,
  Strang2) plus Flory--Huggins, conservative Allen--Cahn and Cahn--Hilliard,
  alongside drop-in demos for MATLAB \texttt{ode15s}, Python/SciPy,
  Chebfun-style \texttt{spin}/\texttt{spin2}, Shenfun-style ETDRK2, and a
  Dedalus~v3 recipe (Section~\ref{sec:ext}). While full-channel MATLAB
  prototypes carry a constant factor ($1.78$--$5.65\times$, Table~\ref{tab:cost}),
  exact-zero-channel short-circuiting or reduced symmetric mode outperforms
  standard FFT in favourable multi-dimensional configurations.
\end{enumerate}

\noindent Symmetry-adapted (isotypic)
decompositions are classical in numerical bifurcation analysis and boundary- or
finite-element solvers, where they decouple linear systems and trace symmetric
branches~\cite{Bossavit1986,Dellnitz1992,AllgowerEtAl1992,GSS1988,GolubitskyStewart2002};
DCT/DST are the reflection transforms of periodic sequences~\cite{Martucci1994},
and fixing a moving reference origin is standard in equivariant dynamical
systems (the method of slices~\cite{RowleyMarsden2000,BeynThuemmler2004}).
Here we apply these tools to a different goal: making the \emph{floating-point}
realization of spectral time-stepping schemes symmetry-preserving without
intervention. Long-time rounding sensitivity has also been highlighted for
chaotic systems in Clean Numerical Simulation (CNS)~\cite{HuLiao2020}, which
executes every time step in multi-precision arithmetic; by contrast, CF-ISM
confines extended precision to $t=0$ only and executes 100\% of the $t>0$ time
stepping in hardware double precision. The stabilized
exponential, IMEX, SAV and deferred-correction integrators tested here are
standard~\cite{CoxMatthews2002,KassamTrefethen2005,DJLQ2021,FuYang2022,TangYang2016,ShenXuYang2019,FengTangYang2015,AscherRuuthWetton1995,AscherRuuthSpiteri1997,DuttGreengardRokhlin2000},
as are the rounding-error models~\cite{Higham2002} and spectral
methods~\cite{Trefethen2000,ShenTangWang2011,FrigoJohnson2005}.

The paper is organized as follows. Section~\ref{sec:setting} sets up the
discrete symmetry framework and proves exact-arithmetic equivariance.
Section~\ref{sec:fp} analyses floating-point symmetry breaking and the
initial-data obstruction. Section~\ref{sec:iso} presents the reflection and
half-period parity channel schemes together with CF-ISM and proves their
properties. Sections~\ref{sec:num1d} and~\ref{sec:num23d} present centred
one-, two- and three-dimensional experiments (including thirteen time
integrators, Wang's 2D benchmark and the Flory--Huggins potential).
Section~\ref{sec:numshift} treats arbitrarily shifted symmetries in 1D, 2D and
3D as well as drop-in extensions to conservative Allen--Cahn and
Cahn--Hilliard equations, and Section~\ref{sec:concl} concludes. Detailed
proofs are collected in the appendices.

%======================================================================
\section{Setting and exact-arithmetic equivariance}\label{sec:setting}

\subsection{Symmetry group and invariant subspaces}
Let $N$ be even, $M=N/2$,
$h=2\pi/N$ and let $\mathcal G_N=\{\mathbf x_{\mathbf j}=h\mathbf j:\
\mathbf j\in\{0,\dots,N-1\}^d\}$ be the periodic grid on the torus
$\mathbb T^d=[-\pi,\pi)^d$. The continuous PDE \eqref{eq:AC} is equivariant under the
torus symmetry group $G=(\mathbb T^d\rtimes\Z_2^d)\times\Z_2$, where
$(T_{\mathbf a}u)(\mathbf x)=u(\mathbf x-\mathbf a)$,
$(R_iu)(\mathbf x)=u(x_1,\dots,-x_i,\dots,x_d)$, and $s\in\{\pm1\}$ acts by
sign flip $u\mapsto s\,u$.

Let $A\le\mathbb T^d\rtimes\Z_2^d$ be a finite Abelian subgroup generated by
$d$ independent commuting involutions $\tau_1,\dots,\tau_d$ (so $|A|=2^d$)
that map $\mathcal G_N$ onto itself by index permutations (hence act
\emph{exactly} in floating-point arithmetic) and preserve $|\mathbf k|^2$ in
Fourier space. Our two primary choices are:
\begin{itemize}
\item the \emph{reflection group} $A_{\mathrm R}=\langle R_1,\dots,R_d\rangle\cong\Z_2^d$,
  with $R_i:j_i\mapsto (N-j_i)\bmod N$ on $\mathcal G_N$ (reflection about $0$),
  or the half-grid reflection $R_i^{(1/2)}:j_i\mapsto N-1-j_i$, which reflects
  about $0\pmod{2\pi}$ on the shifted grid $x_{j_i}=(j_i+\tfrac12)h$ (and
  about $-h/2$ on $\mathcal G_N$, whereas reflection about $+h/2$ on
  $\mathcal G_N$ is $j_i\mapsto (1-j_i)\bmod N$);
\item the \emph{half-period translation group}
  $A_{\mathrm T}=\mathcal T_2=\langle T_{\pi\mathbf e_1},\dots,T_{\pi\mathbf e_d}\rangle\cong\Z_2^d$
  (where $\mathbf e_i$ is the $i$-th standard basis vector and we take $N$
  divisible by $4$ so that $M=N/2$ is even), with
  $T_{\pi\mathbf e_i}:j_i\mapsto (j_i+M)\bmod N$.
\end{itemize}
Set $\Gamma_A=A\times\Z_2$ acting by $(g,s)\cdot u=s\,u\circ g$. For
$\bsig\in\Z_2^d$ the \emph{character}
$\chi_{\bsig}(g)=\prod_{i=1}^d(-1)^{\sigma_ig_i}$ ($g=\tau_1^{g_1}\cdots\tau_d^{g_d}$)
defines the orthogonal isotypic projection
\begin{equation}\label{eq:isoproj}
  \Pi_{\bsig}u = 2^{-d}\sum_{g\in A}\chi_{\bsig}(g)\,u\circ g,
  \qquad u=\sum_{\bsig\in\Z_2^d}u_{\bsig},\qquad u_{\bsig}:=\Pi_{\bsig}u .
\end{equation}
For $A=A_{\mathrm R}$, $u_{\bsig}$ is even in $x_i$ if $\sigma_i=0$ and odd if
$\sigma_i=1$; for $A=\mathcal T_2$, $u_{\bsig}$ is $\pi$-periodic in $x_i$
(even wavenumbers $k_i\in 2\Z$) if $\sigma_i=0$ and $\pi$-anti-periodic
(odd wavenumbers $k_i\in 2\Z+1$) if $\sigma_i=1$. We call $u_{\bsig}$ the
$\bsig$-\emph{channel}. For any subgroup $H\le\Gamma_A$ the fixed-point
subspace $\Fix(H)=\{u:\ h\cdot u=u\ \forall h\in H\}$ is a sum of channels,
\begin{equation}\label{eq:FixH}
  \Fix(H)=V_{\Sigma_H}:=\{u:\ u_{\bsig}=0\ \ \forall\bsig\notin\Sigma_H\},\qquad
  \Sigma_H=\{\bsig\in\Z_2^d:\ \chi_{\bsig}(g)=s\ \ \forall (g,s)\in H\}.
\end{equation}
Examples include: 1D odd functions ($A_{\mathrm R}$, $\Sigma_H=\{1\}$);
$d$-dimensional functions odd in every variable ($A_{\mathrm R}$,
$\Sigma_H=\{(1,\dots,1)\}$); 2D point-odd ($C_2$) functions $u(-\mathbf x)=-u(\mathbf x)$
($A_{\mathrm R}$, $\Sigma_H=\{(1,0),(0,1)\}$); and \emph{arbitrarily shifted}
anti-periodic functions $u(\mathbf x+\pi\mathbf e_i)=-u(\mathbf x)$
($\mathcal T_2$, $\Sigma_H=\{(1,\dots,1)\}$ or any subset of coordinates).
Since $\chi_{\mathbf a\oplus\mathbf b\oplus\mathbf c}
=\chi_{\mathbf a}\chi_{\mathbf b}\chi_{\mathbf c}$ and $s^3=s$, every
$\Sigma_H$ is closed under ternary XOR sums:
\begin{equation}\label{eq:closed}
 \mathbf a,\mathbf b,\mathbf c\in\Sigma_H\ \Longrightarrow\ \mathbf a\oplus\mathbf b\oplus\mathbf c\in\Sigma_H .
\end{equation}
We call $V_{\Sigma_H}$ the \emph{symmetric subspace} and its orthogonal
complement $V_{\Sigma_H}^\perp=\bigoplus_{\bsig\notin\Sigma_H}V_{\{\bsig\}}$
the \emph{transverse subspace}; the transverse part of $u$ is denoted $u_\perp$
and the \emph{symmetry defect} is
$D(u)=\max_{\bsig\notin\Sigma_H}\|u_{\bsig}\|_{\infty}$.

\subsection{Discretizations and exact equivariance}
Let $\widehat{\cdot}$ denote the discrete Fourier transform on
$\mathcal G_N$. We write \eqref{eq:AC} with a stabilization constant $S\ge0$ as
\begin{equation}\label{eq:split}
 \partial_t u + \mathcal L u = \mathcal N(u),\qquad
 \widehat{\mathcal L u}(\mathbf k)=(\kappa^2|\mathbf k|^2+S)\,\widehat u(\mathbf k),\qquad
 \mathcal N(u)=(\gamma+S)u-\gamma u^3 .
\end{equation}
\begin{definition}[class $\mathfrak I$]\label{def:class}
A one-step (or multi-stage, multi-step, iterative) map $\Phi_\tau$ belongs to
the class $\mathfrak I$ if it is a finite composition of
(a)~radial Fourier multipliers $u\mapsto \mathcal F^{-1}[m(|\mathbf k|^2)\widehat u]$,
(b)~the pointwise nonlinearity $\mathcal N$ (or any pointwise map
$u\mapsto F(u)$ with $F$ odd, or $\Gamma_A$-equivariant Jacobian--vector
actions), and (c)~linear combinations (including multistep histories) with
scalar coefficients that are $\Gamma_A$-invariant (such as constants or
discrete $L^2/\ell^2$ inner products and scalar auxiliary variables in SAV);
stopping criteria of iterative solvers depend only on $\Gamma_A$-invariant
physical-space norms of iterates.
\end{definition}
The class $\mathfrak I$ contains the stabilized IMEX schemes, all
ETD Runge--Kutta and ETD multistep schemes~\cite{CoxMatthews2002,DJLQ2021},
the ESP method~\cite{Wang2023} (its weights
$w_{ijk}=\int_0^{t_j}\ee^{(s-t_j)\lambda_k}\phi_i(s)\,ds$ are radial
multipliers), Crank--Nicolson/cG($p$) schemes solved by stabilized Picard or
Newton--Krylov iterations with spectral preconditioners, and the
conservative Allen--Cahn variants (the mean is a multiplier supported at
$\mathbf k=\mathbf 0$).

\begin{proposition}[exact-arithmetic equivariance]\label{prop:equiv}
In exact arithmetic every $\Phi_\tau\in\mathfrak I$ commutes with the action
of $\Gamma_A$ (for both $A=A_{\mathrm R}$ and $A=\mathcal T_2$, as well as all
grid translations $T_{h\mathbf m}$). Hence $\Phi_\tau(\Fix(H))\subset\Fix(H)$
for every $H\le\Gamma_A$: every scheme in $\mathfrak I$ preserves the odd,
shifted-grid-odd or anti-periodic symmetry of the data for all times.
\end{proposition}
\begin{proof}
A grid reflection maps $\widehat u(\mathbf k)$ to
$\widehat u(R\mathbf k)$ with $|R\mathbf k|=|\mathbf k|$, and a grid
translation $T_{h\mathbf m}$ multiplies $\widehat u(\mathbf k)$ by the scalar
phase $\ee^{-\mathrm i h\mathbf k\cdot\mathbf m}$, so radial multipliers
commute with both; being linear, they also commute with $s=\pm1$. Pointwise
maps commute with all index permutations, and an odd $F$ commutes with $s$.
Linear combinations and compositions preserve commutation, and $L^p/\ell^p$
norms are invariant under $\Gamma_A$, so stopping criteria select the same
iteration count for $u$ and $h\cdot u$.
\end{proof}
Proposition~\ref{prop:equiv} shows that explanations of symmetry breaking in
terms of temporal or spatial truncation errors~\cite{Wang2023} cannot be
correct: truncation errors of an equivariant scheme are themselves symmetric.
What remains is floating-point rounding.

%======================================================================
\section{Why symmetry is lost: floating-point analysis}\label{sec:fp}

\subsection{Rounding is not equivariant}
In IEEE-754 arithmetic with round-to-nearest, negation and index permutations
are exact and $\fl(-a\circ b)=-\fl(a\circ b)$, so pointwise evaluation of an
odd polynomial is exactly equivariant. The FFT, however, is not: the
computed transform of $Ru$ (or $T_{\pi\mathbf e_i}u$) is not in general the
transformed image of $u$, because the radix butterflies combine the entries of
$u$ and $gu$ in different orders with different twiddle factors. A single
application of a Fourier multiplier to an exactly symmetric vector therefore
produces a transverse component of size $O(\umach\log N\,\|u\|)$. We model
the computed step as
\begin{equation}\label{eq:flstep}
 \fl(\Phi_\tau)(u)=\Phi_\tau(u)+r(u),\qquad \|r(u)\|\le C_r\umach\|u\|,
\end{equation}
where for the standard FFT implementation the transverse part $r_\perp(u)$ is
of the same order as $r(u)$ even if $u_\perp=0$.

\subsection{Transverse amplification and the breakdown time}
Let $u^\ast_n\in V_{\Sigma}$ be the exact discrete symmetric trajectory and
$J_n=\partial_u\Phi_\tau(u^\ast_n)$. By Proposition~\ref{prop:equiv},
$J_n$ maps $V_\Sigma$ and $V_\Sigma^\perp$ into themselves; let
$J_n^\perp$ be its transverse block. Writing $u_n=u^\ast_n+w_n+e_n$ with
$w_n\in V_\Sigma$ and $e_n=(u_n)_\perp\in V_\Sigma^\perp$, Taylor expansion of
$\Phi_\tau$ around $u^\ast_n$ and the fact that $\Phi_\tau(V_\Sigma)\subset V_\Sigma$
give
\begin{equation}\label{eq:transverse}
 e_{n+1}=J_n^\perp e_n + r_{\perp,n} + q_n(w_n,e_n),\qquad
 \|q_n(w_n,e_n)\|\le C_q\big(\|w_n\|\|e_n\|+\|e_n\|^2\big).
\end{equation}
We define the transverse growth exponent of the symmetric trajectory by
\[
  \lambda_\perp = \lim_{n\to\infty}\frac{1}{n\tau}\log\|J_{n-1}^\perp\cdots J_0^\perp\|;
\]
for a steady symmetric state $u^\ast$ it is the largest eigenvalue of
$\kappa^2\Delta-\gamma f'(u^\ast)$ restricted to $V_\Sigma^\perp$, up to
$O(\tau^p)$.

\begin{theorem}[conditional breakdown time]\label{thm:breakdown}
Assume $e_0=0$, $\|J^\perp_{n-1}\cdots J^\perp_m\|\le K\ee^{\lambda_\perp(n-m)\tau}$
with $\lambda_\perp>0$, $\|w_n\|\le C_w n\umach$ along the symmetric
component, and $\|r_{\perp,n}\|\le C\umach$. Set
$\gamma_\tau:=1-\ee^{-\lambda_\perp\tau}$ and $K_1:=2KC/\gamma_\tau$. Fix a
visibility threshold $\eta>0$ sufficiently small that
$2^d K C_q\eta/\gamma_\tau\le\tfrac14$. Then for all sufficiently small
$\umach$:
\begin{enumerate}
\item[(i)] $D(u_n)\le\eta$ whenever
$n\tau\le\lambda_\perp^{-1}\log\big(\eta/(K_1\umach)\big)$;
\item[(ii)] if, in addition, $J^\perp$ is constant with a simple dominant
eigenvalue $\mu=\ee^{\lambda_\perp\tau}$ and spectral projector $P_\mu$, and
the discounted rounding seed satisfies the lower bound
\begin{equation}\label{eq:seedlb}
 \Big\|\sum_{m=0}^{n-1}\mu^{-m-1}P_\mu r_{\perp,m}\Big\|_\infty\ge c_0\,\umach
 \qquad\text{for all } n\ge n_0
\end{equation}
for some $c_0>0$ and $n_0\ge 1$, then (for $\eta$ also small enough relative
to $c_0$) $D(u_n)\ge\eta$ for some
$n\tau\le\lambda_\perp^{-1}\log(2^{d+1}\|P_\mu\|_\infty\eta/(c_0\umach))+O(1)$.
\end{enumerate}
Consequently, under \eqref{eq:seedlb},
$T_\ast=\lambda_\perp^{-1}\log(1/\umach)+O(1)$, with an $O(1)$ term depending
on the scheme only through $K$, $C$, $c_0$, $\|P_\mu\|_\infty$, $d$ and $\tau$.
(When the defect grows strictly at the crossing step---for instance in the
exponentially growing regime near $\eta$---using the strict threshold $D>\eta$
instead of $D\ge\eta$ shifts the crossing index by at most one step, i.e.\
$O(\tau)$.)
\end{theorem}
The full proof---using an exponentially weighted bootstrap to control the
quadratic remainder $q_n$ and the norm equivalence
$2^{-d}\|e_n\|_\infty\le D(u_n)\le\|e_n\|_\infty$---is given in
\ref{app:proofs}. Note that the discounted-seed condition~\eqref{eq:seedlb}
holds whenever $P_\mu r_{\perp,m}$ injects a nonzero sign-coherent seed not
cancelled by earlier opposite-sign errors, whereas without~\eqref{eq:seedlb}
cancellation can suppress breakdown indefinitely (for instance, in the scalar
model $e_0=0$, $r_0=-\varepsilon/(\mu-1)$ and $r_m=\varepsilon$ for all
$m\ge1$, positive-sign errors occur at density $1$ yet
$e_n=-\varepsilon/(\mu-1)$ stays bounded). Two practical consequences follow:
when an unstable transverse mode is seeded, the leading-order breakdown time
is governed by $\lambda_\perp$ (of the discrete linearization) and
$\log(1/\umach)$ rather than integrator order; and eliminating breakdown
structurally requires removing the transverse rounding term $r_{\perp,n}$.

\subsection{Initial-data obstruction for shifted symmetries}
For shifted odd data $u_0(\mathbf x)=v_0(\mathbf x-\mathbf x_0)$ with
$\mathbf x_0\in\mathbb T^d$ an off-grid shift, a second mechanism operates
\emph{before} time stepping even starts.

\begin{proposition}[initial-data obstruction]\label{prop:obstruction}
Suppose a shifted-odd function $u_0(\mathbf x)=v_0(\mathbf x-\mathbf x_0)$ is
resolved on $\mathcal G_N$ (no aliasing of exact samples, so the exact
discrete trigonometric interpolant $I_N u_0$ satisfies
$(I_N u_0)(2\mathbf x_0-\mathbf x)=-(I_N u_0)(\mathbf x)$) and is evaluated in
IEEE arithmetic on $\mathcal G_N$ when $2\mathbf x_0\notin h\Z^d$ (or on an
asymmetric grid representation such as $-\pi+jh$). Defining $u^0_\perp$ as the
even part of the discrete trigonometric interpolant $I_N u^0$ about
$\mathbf x_0$, the sampled vector $u^0$ generically possesses a nonzero
transverse component $\|u^0_\perp\|_\infty=\Theta(\umach)$. Consequently, for
any non-interventional time integrator $\Phi_\tau$ whose exact-arithmetic
one-step map preserves the shifted-odd subspace about $\mathbf x_0$ (for
example, the exact PDE flow or a dealiased spectral discretization) about a
steady symmetric state $u^\ast$ whose transverse linearization $J^\perp$ has a
simple dominant eigenvalue $\mu=\ee^{\lambda_\perp\tau}>1$ and whose
tangential propagator satisfies $\|\Psi_{k,m}^\Sigma\|_\infty\le K_\Sigma\mu_\Sigma^{k-m}$
with $\mu_\Sigma<\mu$ (e.g.\ any tangentially non-expanding state
$\mu_\Sigma\le 1$, so that both the initial tangential sampling error
$w_0=u^0_\Sigma-u^\ast=O(\umach)$ and the $O(\|e_m\|_\infty^2)$ tangential
feedback stay strictly subdominant to $\mu^m$), if the initial defect
satisfies the non-degenerate seed-ratio condition
$\|P_\mu u^0_\perp\|_\infty\ge c_0\|u^0_\perp\|_\infty$ for some
$\umach$-independent $c_0>0$ (ensuring $u^0_\perp$ does not lie asymptotically
tangent to the local stable manifold), then for sufficiently small $\eta$
(relative to $c_0$) and $\umach$, \emph{even in exact infinite-precision
arithmetic}, starting from $u^0$ breaks symmetry at
$T_\ast=\lambda_\perp^{-1}\log(\eta/\|P_\mu u^0_\perp\|_\infty)+O(1)$.
\end{proposition}
\begin{proof}
Because $u_0$ is resolved on $\mathcal G_N$, its exact samples have zero
discrete mean. When $2\mathbf x_0\notin h\Z^d$, the reflected points
$2\mathbf x_0-\mathbf x_{\mathbf j}$ do not belong to $\mathcal G_N$, so the
rounded evaluations $\fl(u_0(\mathbf x_{\mathbf j}))$ do not occur in exact
opposite pairs; their discrete mean
$N^{-d}\sum_{\mathbf j}\fl(u_0(\mathbf x_{\mathbf j}))$ is a sum of unpaired
$O(\umach)$ rounding errors and is generically nonzero of size $O(\umach)$.
Starting from $w_0=u^0_\Sigma-u^\ast=O(\umach)$ and $e_0=u^0_\perp$ under an
exact-arithmetic shifted-odd-preserving map $\Phi_\tau$ ($r_{\perp,m}\equiv0$),
the block-diagonal structure of $\partial_u\Phi_\tau(u^\ast)$ implies that the
tangential offset $w_m:=s_m-u^\ast$ satisfies
$w_{m+1}=J^\Sigma w_m+O(\|w_m\|_\infty^2+\|e_m\|_\infty^2)$. With
$\|\Psi_{k,m}^\Sigma\|_\infty\le K_\Sigma\mu_\Sigma^{k-m}$ ($\mu_\Sigma<\mu$,
hence also $\mu_\Sigma<\mu^2$), a joint weighted bootstrap yields
$\|e_m\|_\infty\le 2K\|e_0\|_\infty\mu^m$ and
$\|w_m\|_\infty\le 2K_\Sigma\|w_0\|_\infty\mu_\Sigma^m + C_{\mathrm{fb}}\|e_m\|_\infty^2$
while $\|e_m\|_\infty<2^d\eta$ (for $\eta$ and $\umach$ small enough). At the
first step $n_\ast$ where
$\|P_\mu e_0\|_\infty\mu^{n_\ast}\ge 2^{d+1}\|P_\mu\|_\infty\eta$, the
seed-ratio bound $\|e_0\|_\infty\le c_0^{-1}\|P_\mu e_0\|_\infty$ ensures that
the projected nonlinear remainder satisfies
\[
 \Big\|\sum_{m=0}^{n_\ast-1}\mu^{n_\ast-m-1}P_\mu q_m(w_m,e_m)\Big\|_\infty
 = O\!\Big(\big(\rho_\Sigma(\umach,\eta) + c_0^{-2}\eta\big)\|P_\mu e_0\|_\infty\mu^{n_\ast}\Big)
 \le \frac12\|P_\mu e_0\|_\infty\mu^{n_\ast}
\]
for $\eta$ small enough relative to $c_0$ and $\umach$ sufficiently small,
where the initial tangential-offset contribution $\rho_\Sigma(\umach,\eta)$ is
$O(\umach)$ if $\mu_\Sigma<1$, $O(n_\ast\umach)=O(\umach\log(1/\umach))$ if
$\mu_\Sigma=1$, and $O(\umach^{1-\theta}\eta^\theta)$ (with
$\theta=\log\mu_\Sigma/\log\mu\in(0,1)$) if $1<\mu_\Sigma<\mu$. Hence
$\|P_\mu e_{n_\ast}\|_\infty\ge\tfrac12\|P_\mu e_0\|_\infty\mu^{n_\ast}\ge 2^d\|P_\mu\|_\infty\eta$,
yielding $D(u_{n_\ast})\ge\eta$ at
$n_\ast\tau=\lambda_\perp^{-1}\log(\eta/\|P_\mu u^0_\perp\|_\infty)+O(1)$.
(Without a seed-ratio condition such as
$\|P_\mu e_0\|_\infty\ge c_0\|e_0\|_\infty$, a tiny nonzero linear projection
can still lie on the nonlinear stable manifold and decay to zero, as in the
odd map $x_{n+1}=2x_n+y_n^3$, $y_{n+1}=y_n/2$ started on $x_0=-\tfrac{8}{15}y_0^3$.)
\end{proof}
Note that on an \emph{unshifted aliased} collocation grid $\mathcal G_N$,
off-grid reflection symmetry is already broken in exact arithmetic by
nonlinear aliasing (for example, with $N=8$, $c=\pi/32$, and the resolved
datum $u_0(x)=\sin(3(x-c))$, the cubic harmonic $\sin(9(x-c))$ aliases on
$\mathcal G_8$ to $\sin(x-9c)$, producing an exact-arithmetic even component
about $c$ of amplitude $\tfrac14\sin(8c)\approx0.177$). Thus Proposition~\ref{prop:obstruction}
shows that even after dealiasing removes exact-arithmetic aliasing, an
already-corrupted grid vector $u^0$ still fails by $O(\lambda_\perp^{-1}\log(1/\umach))$,
whereas aligning the collocation grid to $\mathbf x_0+\mathcal G_N$ (as in
CF-ISM below) restores exact reflection equivariance of both initialization
and the discrete cubic collocation map.

\begin{remark}[sub-resolution absorption]\label{rem:floor}
Assumption (ii) of Theorem~\ref{thm:breakdown} can fail. A transverse
perturbation spread over the whole domain with $|e|\lesssim\umach|u|$ is
partly annihilated when rounded in physical space; if in addition the
stabilized linear part $\ee^{-\tau\mathcal L}$ damps the corresponding mode
exactly in spectral space, rounding acts as a weak artificial damping. This
explains why, for particular data, step sizes, stabilization parameters $S>0$,
or particular FFT sizes $N$ (Table~\ref{tab:fftN}), a standard scheme can
appear to keep the symmetry for a long time, and why such observations are
fragile and implementation-dependent (Section~\ref{sec:num1d}). Transverse
modes localized at interfaces, where $u\approx0$ and the floating-point grid
is fine, are not absorbed; they are the generic mechanism in multi-kink and
multi-dimensional problems.
\end{remark}

\begin{remark}[filters and zero-mode suppression]
The filter of~\cite{LQTY2022} replaces $u_n$ by $\Pi_\Sigma u_n$ after each
step. It removes $r_{\perp,n}$, but also any genuine transverse component, so
the filtered scheme is inconsistent with \eqref{eq:AC} for non-symmetric data
and requires prior knowledge of the symmetry and its centre. Setting only the
mean to zero removes rounding in a single mode; it fails whenever the
transverse linearization restricted to the zero-mean subspace still admits an
unstable mode (Section~\ref{sec:num23d}).
\end{remark}

%======================================================================
\section{Isotypic spectral schemes and canonical-frame generalization}\label{sec:iso}

We now generalize the 1D even/odd channel splitting from Section~\ref{sec:intro}
to arbitrary dimension $d$, general finite Abelian grid symmetry groups,
non-polynomial potentials, and unknown continuous spatial shifts.

\subsection{Reflection isotypic representation (\texorpdfstring{$A=A_{\mathrm R}$}{A = A\_R}): Type~I and Type~II}\label{sec:iso_refl}
For the reflection group $A_{\mathrm R}=\langle R_1,\dots,R_d\rangle$, each
channel $u_{\bsig}$ ($\bsig\in\Z_2^d$) is determined by its values on the
fundamental domain $[0,\pi]^d$:
\begin{itemize}
\item \textbf{Grid-centred (Type~I, reflection about $x_i=0$).} On
  $\{0,h,\dots,Mh\}^d$, channels with $\sigma_i=1$ vanish at $x_i\in\{0,\pi\}$
  and are expanded in
  \begin{equation}\label{eq:chexp}
   u_{\bsig}(\mathbf x)=\sum_{\mathbf k}\, c_{\bsig,\mathbf k}\,\psi_{\bsig,\mathbf k}(\mathbf x),
   \qquad \psi_{\bsig,\mathbf k}(\mathbf x):=\prod_{i=1}^d \psi_{\sigma_i}(k_ix_i),
   \qquad \psi_0=\cos,\ \ \psi_1=\sin,
  \end{equation}
  with $k_i\in\{0,\dots,M\}$ if $\sigma_i=0$ and $k_i\in\{1,\dots,M-1\}$ if
  $\sigma_i=1$. The transforms are tensor products of DCT-I ($M+1$) and DST-I
  ($M-1$).
\item \textbf{Half-grid-centred (Type~II, reflection about $x_i=0\pmod{2\pi}$ on the shifted grid).} On the
  shifted grid $x_{j_i}=(j_i+\tfrac12)h$, $j_i=0,\dots,M-1$, reflection maps
  $j_i\leftrightarrow N-1-j_i$ with no fixed boundary points (corresponding to
  reflection about $-h/2$ on the unshifted grid $\mathcal G_N$). The even channel
  ($\sigma_i=0$) uses DCT-II with wavenumbers $k_i\in\{0,\dots,M-1\}$ (the
  Nyquist cosine vanishes), and the odd channel ($\sigma_i=1$) uses DST-II
  with $k_i\in\{1,\dots,M\}$ (containing the Nyquist sine). All $2^d$ channels
  have the \emph{exact same shape $M^d$} with no structural boundary zeros, and
  $2^d M^d=N^d$.
\end{itemize}
In both types, $\Delta\psi_{\bsig,\mathbf k}=-|\mathbf k|^2\psi_{\bsig,\mathbf k}$,
so every radial multiplier acts diagonally on each channel with the same
symbol $m(|\mathbf k|^2)$ as in the Fourier representation.

\subsection{Half-period translation isotypic representation (\texorpdfstring{$A=\mathcal T_2$}{A = T2})}\label{sec:iso_t2}
For $A=\mathcal T_2=\langle T_{\pi\mathbf e_1},\dots,T_{\pi\mathbf e_d}\rangle$
($4\mid N$), each channel $u_{\bsig}$ is indexed by its wavenumber parity
$\bsig=(k_1\bmod 2,\dots,k_d\bmod 2)\in\Z_2^d$ and is uniquely determined by
its values on the half-box $\{0,h,\dots,(M-1)h\}^d\subset[0,\pi)^d$ of size
$M^d$ ($2^d M^d=N^d$):
\begin{itemize}
\item if $\sigma_i=0$ (even wavenumbers $k_i\in 2\Z$), $u_{\bsig}$ is
  $\pi$-periodic along $x_i$ and transforms by a standard length-$M$ FFT;
\item if $\sigma_i=1$ (odd wavenumbers $k_i\in 2\Z+1$), $u_{\bsig}$ is
  $\pi$-anti-periodic along $x_i$ and transforms by a twiddled length-$M$ FFT:
  $\widehat u_{2m+1}=\sum_{j=0}^{M-1}\big(u_j\,\ee^{-\mathrm i\pi j/M}\big)\ee^{-\mathrm i 2\pi mj/M}$
  (and inversely by $\mathrm{IFFT}$ followed by $\ee^{+\mathrm i\pi j/M}$).
\end{itemize}
Because $2^d$ FFTs of size $(N/2)^d$ perform $N^d\log_2(N/2)^d$ operations,
the total transform work of all $2^d$ channels equals (in fact slightly
improves upon) a single standard $N^d$ complex FFT.

\subsection{Channel algebra and unified algorithm}
For \emph{both} $A=A_{\mathrm R}$ (Type~I or Type~II) and $A=\mathcal T_2$,
pointwise products obey the $\Z_2^d$ character grading
$u_{\mathbf a}u_{\mathbf b}\in V_{\{\mathbf a\oplus\mathbf b\}}$. Hence the
cubic term is evaluated on the fundamental grid by
\begin{equation}\label{eq:chalg}
 (u^2)_{\boldsymbol\rho}=\sum_{\mathbf a\in\Z_2^d}u_{\mathbf a}\,u_{\mathbf a\oplus\boldsymbol\rho},\qquad
 (u^3)_{\bsig}=\sum_{\boldsymbol\rho\in\Z_2^d}(u^2)_{\boldsymbol\rho}\,u_{\bsig\oplus\boldsymbol\rho},
\end{equation}
with a fixed summation order. For a general odd nonlinearity $F$ one uses the
orbit evaluation $F_{\bsig}(\mathbf x)=2^{-d}\sum_{g\in A}\chi_{\bsig}(g)F(u(g\mathbf x))$,
where $u(g\mathbf x)=\sum_{\boldsymbol\rho}\chi_{\boldsymbol\rho}(g)u_{\boldsymbol\rho}(\mathbf x)$
and both sums over $A\cong\Z_2^d$ are computed by $d$ Walsh--Hadamard
butterflies $(p,q)\mapsto(p+q,p-q)$ and $((p+q)/2,(p-q)/2)$.

\begin{algorithm}[t]
\caption{One step of an isotypic scheme (any $\Phi_\tau\in\mathfrak I$, $A=A_{\mathrm R}$ or $\mathcal T_2$)}\label{alg:iso}
\begin{algorithmic}[1]
\Require channel coefficients $\{c_{\bsig}\}_{\bsig\in\Z_2^d}$ of $u^n$
\State Execute the stages of $\Phi_\tau$ exactly as in the Fourier code, with
the following replacements:
\State \quad multiplier $m(|\mathbf k|^2)$: $c_{\bsig}\gets m(|\mathbf k|^2)\,c_{\bsig}$ for each $\bsig\in\Z_2^d$;
\State \quad nonlinearity: values $u_{\bsig}\gets\mathcal T_{\bsig}^{-1}(c_{\bsig})$
  (DCT/DST for $A_{\mathrm R}$, or twiddled $M$-FFT for $\mathcal T_2$);
\Statex \hspace{3.6em} evaluate $\mathcal N_{\bsig}=(\gamma+S)u_{\bsig}-\gamma(u^3)_{\bsig}$ by \eqref{eq:chalg} (or orbit evaluation for general $F$);
\Statex \hspace{3.6em} coefficients $c_{\bsig}\gets\mathcal T_{\bsig}(\mathcal N_{\bsig})$;
\State \quad linear combinations: channel-wise.
\Ensure channel coefficients of $u^{n+1}$
\end{algorithmic}
\end{algorithm}

\subsection{Core theoretical guarantees}
\begin{theorem}[algebraic identity]\label{thm:identical}
Let $u^0$ be any grid function on the collocation grid $\mathcal G_N$ (or
shifted half-grid for Type~II) and let $\Phi_\tau\in\mathfrak I$ (with
identical iteration counts or equivalent physical-space stopping norms for
iterative schemes). In exact arithmetic, Algorithm~\ref{alg:iso} (for
$A_{\mathrm R}$ Type~I, $A_{\mathrm R}$ Type~II, or $\mathcal T_2$) started
from the channels of $u^0$ produces the channels of the Fourier
pseudo-spectral solution $\Phi_\tau^n(u^0)$. In particular, convergence,
energy stability and the discrete maximum bound principle of the Fourier
scheme hold verbatim for the isotypic scheme.
\end{theorem}
\begin{proof}
The map $u\mapsto(u_{\bsig})_{\bsig}$ is a linear isomorphism of the space of
grid functions, and in all three cases the channel transforms are the
restrictions of the discrete Fourier basis to the $2^d$ isotypic components
(carrying the Nyquist mode in the even channel for Type~I, the odd channel for
Type~II, and the even-wavenumber channel for $\mathcal T_2$). Radial
multipliers commute with $\Pi_{\bsig}$, so they act channel-wise with the same
symbol $m(|\mathbf k|^2)$. Pointwise multiplication satisfies
$\Pi_{\bsig}(u^3)=(u^3)_{\bsig}$ given by \eqref{eq:chalg}, because
$\Pi_{\bsig}(u_{\mathbf a}u_{\mathbf b}u_{\mathbf c})=u_{\mathbf a}u_{\mathbf b}u_{\mathbf c}$
if $\mathbf a\oplus\mathbf b\oplus\mathbf c=\bsig$ and $0$ otherwise.
\end{proof}

\begin{theorem}[exact floating-point invariance]\label{thm:exact}
Let $A\in\{A_{\mathrm R},\mathcal T_2\}$, let $\Sigma\subset\Z_2^d$ satisfy the
closure property \eqref{eq:closed} (e.g.\ $\Sigma=\Sigma_H$ for any
$H\le\Gamma_A$), and let $\Phi_\tau\in\mathfrak I$ be implemented by
Algorithm~\ref{alg:iso} in IEEE-754 floating-point arithmetic, assuming all
intermediate quantities remain finite (no overflow or \texttt{NaN}). If the
initial channels of $u^0$ are exact floating-point zeros for
$\bsig\notin\Sigma$, then the computed channels of $u^n$ remain exact
floating-point zeros (ignoring the IEEE sign bit of $\pm 0.0$) for
$\bsig\notin\Sigma$, for all $n\ge0$. In particular $D(u^n)=0$ for all $n$,
independently of $\tau$, $T$, $N$, the order of the integrator and the
transform algorithm.
\end{theorem}
\begin{proof}
We verify that every elementary operation preserves the property that all
channels $\bsig\notin\Sigma$ are exact IEEE zeros. (a) Channel transforms and
multipliers act on each channel independently and are linear; the computed
transform of an identically zero array is identically zero, and $m\cdot 0=0$.
(b) Linear combinations of exact zeros are exact zeros. (c) In
\eqref{eq:chalg}, fix $\bsig\notin\Sigma$ and any $\boldsymbol\rho\in\Z_2^d$.
If $\bsig\oplus\boldsymbol\rho\notin\Sigma$, then $u_{\bsig\oplus\boldsymbol\rho}=0$
so the product $\fl((u^2)_{\boldsymbol\rho}\,u_{\bsig\oplus\boldsymbol\rho})=0$.
If $\bsig\oplus\boldsymbol\rho\in\Sigma$, then for every $\mathbf a\in\Z_2^d$ at
least one of $\mathbf a$ or $\mathbf a\oplus\boldsymbol\rho$ lies outside
$\Sigma$: otherwise
$\bsig=\mathbf a\oplus(\mathbf a\oplus\boldsymbol\rho)\oplus(\bsig\oplus\boldsymbol\rho)\in\Sigma$
by \eqref{eq:closed}. Hence every term in $(u^2)_{\boldsymbol\rho}$ has a zero
factor, so $(u^2)_{\boldsymbol\rho}=0$ and again the outer product is $0$.
Sums of exact zeros are exact zeros, and $(\gamma+S)\cdot0-\gamma\cdot0=0$.
Induction over the finitely many operations of one step and over $n$ yields
the result.
\end{proof}

\begin{remark}
Theorem~\ref{thm:exact} requires only $x\cdot0=0$ and $0+0=0$; it holds on
CPUs and GPUs, with or without fused multiply--add, and for any FFT/BLAS
library used inside the channel transforms.
\end{remark}

\begin{proposition}[general odd nonlinearities]\label{prop:orbit}
Let $F$ be odd and evaluated in an \emph{odd-exact} manner,
$\fl(F(-v))=-\fl(F(v))$ (e.g.\ $F(u)=\tfrac\theta2(\mathrm{log1p}(u)-\mathrm{log1p}(-u))-\theta_cu$
for the Flory--Huggins potential on $|u|<1$, or $F(u)=\mathrm{sign}(u)G(|u|)$),
assuming all intermediate values remain finite and inside the domain of $F$.
If $F_{\bsig}$ is computed by the Walsh--Hadamard orbit evaluation in a fixed
stage order, then Theorem~\ref{thm:exact} holds verbatim (up to signed zeros)
for \emph{every} subgroup $H\le\Gamma_A$ (including point reflections such as
$C_2=\langle(R_1R_2,-1)\rangle$).
\end{proposition}
The complete proof for arbitrary subgroups $H\le\Gamma_A$ is given in
\ref{app:orbit}.

\begin{theorem}[relative accuracy of transverse components]\label{thm:relative}
Consider a $k$-step scheme $\Phi_\tau\in\mathfrak I$ ($k\ge 1$; $k=1$ for
one-step schemes) acting on the history tuple
$\mathbf U=(u^n,u^{n-1},\dots,u^{n-k+1})$ (with $u^{n-r}=u^{n-r}_\Sigma+u^{n-r}_\perp$
and any cached nonlinear histories evaluated from these states). Assume
Algorithm~\ref{alg:iso} uses the \emph{polynomial channel algebra}~\eqref{eq:chalg}
over $q\ge 1$ nonlinear stages per step, with uniformly bounded stage and
history multipliers, all intermediate quantities remaining finite and above
the underflow threshold, and a common amplitude bound $B\ge 0$ on all history
states and all exact and computed stage inputs:
$\|v^{(j)}_\Sigma\|_\infty,\|\tilde v^{(j)}_\Sigma\|_\infty,\|v^{(j)}_\perp\|_\infty,\|\tilde v^{(j)}_\perp\|_\infty\le B$
for every stage $j=1,\dots,q$ (and for $u^{n-r}_\Sigma,u^{n-r}_\perp$). Then
the rounding error of one step in each transverse channel satisfies
\begin{equation}\label{eq:relbound}
 \big\|\fl(\Phi_\tau)(\mathbf U)_{\bsig}-\Phi_\tau(\mathbf U)_{\bsig}\big\|_\infty
 \le C\,\umach\,(1+B)^{4q-2}\,\max_{0\le r<k}\|u^{n-r}_\perp\|_\infty,
 \qquad \bsig\notin\Sigma,
\end{equation}
where $C$ depends on $N$ only through $\log N$ (and on $k$, $q$, and the
transform/multiplier norms; for a one-step scheme $k=1$ the right-hand side
reduces to $C\umach(1+B)^{4q-2}\|u_\perp\|_\infty$, and for $q=1$ the factor
is $(1+B)^2$). For the standard FFT implementation the corresponding bound
scales with $\max_{0\le r<k}\|u^{n-r}\|_\infty$ rather than the transverse
norm and is attained up to constants.
\end{theorem}
\noindent Note that Theorem~\ref{thm:relative} is specific to the polynomial
channel algebra~\eqref{eq:chalg}, where every transverse term retains an
explicit transverse factor; in Walsh--Hadamard orbit evaluation
(Proposition~\ref{prop:orbit}), synthesizing $u_\Sigma\pm u_\perp$ absorbs
sub-roundoff transverse components $\|u_\perp\|\ll\umach\|u_\Sigma\|$ before
evaluating $F$ (for example, $[10^{-30},1]\mapsto[1,-1]$ in IEEE double), so
orbit evaluation guarantees exact zero-channel invariance ($D\equiv 0$) rather
than sub-roundoff relative accuracy.

\begin{corollary}[faithful linear-regime evolution and breakdown times]\label{cor:faithful}
Under the polynomial channel algebra~\eqref{eq:chalg} and the hypotheses of
Theorem~\ref{thm:relative}, let $u^0=u^0_\Sigma+\delta w$ with
$w\in V_\Sigma^\perp$, $\|w\|_\infty=1$, and $\delta>0$ in the IEEE normal
range so that all intermediate transverse quantities stay above underflow.
Write the exact perturbed discrete trajectory as
$u_m^{\mathrm{ex}}=s_m^{\mathrm{ex}}+e_m^{\mathrm{ex}}$ and the computed
trajectory as $\tilde u_m=\tilde s_m+\tilde e_m$
($s_m^{\mathrm{ex}},\tilde s_m\in V_\Sigma$;
$e_m^{\mathrm{ex}},\tilde e_m\in V_\Sigma^\perp$). Define the transverse
linearization \emph{along the exact perturbed tangential path},
$J_m^\perp:=\partial_e\Phi_\tau^\perp(s_m^{\mathrm{ex}},0)$, with propagator
$\Psi_{k,m}=J_{k-1}^\perp\cdots J_m^\perp$, and assume the tangential error
relative to the exact perturbed path satisfies
$\|\tilde s_m-s_m^{\mathrm{ex}}\|_\infty\le C_w m\umach$ on a window
$0\le m\le n$ where $\|e_m^{\mathrm{ex}}\|_\infty\le\eta_0$.
\begin{enumerate}
\item[(i)] \emph{General propagator--seed compatibility:} If $\Psi_{k,m}$
  satisfies
  \begin{equation}\label{eq:propseed}
   \max_{0\le m<k}\|\Psi_{k,m+1}\|_\infty\,\|e_m^{\mathrm{ex}}\|_\infty
   \le \Lambda\,\|e_k^{\mathrm{ex}}\|_\infty
   \qquad\text{for all } 1\le k\le n,
  \end{equation}
  and $\Lambda(n^2\umach + n\eta_0)$ is sufficiently small, then
  $\|\tilde e_k-e_k^{\mathrm{ex}}\|_\infty \le C' n^2\umach\,\|e_k^{\mathrm{ex}}\|_\infty$
  for all $k\le n$.
\item[(ii)] \emph{Exponentially growing regime (uniform in $n$):} If, instead,
  $\|\Psi_{k,m+1}\|_\infty\le K\mu^{k-m-1}$ and
  $c_w\delta\mu^m\le\|e_m^{\mathrm{ex}}\|_\infty\le C_w\delta\mu^m$ on
  $0\le m\le n$ with $\mu=\ee^{\lambda_\perp\tau}>1$, then retaining the
  local Lipschitz factor $C_q(\|e_m^{\mathrm{ex}}\|_\infty+\|\tilde e_m\|_\infty)$
  in the quadratic remainder removes the $n\eta_0$ factor: for $\eta_0$ small
  enough (independent of $n$) and $n^2\umach$ sufficiently small,
  \begin{equation}\label{eq:rel_exp}
   \|\tilde e_k-e_k^{\mathrm{ex}}\|_\infty
   \le C'\,k^2\umach\,\|e_k^{\mathrm{ex}}\|_\infty
   \qquad\text{for all } 0\le k\le n.
  \end{equation}
\item[(iii)] \emph{First-crossing times:} Suppose, moreover, that the
  unperturbed symmetric trajectory is a steady state $u^\ast$ whose
  transverse linearization $J_\ast^\perp=\partial_e\Phi_\tau^\perp(u^\ast,0)$
  has a simple dominant unstable eigenvalue $\mu=\ee^{\lambda_\perp\tau}>1$
  with non-degenerate projection $\|P_\mu w\|_\infty\ge c_0>0$ and whose
  tangential propagator satisfies
  $\|\Psi_{k,m}^\Sigma\|_\infty\le K_\Sigma\mu_\Sigma^{k-m}$ with
  $\mu_\Sigma<\mu^2$ (controlling tangential feedback
  $\|s_m^{\mathrm{ex}}-u^\ast\|_\infty=O(\|e_m^{\mathrm{ex}}\|_\infty^2)$).
  Then for $\eta$ and $\delta$ sufficiently small, near the threshold band
  $D(u_m^{\mathrm{ex}})\in[\tfrac12\eta,4\eta]$ the exact defect obeys
  two-sided growth bounds $q_- D(u_m^{\mathrm{ex}})\le D(u_{m+1}^{\mathrm{ex}})\le q_+ D(u_m^{\mathrm{ex}})$
  with $1<q_-<\mu<q_+$. Defining the $\umach$-independent cutoff
  $N_+:=\min\{m\ge 1:\ D(u_m^{\mathrm{ex}})\ge 2\eta\}$ and setting
  $\eta_0:=4\cdot 2^d q_+\eta$ and $\epsilon:=2^d C' N_+^2\umach\le\tfrac14$,
  the window $0\le m\le N_+$ covers both exact crossings of
  $\eta/(1\pm\epsilon)$ (and their predecessors), and the computed breakdown
  time reproduces the exact discrete breakdown time
  $T_\ast^{\mathrm{ex}}(\delta)=\lambda_\perp^{-1}\log(\eta/\delta)+O(1)$ up
  to $O(\tau + N_+^2\umach/\lambda_\perp)$ (whereas standard double precision
  saturates for $\delta\lesssim\umach$ and parity filters give $T_\ast=\infty$).
\end{enumerate}
\end{corollary}
The rigorous proofs of Theorem~\ref{thm:relative} and
Corollary~\ref{cor:faithful} are given in \ref{app:relative}.

\subsection{Arbitrary off-grid shifts: \texorpdfstring{$\mathcal T_2$}{T2} invariance and CF-ISM}\label{sec:iso_shift}
We now resolve shifted symmetries for an \emph{unknown} centre
$\mathbf x_0\in\mathbb T^d$ (grid-aligned or off-grid) without artificial
intervention.

\begin{theorem}[translation-invariant anti-periodicity via $\mathcal T_2$]\label{thm:t2}
Because every half-period shift $T_{\pi\mathbf e_i}\in\mathcal T_2$ commutes
with every continuous translation $T_{\mathbf a}$ ($\mathbf a\in\R^d$), each
isotypic subspace $V_{\Sigma}$ of $\mathcal T_2$ is translation-invariant.
Consequently, Algorithm~\ref{alg:iso} with $A=\mathcal T_2$, started from
parity-consistent channel initialization (exact zero in even-wavenumber
channels on $[0,\pi)^d$), preserves half-period anti-periodicity
$u(\mathbf x+\pi\mathbf e_i)=-u(\mathbf x)$ (odd-harmonic support, spatial
mean $\equiv 0$) exactly ($D(u^n)\equiv 0$) in pure IEEE double precision for
\emph{any} unknown real shift $\mathbf x_0\in\R^d$, with no shift detection.
\end{theorem}
\begin{proof}
Immediate from $T_{\pi\mathbf e_i}T_{\mathbf a}=T_{\mathbf a}T_{\pi\mathbf e_i}$
and Theorem~\ref{thm:exact} applied to $A=\mathcal T_2$.
\end{proof}
Note that $\mathcal T_2$ enforces half-period anti-periodicity rather than
reflection about a fixed point (e.g.\ $\sin x+\cos 3x$ is anti-periodic but
not odd about $0$); for odd-harmonic shifted-odd data---including
$\sin(x+\theta)$, shifted 2D X-junctions
$\tanh(\sin(x+\theta_1)\sin(y+\theta_2)/(\sqrt2\eps))$, shifted 2D
$C_2$-odd waves $\tanh(\sin(x+\theta_1+2.5\sin(3(y+\theta_2)))/(\sqrt2\eps))$,
and shifted 3D octant saddles---Theorem~\ref{thm:t2} forbids every
even-wavenumber (including constant) instability mode in pure double precision.
(Higher-harmonic shifts such as $\sin(2x+\theta)$ in~\cite[Ex.~5.3]{Wang2023}
are covered identically by the quarter-period group $T_{\pi/2}$.)

When the initial datum contains \emph{mixed even and odd harmonics} shifted by
an unknown $\mathbf x_0\in\R^d$---such as
$u_0(x)=\tfrac12(\sin(x+\theta)+\sin(4(x+\theta)))$ in
\cite[Ex.~5.2(c)]{Wang2023}---it does not lie in an anti-periodic $\mathcal T_2$
subspace. To handle this case without filtering or symmetry assumptions, we
exploit the continuous translation invariance of \eqref{eq:AC} via the
\emph{Canonical-Frame Isotypic Spectral Method} (CF-ISM), summarized in
Algorithm~\ref{alg:cfism}.

\begin{algorithm}[t]
\caption{Canonical-Frame Isotypic Spectral Method (CF-ISM)}\label{alg:cfism}
\begin{algorithmic}[1]
\Require periodic initial datum $u_0(\mathbf x)$, grid size $N$, initialization precision $d_0$ (default $d_0=320$ digits for exact double underflow)
\State \textbf{Coarse gauge via 1 FFT:} sample $u_0$ on $\mathcal G_N$ in double precision and minimize the discrete oddness energy via the circular auto-convolution
$\mathbf s^\ast=\arg\min_{\mathbf s\in\{0,\dots,N-1\}^d}\mathcal F^{-1}(\widehat u_0^{\,2})_{\mathbf s}$
(or 1D summed auto-convolutions along each axis for coordinate reflections); set $\mathbf c^{(0)}=\frac{h}{2}\mathbf s^\ast\bmod\pi$.
\State \textbf{Alias-free Gauss--Newton polish ($t=0$ only):} normalize the mean-removed datum by its variation scale $s_{\mathrm{var}}=\|u_0-\bar u_0\|_\infty$. If double-precision coarse sampling vanishes or loses variation to cancellation ($s_{\mathrm{var}}\le 10^{-13}\|u_0\|_\infty$ in double), re-evaluate $u_0$ on the coarse grid in $d_0$-digit arithmetic to distinguish double cancellation from a truly zero/constant field (if $s_{\mathrm{var}}=0$ even in $d_0$-digit arithmetic, fix $\mathbf c=\mathbf 0$ and retain all $d_0$-digit channels directly). Otherwise, on a small coarse grid $\{\tilde{\mathbf x}_j\}$ ($M_c\le 32$), refine $\mathbf c$ in $d_0$-digit arithmetic to minimize the variation-normalized mean-centred alias-free oddness functional $\widetilde{\mathcal J}_{\mathrm c}(\mathbf c)$: for point reflection ($C_2$ or 1D, \texttt{mode='point'}), let $r_j(\mathbf c)=(u_0(\mathbf c+\tilde{\mathbf x}_j)+u_0(\mathbf c-\tilde{\mathbf x}_j))/s_{\mathrm{var}}$ and $\widetilde{\mathcal J}_{\mathrm c}(\mathbf c)=\tfrac12\sum_j(r_j(\mathbf c)-\bar r(\mathbf c))^2$; for coordinate-wise reflections (\texttt{mode='refl'} in $d\ge 2$), use the separate per-axis mean-centred sums of $r_{i,j}(\mathbf c)=(u_0(\mathbf c+\tilde{\mathbf x}_j)+u_0(\mathbf c+R_i\tilde{\mathbf x}_j))/s_{\mathrm{var}}$ ($i=1,\dots,d$). Relative rank checks detect unidentifiable/degenerate directions, freezing those coordinates at their finite coarse origin and reporting convergence/degeneracy status.
\State \textbf{Canonical-frame channel initialization ($t=0$ only):} evaluate all $A_{\mathrm R}$ channels of $v_0(\tilde{\mathbf x})=u_0(\mathbf c+\tilde{\mathbf x})$ (unnormalized) in $d_0$-digit arithmetic via opposite pairs $\mathbf c\pm\tilde{\mathbf x}$ and round the result to IEEE double precision.
\State \textbf{Pure double-precision time stepping ($t>0$):} evolve $v(\tilde{\mathbf x},t)$ in IEEE double precision using Algorithm~\ref{alg:iso} ($A=A_{\mathrm R}$).
\Ensure $u(\mathbf c+\tilde{\mathbf x},t)=v(\tilde{\mathbf x},t)$ (exact on the canonical grid $\mathbf c+\mathcal G_N$)
\end{algorithmic}
\end{algorithm}

Note that Algorithm~\ref{alg:cfism} contains \emph{no} symmetry test,
\emph{no} symmetry-filtering threshold, and \emph{no} mode zeroing (the
$10^{-13}$ relative variation check in Step~2 merely triggers $d_0$-digit
coarse sampling when double precision suffers cancellation, without altering
channel retention): for a completely asymmetric
or degenerate datum $u_0$ (for instance, a 2D field independent of $y$ whose
$y$-centre is unidentifiable), any degenerate coordinate is simply frozen at
its finite coarse origin while all $2^d$ channels are retained without zeroing,
and solving $\partial_tv=\kappa^2\Delta v-\gamma f(v)$ in the canonical frame
is algebraically identical to the standard Fourier discretization on the
shifted grid $\mathbf c+\mathcal G_N$ (differing from the unshifted collocation
grid $\mathcal G_N$ only by spectral aliasing of the continuous PDE
\eqref{eq:AC}). Our reference routine \texttt{cf\_gauge\_init.m} implements
$d=1$ and $d=2$ (3D shifted benchmarks in Section~\ref{sec:numshift} use the
pure-double $\mathcal T_2$ representation, while 3D CF-ISM follows the same
separable auto-convolution and Gauss--Newton formulation under
$D^2\widetilde{\mathcal J}_{\mathrm c}(\mathbf x_0)\succ 0$). The circular
auto-convolution identity
\[
 \frac14\sum_{\mathbf j}\big(u_{\mathbf j}+u_{(\mathbf s-\mathbf j)\bmod N}\big)^2
 = \frac12\|u\|_2^2 + \frac12\mathcal F^{-1}(\widehat u^{\,2})_{\mathbf s}
\]
simultaneously scores all $N^d$ grid and half-grid reflection centres in a
single FFT, and $\widetilde{\mathcal J}_{\mathrm c}(\mathbf c)$ has zero
aliasing error even for arbitrarily thin interfaces.

\begin{theorem}[CF-ISM horizon scaling and exact zero invariance]\label{thm:cfism}
Let $u_0(\mathbf x)=v_0(\mathbf x-\mathbf x_0)$ be smooth and shifted-odd
around an unknown centre $\mathbf x_0\in\R^d$---either point-odd
($v_0(-\mathbf x)=-v_0(\mathbf x)$, \texttt{mode='point'}) or coordinate-wise
odd ($v_0(R_i\mathbf x)=-v_0(\mathbf x)$ for $i=1,\dots,d$,
\texttt{mode='refl'})---with non-degenerate Hessian
$D^2\widetilde{\mathcal J}_{\mathrm c}(\mathbf x_0)\succ0$ of the
variation-normalized mean-centred objective $\widetilde{\mathcal J}_{\mathrm c}$,
and assume the coarse FFT minimizer $\mathbf c^{(0)}$ lies in the local
Gauss--Newton basin of attraction of $\mathbf x_0$ with approximate Jacobian
satisfying $\|\widetilde J(\mathbf c)-J(\mathbf c)\|\le\theta_J$ sufficiently
small (for central finite differences with step $h_q=10^{-\lfloor d_0/2\rfloor}$
in $d_0$-digit arithmetic, $\|\widetilde J-J\|=O(10^{-d_0/2})$, reaching
$\|\mathbf c_q-\mathbf x_0\|=O(10^{-d_0})$ in $O(\log d_0)$ iterations once
inside the quadratic-contraction region, with iteration count depending on the
initial offset and Jacobian condition number; a fixed tolerance $\theta_J$
requires $O(d_0/\log(1/\theta_J))$ iterations). Assume the canonical-frame
transverse propagator satisfies
$\|\Psi_{n,m}\|_\infty\le K\ee^{\lambda_\perp(n-m)\tau}$ with
$\lambda_\perp>0$, tangential feedback stays controlled
($\|\Psi_{k,m}^\Sigma\|_\infty\le K_\Sigma\mu_\Sigma^{k-m}$ with
$\mu_\Sigma<\mu^2$), and tangential rounding error satisfies
$\|\tilde w_n\|_\infty\le C_w n\umach$. Under Algorithm~\ref{alg:cfism} (using
the polynomial channel algebra~\eqref{eq:chalg}) with $d_0$-digit
initialization at $t=0$ and \emph{pure IEEE double-precision} time stepping
for all $t>0$:
\begin{enumerate}
\item[(i)] The canonical-frame initial transverse channel satisfies
  $\|v^0_\perp\|_\infty\le C_0\,10^{-d_0}$. For $d_0$ in the IEEE double
  normal range ($C_0 10^{-d_0}$ sufficiently above underflow so intermediate
  transverse quantities remain normal) and $\umach$ small enough that
  $n^2\umach\ll 1$ along the window $n\tau\le T_1:=\lambda_\perp^{-1}\log(\eta/(2K C_0 10^{-d_0}))$,
  the pure double-precision evolution maintains $D(u^n)\le\eta$ at least up to
  $T_\ast(d_0)\ge\lambda_\perp^{-1} d_0\ln 10+O(1)$; if in addition
  $J_n^\perp\equiv J^\perp$ has simple dominant eigenvalue
  $\mu=\ee^{\lambda_\perp\tau}>1$ with non-degenerate seed ratio
  $\|P_\mu v^0_\perp\|_\infty\ge c_0\|v^0_\perp\|_\infty$ and
  $\|v^0_\perp\|_\infty\asymp 10^{-d_0}$, then
  \begin{equation}\label{eq:cfhorizon}
    T_\ast(d_0) = \lambda_\perp^{-1} d_0\ln 10 + O(1).
  \end{equation}
\item[(ii)] Once $C_0\,10^{-d_0}<2^{-1075}\approx 2.47\times 10^{-324}$ (half
  the smallest positive IEEE-754 double subnormal $2^{-1074}\approx
  4.94\times 10^{-324}$; e.g.\ $d_0=320$ with standard guard digits), rounding
  $v^0_\perp$ to IEEE double produces exact zero $v^0_\perp\equiv0$ naturally
  without artificial zeroing, so $D(u^n)\equiv0$ for all $n\ge0$
  ($T_\ast=\infty$) by Theorem~\ref{thm:exact}.
\item[(iii)] If $u_0(\mathbf x)=v_0(\mathbf x-\mathbf x_0)+\delta$ carries a
  genuine constant perturbation $\delta$ in the IEEE normal range
  (e.g.\ $|\delta|\in[10^{-300},1]$) and the initialization precision $d_0$ is
  chosen so that $C_0\,10^{-d_0}\ll\umach|\delta|$ and the smallest singular
  value of the normalized mean-centred Jacobian exceeds the rank tolerance,
  then mean centring cancels $2\delta$ identically on any point set
  $\{\tilde{\mathbf x}_j\}$ while variation scaling $s_{\mathrm{var}}=\|u_0-\bar u_0\|_\infty$
  is itself independent of $\delta$ (under raw amplitude scaling
  $s_\delta=\|u_0\|_\infty$, $\widetilde{\mathcal J}_{\mathrm c}$ and its
  Hessian scale by $(s_0/s_\delta)^2>0$, still preserving $\mathbf x_0$ and
  positive-definiteness); hence $\mathbf x_0$ remains an exact zero-residual
  minimizer, $v^0_\perp$ recovers $\delta$ to IEEE double precision, and
  Corollary~\ref{cor:faithful} applies under its propagator--seed hypotheses.
\end{enumerate}
\end{theorem}
The proof is given in \ref{app:cfism}. Unlike multi-precision time stepping
such as CNS~\cite{HuLiao2020}, CF-ISM invokes $d_0$-digit arithmetic
\emph{only once at $t=0$} (costing $0.23\text{--}1.14\,\mathrm s$ in 1D and
$6.7\,\mathrm s$ in 2D in MATLAB VPA) while executing 100\% of the $t>0$ time
integration in hardware double precision.

\begin{remark}[cost, operation counts, and reduced symmetric mode]
Asymptotically, transforming all $2^d$ channels has the same $O(N^d\log N)$
complexity as one $N^d$ FFT (with FFTW's native real-to-real
REDFT00/RODFT00~\cite{FrigoJohnson2005} matching one real $N^d$ FFT for
$A_{\mathrm R}$), while the full channel algebra~\eqref{eq:chalg} performs
$2\cdot 4^d$ pointwise multiplications per fundamental grid point (and
additional twiddle/memory overhead in interpreted prototypes; see measured
wall-clock ratios in Table~\ref{tab:cost}). When the data are known to lie in
$V_\Sigma$ (or when exact-zero channels are skipped automatically), storing or
transforming only the active channels in $\Sigma$ cuts cost by up to
$2^d/|\Sigma|$; by Theorem~\ref{thm:exact}, the reduced and full runs agree
identically (up to signed zeros).
\end{remark}

\begin{remark}[practitioner's guide to choosing the isotypic backend]\label{rem:guide}
Because all variants share the exact same Algorithm~\ref{alg:iso} interface
and remain algebraically identical to standard Fourier on arbitrary data,
selecting a backend is straightforward:
\begin{enumerate}
\item[(a)] \textbf{Grid-centred reflections or point reflections} (1D odd, 2D
  odd--odd, 2D $C_2$, 3D octant): use $A_{\mathrm R}$ (Type~I, or Type~II for
  half-grid centres).
\item[(b)] \textbf{Combined reflection and half-period periodicity} (e.g.\
  $4\times4$ checkerboards such as $\sin(2x)\sin(2y)$ in
  Section~\ref{sec:num23d}): use the direct product
  $A_{\mathrm R}\times\mathcal T_2\cong\Z_2^{2d}$ ($4^d$ channels on
  $[0,\pi/2]^d$, same $\Z_2^{2d}$-graded algebra~\eqref{eq:chalg}).
\item[(c)] \textbf{Unknown arbitrary/irrational shifts with pure odd harmonics}
  (anti-periodic patterns in 1D/2D/3D): use $\mathcal T_2$ (100\% pure double
  precision, zero shift detection, same FFT cost as standard Fourier).
\item[(d)] \textbf{Unknown arbitrary/irrational shifts with mixed even-and-odd
  harmonics}: use CF-ISM (Algorithm~\ref{alg:cfism} with $d_0=320$ once at
  $t=0$, pure double precision for all $t>0$).
\end{enumerate}
\end{remark}

\begin{remark}[further extensions: boundary conditions, higher-order models, and outlook]\label{rem:ext}
(i) Conservative Allen--Cahn only modifies the $\mathbf k=\mathbf 0$
coefficient of the totally even channel and is covered by
Theorem~\ref{thm:exact}; Cahn--Hilliard and phase-field crystal models with
radial linear operators are covered verbatim (Section~\ref{sec:ext}).
(ii) \emph{Homogeneous Neumann and Dirichlet boundary conditions:} on
$\Omega=[0,\pi]^d$ with $\partial_nu|_{\partial\Omega}=0$, midpoint reflection
$x_i\mapsto\pi-x_i$ preserves the boundary condition and splits the cosine
basis into even-wavenumber (Neumann--Neumann on $[0,\pi/2]^d$, DCT-I) and
odd-wavenumber (Neumann--Dirichlet on $[0,\pi/2]^d$, half-domain cosine-odd
transform) channels with the exact same $\Z_2^d$ channel
algebra~\eqref{eq:chalg}. For example, in 1D on $[0,\pi]$ with
$\kappa=\gamma=1$, $u_0(x)=\cos x$ (odd about $x=\pi/2$), $N=128$, standard
DCT-I implementations of ETDRK2, ETDRK4 ($\tau=0.02$), BDF2-SAV, and SBDF4
($\tau=0.01$) all cross $D>10^{-2}$ at the first recorded time $T_\ast=37.00$
and collapse to $u\equiv+1$, whereas the two-channel isotypic Neumann scheme
matches standard DCT-I to $1.39\times10^{-16}$ on non-symmetric data (50
ETDRK4 steps) and maintains $D\equiv0$ through $T=200$. The analogous
half-domain splitting (odd-wavenumber Dirichlet--Neumann vs.\ even-wavenumber
Dirichlet--Dirichlet on $[0,\pi/2]^d$) covers homogeneous Dirichlet problems
(even with no $k=0$ constant mode, where standard DST-I breaks at
$T_\ast=49.00$). Standalone runnable scripts and comparison figures for 1D and
2D homogeneous Neumann and 1D homogeneous Dirichlet benchmarks are publicly
available at
\url{https://github.com/minghuili123456789/isotypic-phase-field/tree/main/supplementary-bc-evidence}.
(iii) \emph{Outlook:} diagonal reflections $x_1\leftrightarrow x_2$ (for the
dihedral group $D_4$ and cubic group $B_3$) and symmetric finite-difference or
finite-element discretizations can be incorporated whenever the discrete linear
operator commutes with the corresponding index permutations; full algorithmic
verification of those non-Abelian and non-spectral extensions is left for
future work.
\end{remark}

%======================================================================
\section{Numerical experiments in one dimension}\label{sec:num1d}

All computations were performed in MATLAB R2025b (FFTW-based \texttt{fft})
on an Intel Core i9-12900K workstation, in IEEE double precision unless
stated otherwise. ``Standard'' refers to the usual full-grid Fourier
pseudo-spectral implementation; ``isotypic'' to Algorithm~\ref{alg:iso}, in
which the DCT-I/DST-I are realized by FFTs of the symmetric extensions of
each channel. The same integrator code is used for both representations;
only the transforms and the cubic term differ. To give the standard
implementation the best possible chance, symmetric initial data are
\emph{assembled from the channels} and are therefore exactly
(bitwise) symmetric for both representations. We use the defect
$D(u_n)=\max_{\bsig\notin\Sigma}\|(u_n)_{\bsig}\|_\infty$ and define the
breakdown time $T_\ast$ as the first \emph{recorded} time at which $D>10^{-2}$
(or $10^{-1}$ where stated), with diagnostic sampling interval
$\Delta t_{\mathrm{out}}=0.05$ for Table~\ref{tab:law} (and every macrostep
$T_e$ for ESP).

\subsection{Verification}
For random-looking non-symmetric data in $d=1,2,3$ ($N=64,32,16$,
$\kappa^2=0.1$, $\gamma=S=1$, $\tau=0.05$, 20 steps), the isotypic and
standard solutions differ by $3.3\times10^{-16}$ to $1.3\times10^{-15}$ across
IMEX1, ETDRK2, ETDRK4, and Crank--Nicolson (with unified physical-space
$\ell_\infty$ Picard stopping tolerance $10^{-12}$, up to 200 sweeps), in
agreement with Theorem~\ref{thm:identical}.
For symmetric data the transverse channels of the isotypic solution were
exactly zero in \emph{every} computation reported below (several million
time steps in total), in agreement with Theorem~\ref{thm:exact}.

\subsection{The benchmark of~\texorpdfstring{\cite{LQTY2022,Wang2023}}{LQTY2022 and Wang2023}: \texorpdfstring{$\kappa=1$, $u_0=\sin x$}{kappa = 1, u0 = sin x}}
We take $\kappa=1$, $\gamma=1$, $N=128$ as in~\cite[Ex.~5.1]{Wang2023}.
The symmetric solution decays like $a(t)\sin x$ with
$a(t)\sim(3t/2)^{-1/2}$; the time-varying transverse (even) linearization has
the asymptotic leading eigenvalue $1-\tfrac32a(t)^2+O(a(t)^4)\to1$ (with the
constant mode as the limiting eigenfunction as $a(t)\to0$), giving a limiting
growth rate $\lambda_\perp=1$ and leading-order breakdown scale
$\log(1/\umach)\approx37$ (with the slow $a(t)^2\sim 2/(3t)$ decay introducing
a subleading $\log\log(1/\delta)$ shift beyond the constant-coefficient
setting of Theorem~\ref{thm:breakdown}(ii)).

\begin{figure}[t]
\centering
\includegraphics[width=\textwidth]{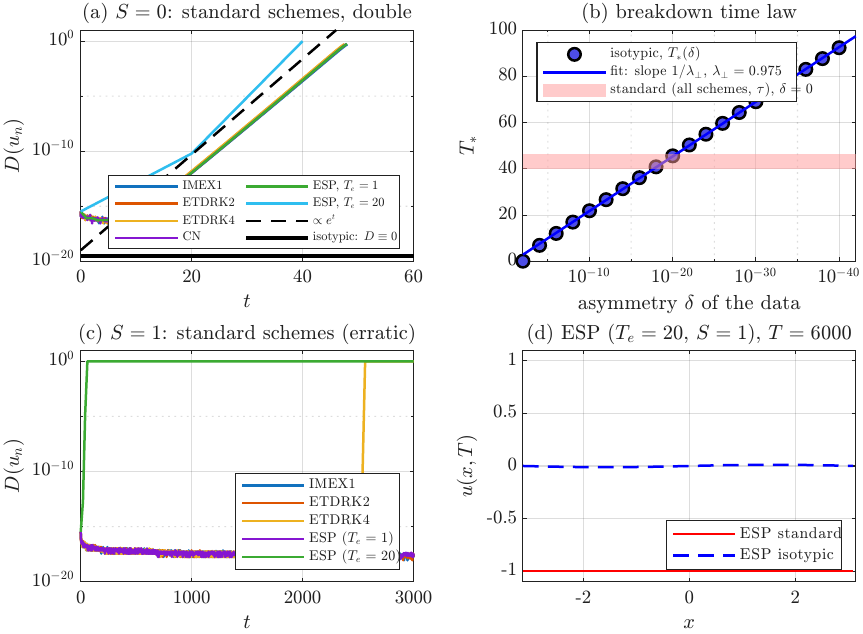}
\caption{$\kappa=1$, $u_0=\sin x$, $N=128$. (a) Defect of standard schemes, $S=0$,
$\tau=10^{-2}$ (ESP: $T_e=1,20$); the isotypic schemes give $D\equiv0$.
(b) Breakdown time of the isotypic ETDRK4 scheme versus a genuine asymmetry $\delta$
of the data (circles), affine fit, and the band of all standard runs of
Table~\ref{tab:law}. (c) Standard schemes with $S=1$ up to $T=3000$.
(d) Final states of ESP ($T_e=20$, $S=1$) at $T=6000$: standard (converged to the
wrong state $u\equiv-1$) and isotypic (odd).}\label{fig:wang}
\end{figure}

\paragraph{Without stabilization ($S=0$)} Table~\ref{tab:law} lists the breakdown
times of the standard implementation for four integrators with three step
sizes each and for ESP with $N_s=24$, $N_{\rm Iter}=20$ and $T_e=1,5,20$.
All 15 runs cross $D>10^{-2}$ at first recorded times $T_\ast\in[40,46.4]$
(for ESP with $T_e=20$, the threshold is already exceeded at the macrostep
$t=40$) and converge to the wrong steady state $u\equiv\pm1$, see
Figure~\ref{fig:wang}(a): the breakdown time is nearly independent of the
integrator, of its order, of the step size (over two orders of magnitude) and
of the ESP parameters, consistent with Theorem~\ref{thm:breakdown}. The
isotypic implementations of these integrators keep $D\equiv0$ throughout.

\begin{table}[t]
\centering
\caption{Breakdown times $T_\ast$ ($D>10^{-2}$) of the standard implementation,
$\kappa=1$, $u_0=\sin x$, $N=128$, $S=0$, double precision. The isotypic
implementation gives $D\equiv0$ up to $T=300$ in all cases.}\label{tab:law}
\begin{tabular}{lccc}
\toprule
scheme & $\tau=10^{-3}$ & $\tau=10^{-2}$ & $\tau=10^{-1}$\\
\midrule
IMEX1 (stabilized semi-implicit) & 43.60 & 43.90 & 46.40\\
ETDRK2 & 43.60 & 43.45 & 44.40\\
ETDRK4 & 43.50 & 43.35 & 43.80\\
Crank--Nicolson (Picard, tol $10^{-13}$) & 43.55 & 43.75 & 42.10\\
\midrule
 & $T_e=1$ & $T_e=5$ & $T_e=20$\\
ESP, $N_s=24$, $N_{\rm Iter}=20$ & 44.0 & 45.0 & $\le 40$\\
\bottomrule
\end{tabular}
\end{table}

\paragraph{The breakdown-time law} To measure the effective transverse growth
independently we add a genuine constant perturbation $\delta$ to the data and
integrate with the isotypic ETDRK4 scheme ($S=0$, $\tau=10^{-2}$, threshold
$D>10^{-2}$). Figure~\ref{fig:wang}(b) shows that $T_\ast(\delta)$ scales
nearly linearly with $-\log\delta$ over the finite window
$10^{-40}\le\delta\le10^{-4}$, with effective slope
$1/\lambda_\perp^{\mathrm{eff}}$ ($\lambda_\perp^{\mathrm{eff}}=0.975$).
Asymptotically, because $u^\ast(x,t)\approx a(t)\sin x$ decays slowly with
$a(t)^2\approx 2/(3t)$, averaging the linearized operator
$\partial_{xx}+1-3a(t)^2\sin^2 x$ against the limiting constant mode gives an
instantaneous growth rate $1-\tfrac32 a(t)^2+O(a^4)\approx 1-1/t$, whose
integral $t-\log t$ introduces a subleading $\log\log(1/\delta)$ correction to
the constant-coefficient law of Theorem~\ref{thm:breakdown}. Consistent with
this asymptotic correction, on the separate extended faithfulness scan of
Section~\ref{sec:faithful} ($S=1$, $\tau=0.05$, $D>0.1$, across the 17 samples
with $\delta\le 10^{-4}$ up to $\delta=10^{-200}$), replacing a free affine
fit in $x=-\log\delta$ (RMS residual $0.47$) by $x+\log x+c$ reduces the RMS
residual by more than an order of magnitude to $0.04$. The standard runs of
Table~\ref{tab:law} fall into the band $T_\ast\approx43.5\pm3$, corresponding
to an \emph{effective} seed $\delta_{\rm eff}\approx10^{-19}$, i.e.\ a
fraction of $\umach$ of the local values near the zeros $x=0,\pi$ of $u$,
where the floating-point grid is finest.

\paragraph{With stabilization ($S=1$, the setting of~\cite{Wang2023})} The
picture changes in an instructive way (Figure~\ref{fig:wang}(c)). The standard
IMEX1, ETDRK2 and ESP ($T_e=1$) runs keep $D\approx3\times10^{-16}$ up to
$T=3000$, whereas ETDRK4 suddenly breaks at $t\approx2560$ and ESP with
$T_e=20$---the large-step setting advertised in~\cite[Fig.~9(b)]{Wang2023}---breaks
at $t\approx60$ and converges to $u\equiv-1$ (Figure~\ref{fig:wang}(d)).
This is the sub-resolution absorption of Remark~\ref{rem:floor}: with
$S>0$ the zero mode is damped exactly in spectral space, while its
amplification, which enters only through the nonlinear term evaluated in
physical space, is lost as long as the mode is below the resolution of the
floating-point grid near $|u|$. Whether a run escapes from this regime depends
on fine details of the rounding; ``symmetry preservation'' by a standard
scheme is therefore a fragile, implementation-dependent event and not a
property of the scheme. In particular, the apparent symmetry preservation of
ESP reported in~\cite{Wang2023} is not robust with respect to $T_e$. We also remark that in single precision none of the standard $S=0$ runs
exceeded $D>10^{-2}$ within the observation window $T=300$, although
Theorem~\ref{thm:breakdown}(i) allows an earlier breakdown: here the rounding
does not satisfy the discounted-seed condition~\eqref{eq:seedlb}. The isotypic
scheme, in contrast, behaves identically in all these situations: $D\equiv0$.

\subsection{Universality across seven additional high-order time integrators}\label{sec:highorder}
To verify directly that higher temporal order does not prevent long-time
symmetry distortion---and that the isotypic representation acts as a
drop-in backend for any scheme---we test seven further widely used high-order
integrators on the same 1D benchmark ($\kappa=\gamma=1$, $S=0$, $u_0=\sin x$,
$N=128$; summary order/equivalence tests in \texttt{exp\_highorder.m} with
two-step-size order estimates at $t=0.5$ and $T=300$, plus long-time
trajectories to $T=1000$ in \texttt{exp\_highorder\_evol.m}): 2nd-order
BDF2-SAV~\cite{ShenXuYang2019} ($\tau=0.01$, $C_0=10$, measured order $1.99$), 3rd-order 4-stage $L$-stable IMEX Runge--Kutta
ARS(4,4,3)~\cite{AscherRuuthSpiteri1997} ($\tau=0.02$, order $2.91$), 3rd- and
4th-order IMEX-BDF (SBDF3, SBDF4)~\cite{AscherRuuthWetton1995} ($\tau=0.01$,
orders $2.92$ and $3.78$), 4th-order integrating-factor Runge--Kutta (IFRK4,
$\tau=0.02$, order $4.03$), 4th-order Gauss--Lobatto spectral deferred
correction (SDC4, 3 nodes, 5 sweeps, $\tau=0.05$, order
$4.00$)~\cite{DuttGreengardRokhlin2000,FengTangYang2015}, and 2nd-order Strang
operator splitting using the exact odd nonlinear flow
$\Phi^{\rm NL}_\tau(u)=u/\sqrt{u^2+(1-u^2)\ee^{-2\gamma\tau}}$ via
Walsh--Hadamard orbit evaluation ($\tau=0.02$, order $1.98$). On non-symmetric
1D data, every drop-in isotypic implementation reproduces its standard FFT
counterpart to $1.1\times10^{-16}$--$2.3\times10^{-15}$ (50 steps).

\begin{figure}[t]
\centering
\includegraphics[width=\textwidth]{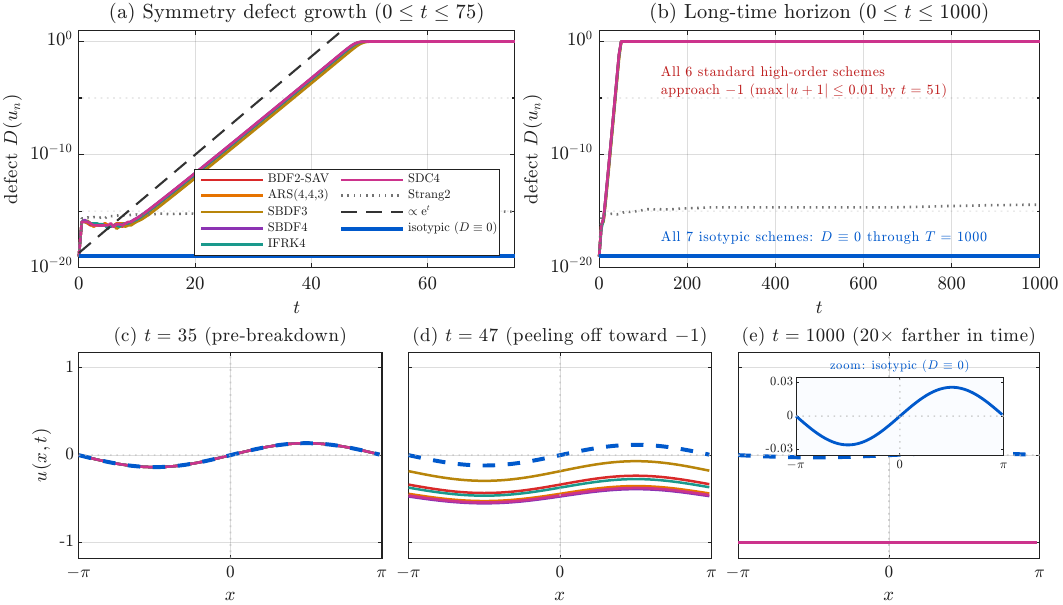}
\caption{Time-evolving departure of seven widely used high-order integrators
on 1D Allen--Cahn ($\kappa=\gamma=1$, $S=0$, $u_0=\sin x$, $N=128$,
$0\le t\le 1000$): BDF2-SAV, ARS(4,4,3), SBDF3, SBDF4, IFRK4, SDC4 and Strang2.
(a)~Defect growth over $0\le t\le 75$: all six standard polynomial high-order
schemes grow like $\ee^t$ and break at $T_\ast\in[43.0,44.1]$.
(b)~Long-time horizon $0\le t\le 1000$: all seven drop-in isotypic schemes stay
at $D\equiv0$ exactly.
(c)--(e)~Solution snapshots at $t=35$, $t=47$ (all six standard high-order
curves visibly peeling off toward $-1$), and $t=1000$ ($20\times$ farther in
time, with zoom inset confirming the exact odd isotypic profile).}\label{fig:highorder}
\end{figure}

Figure~\ref{fig:highorder} displays the time-evolving departure up to $T=1000$.
All six standard polynomial-nonlinearity high-order schemes amplify rounding
defects in lockstep like $\ee^t$ and break at $T_\ast\in[43.0,44.1]$
(BDF2-SAV $43.4$, ARS(4,4,3) $43.1$, SBDF3 $44.1$, SBDF4 $43.0$, IFRK4 $43.3$,
SDC4 $43.0$), visibly peeling away from the odd wave by $t=47$. By $t=48$,
all six standard schemes have negative values at every grid point (with
$\max_j|u_j+1|\in[0.20,0.64]$), and their first recorded times within
$10^{-2}$ of the wrong steady state $u\equiv -1$ in the maximum norm range
from $t=50.0$ to $51.0$ (while standard Strang2 still generates nonzero
rounding defects of $2.2\times10^{-15}$). By
contrast, plugging our isotypic spectral backend underneath \emph{every single
one} of these seven high-order integrators keeps $D(u_n)\equiv 0$ exactly all
the way to $T=1000$---more than $20\times$ beyond the standard breakdown
horizon---preserving the $\sin x$-dominated odd profile (peak amplitude
$\approx 0.0257\text{--}0.0258$ at $T=1000$, with small odd higher harmonics
from the cubic term) with zero symmetry drift.

\subsection{Faithfulness: genuine sub-roundoff asymmetries}\label{sec:faithful}
\begin{figure}[t]
\centering
\includegraphics[width=\textwidth]{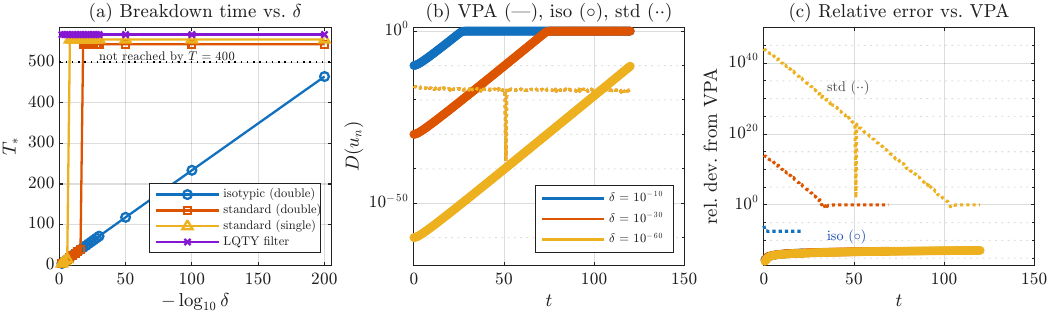}
\caption{Faithfulness to genuine sub-roundoff asymmetries ($\kappa=1$,
$u_0=\sin x+\delta$). (a)~Computed breakdown time $T_\ast$ ($D>0.1$, ETDRK4,
$S=1$) versus $\delta\in[10^{-200},10^{-2}]$ for isotypic double, standard
double/single, and the parity filter~\cite{LQTY2022}. (b)~Defect histories
against a 120-digit VPA reference ($N=32$, ETDRK2, $\tau=0.05$) for
$\delta=10^{-10},10^{-30},10^{-60}$. (c)~Relative deviation from VPA in the
linear regime ($D<10^{-2}$): isotypic double stays near $10^{-14}$ while
standard double is off by up to $10^{44}$.}\label{fig:faithful}
\end{figure}

A good symmetry-preserving method should not only keep exact symmetries but
also resolve genuine small asymmetries correctly. We perturb $u_0=\sin x$ by
$\delta$ ($\delta=10^{-2},\dots,10^{-200}$) and integrate with ETDRK4, $S=1$,
$\tau=0.05$ using (i) the isotypic scheme, (ii) the standard scheme in double and
(iii) single precision, and (iv) the standard scheme with the parity filter of
\cite{LQTY2022}, which replaces $\widehat u$ by $\mathrm i\,\mathrm{Im}\,\widehat u$ after
every step. Figure~\ref{fig:faithful}(a) shows the first recorded breakdown times
($D>0.1$, main scan horizon $T=400$, with \texttt{exp\_delta.m} also running
an extended $T=500$ continuation for isotypic at $\delta=10^{-200}$).
The isotypic scheme follows the logarithmic law over the whole range, including
$\delta=10^{-100}$ ($T_\ast=233.7$) and the empirical record at
$\delta=10^{-200}$ ($T_\ast=464.7$; note that while $\delta^3=10^{-600}$
underflows benignly to zero outside the formal no-underflow hypothesis of
Theorem~\ref{thm:relative}, the dominant linear term $u_\Sigma^2 u_\perp$
remains well within the IEEE double normal range). Under these parameters and
sampling, the standard double-precision scheme reproduces $T_\ast(\delta)$
only down to $\delta\ge10^{-16}$ and does not reach $D>0.1$ within $T=400$ for
$\delta\le10^{-18}$ (right-censored in Figure~\ref{fig:faithful}(a)); in
single precision the corresponding observed cutoff is $\delta\approx10^{-6}$.
The filtered scheme never breaks by construction---erasing every genuine
$\delta>0$.

To validate the isotypic results independently, we computed reference
solutions for $N=32$, ETDRK2, $\tau=1/20$, $u_0=\sin x+\delta$ with MATLAB's
variable-precision arithmetic, with \emph{all} quantities (including $\pi$,
$\tau$, $\delta=10^{-p}$, ETD multipliers and the DFT matrix) constructed
symbolically with 120 significant digits. Figure~\ref{fig:faithful}(b)--(c)
shows that, in the linear regime $D<10^{-2}$, the isotypic double-precision
defect agrees with the 120-digit reference to a relative accuracy of
$2.9\times10^{-14}$ for $\delta=10^{-10}$, $8.7\times10^{-14}$ for
$\delta=10^{-30}$ and $1.5\times10^{-13}$ for $\delta=10^{-60}$, while the
standard double-precision computation is off by factors $10^{14}$ and $10^{44}$
for $\delta=10^{-30}$ and $10^{-60}$. This confirms Theorem~\ref{thm:relative}
and Corollary~\ref{cor:faithful}: the transverse dynamics are computed with
\emph{relative} accuracy, independently of the size of the transverse part.

\subsection{The adaptive-cG benchmark (\texorpdfstring{$\eps=1$, $N=4096$}{eps = 1, N = 4096}, Crank--Nicolson) and multi-kink dynamics}
\begin{figure}[t]
\centering
\includegraphics[width=\textwidth]{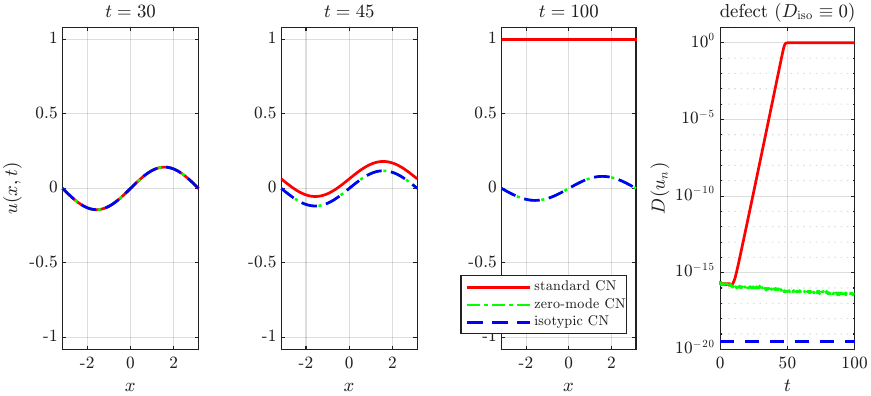}
\caption{Crank--Nicolson ($\tau=10^{-2}$, $N=4096$, $\kappa=\gamma=1$,
$u_0=\tfrac12(\sin x+\sin 4x)$): snapshots at $t=30,45,100$ and defect
histories for standard CN, zero-mode-suppressed CN, and isotypic CN
($D\equiv0$, drawn at the bottom of the axis).}\label{fig:draft57}
\end{figure}

In our companion work on adaptive continuous-Galerkin time stepping, the
one-dimensional test $\partial_tu=\Delta u-(u^3-u)$ ($\kappa=\gamma=1$,
$S=0$) with $u_0(x)=\tfrac12(\sin x+\sin 4x)$, $N=4096$, Crank--Nicolson
($\tau=10^{-2}$, fixed-step Picard iteration up to 200 sweeps, physical-space
$\ell_\infty$ stopping tolerance $10^{-10}$, monitored every $0.1$ time units)
and $T=100$ was stabilized by an a~posteriori suppression of the zero Fourier
mode. Figure~\ref{fig:draft57} compares the standard CN scheme,
zero-mode-suppressed CN, and the isotypic CN scheme (with no zero-mode
suppression). The standard CN solution develops a defect $D=2.9\times10^{-8}$
by $t=30$, lifts visibly above the odd sine wave at $t=45$
($D=6.4\times10^{-2}$), and collapses to the wrong steady state $u\equiv 1$
well before $t=100$. Zero-mode suppression keeps the odd wave visually intact
with residual defect $D\approx 2.7\times10^{-16}$. The isotypic CN scheme
achieves $D(u_n)\equiv0$ for all $10^4$ steps without modifying any Fourier
mode, preserving odd reflection symmetry with final range $[-0.0806,0.0806]$
at $T=100$.

In the thin-interface multi-kink regime ($\partial_tu=\Delta u-\eps^{-2}f(u)$,
$\eps=0.05$, $N=1024$, $S=2\eps^{-2}$, $\tau=10^{-3}$, $T=100$), standard
IMEX1, ETDRK2, ETDRK4 and Crank--Nicolson all exhibit small neutral-mode
drifts (reaching $\max D\approx 3.0\text{--}6.0\times10^{-12}$ by $T=100$ in
\texttt{explore1d.m}), whereas the isotypic schemes maintain $D\equiv0$.
Similarly, in the thin-interface benchmark of~\cite[Ex.~5.2(d)]{Wang2023}
($\kappa^2=10^{-3}$, $N=512$, $S=1$, $\tau=0.05$ for ETDRK4 and $T_e=10$ for
ESP), the defect peaks transiently at $3.5\times10^{-11}$ during initial
coarsening, decays as the interfaces settle, and then drifts slowly from a
lower baseline to $1.4\times10^{-13}$ (ETDRK4) and $1.8\times10^{-13}$ (ESP)
at $T=5000$, whereas the isotypic schemes maintain $D\equiv0$.

%======================================================================
\section{Numerical experiments in two and three dimensions}\label{sec:num23d}

In this section we use the scaling
$\partial_tu=\Delta u-\eps^{-2}f(u)$ on $[-\pi,\pi)^d$ with ETDRK2 and the
stabilization parameter $S=2\eps^{-2}$~\cite{DJLQ2021}, and compare three
implementations: the standard Fourier scheme, the standard scheme with the
zero-mode suppression $\widehat u(\mathbf 0)\gets0$ after every step, and the
isotypic scheme.

\subsection{Is the standard FFT equivariant? It depends on \texorpdfstring{$N$}{N}}
Table~\ref{tab:fftN} reports the largest defect produced by the standard
implementation during the first five steps ($\tau=10^{-3}$, $\eps=0.05$) for
bitwise-symmetric data in four symmetry classes (3D lengths $N\ge 160$ were
not sampled). In one dimension a defect of size $\umach$ is created for every
$N$. In two and three dimensions, however, for these tested inputs, FFTW plans
and five-step runs, MATLAB's FFTW happens to yield zero defect for
$N\in\{32,64,128\}$, and nonzero $O(\umach)$ defects for
$N\in\{48,96,160,192,256\}$. Whether a standard multi-dimensional computation
preserves a symmetry thus depends sensitively on the factorization of $N$
chosen by the FFT planner---a striking illustration of the fragility discussed
in Remark~\ref{rem:floor} and of the need for a structural solution. All
isotypic computations give $D\equiv0$ irrespective of $N$.

\begin{table}[t]
\centering
\caption{Largest symmetry defect created by the standard FFT implementation in the
first five ETDRK2 steps, for bitwise symmetric data ($\eps=0.05$,
$\tau=10^{-3}$). ``0'' means exactly zero.}\label{tab:fftN}
\small
\begin{tabular}{lcccccccc}
\toprule
symmetry class & $N=32$ & 48 & 64 & 96 & 128 & 160 & 192 & 256\\
\midrule
1D odd & 3.3e-16 & 3.3e-16 & 5.3e-16 & 4.4e-16 & 5.0e-16 & 4.0e-16 & 4.4e-16 & 3.3e-16\\
2D odd--odd & 0 & 1.1e-15 & 0 & 1.4e-15 & 0 & 1.7e-15 & 1.6e-15 & 6.9e-16\\
2D point-odd ($C_2$) & 0 & 7.6e-16 & 0 & 5.8e-16 & 0 & 6.7e-16 & 8.9e-16 & 7.8e-16\\
3D odd$^3$ & 0 & 1.9e-15 & 0 & 2.0e-15 & 0 & -- & -- & --\\
\bottomrule
\end{tabular}
\end{table}

\subsection{Two dimensions: an X-junction}
\begin{figure}[t]
\centering
\includegraphics[width=\textwidth]{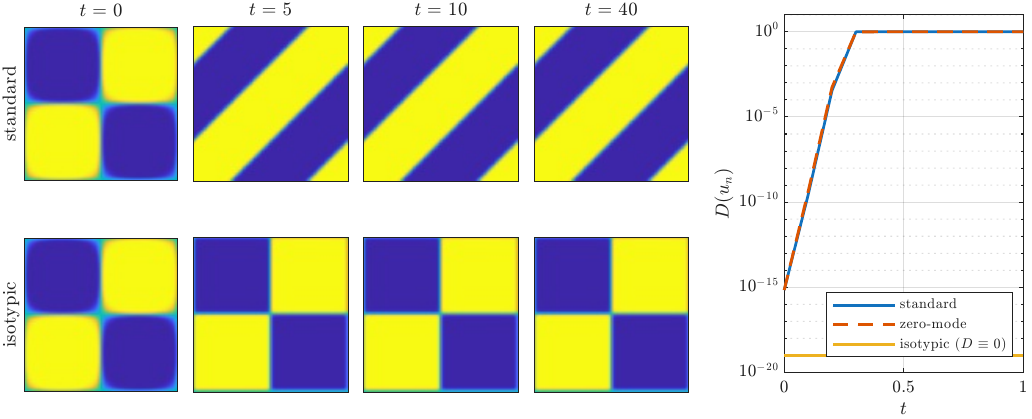}
\caption{X-junction, $\eps=0.05$, $N=256$, ETDRK2, $\tau=10^{-3}$. Snapshots of the
standard (top) and isotypic (bottom) solutions, and defect histories (right; the
isotypic defect is identically zero and is drawn at the bottom of the axis).}\label{fig:X2d}
\end{figure}

The datum $u_0=\tanh\big(\sin x\sin y/(\sqrt2\eps)\big)$, $\eps=0.05$, is odd in
$x$ and in $y$ ($\Sigma=\{(1,1)\}$). The symmetric evolution keeps the four straight
interfaces $x=0$, $x=\pi$, $y=0$, $y=\pi$, which meet in four X-junctions
at $(0,0)$, $(0,\pi)$, $(\pi,0)$ and $(\pi,\pi)$; the junctions cannot move by
symmetry, so that the four-phase checkerboard is a symmetric steady state. It is, however,
strongly unstable with respect to reconnection of the interfaces at the
junctions. With $N=256$, $\tau=10^{-3}$, the standard computation creates a
defect $7\times10^{-16}$ already in the first transform, which grows to
$2.6\times10^{-10}$ at $t=0.1$ and $3.5\times10^{-4}$ at $t=0.2$
(transverse growth rate $\lambda_\perp\approx1.4\times10^{2}$), crosses
$D>10^{-2}$ shortly after $t=0.22$, and undergoes macroscopic reconnection by
$t\approx0.3$ into diagonal stripes, a state that is not odd in $x$ or $y$
(Figure~\ref{fig:X2d}). The zero-mode suppression is of no help---the
transverse linearization restricted to zero-mean perturbations remains
strongly unstable at the junctions---and breaks at the same time. The
isotypic scheme keeps $D\equiv0$ and the symmetric checkerboard up to $T=40$
($4\times10^4$ steps).

\subsection{Head-to-head 2D comparison on~\texorpdfstring{\cite[Ex.~4.1]{Wang2023}}{Wang2023 Ex. 4.1}: identical ESP parameters}
\begin{figure}[t]
\centering
\includegraphics[width=\textwidth]{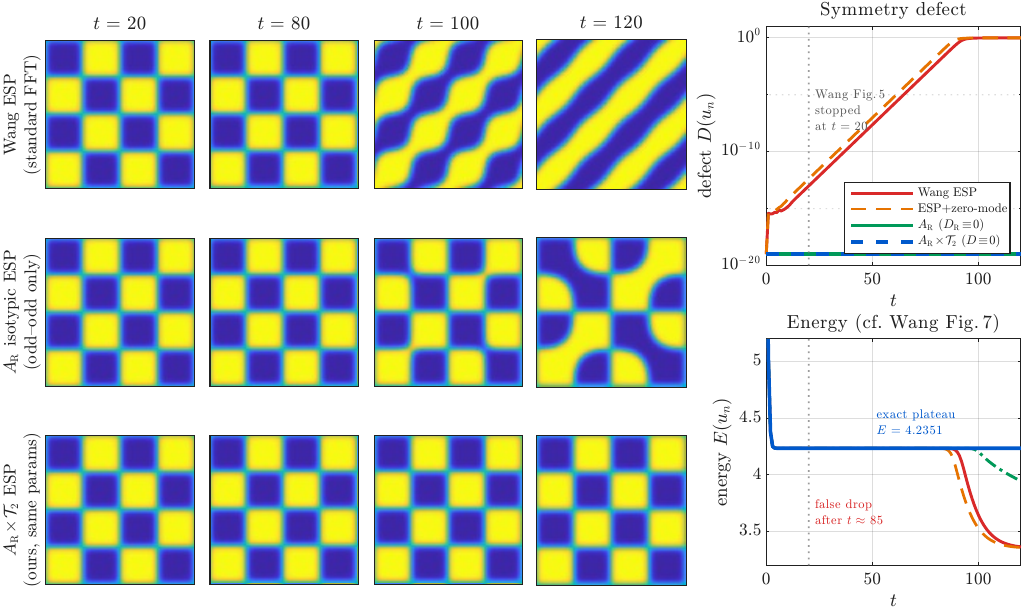}
\caption{Head-to-head 2D comparison on~\cite[Ex.~4.1]{Wang2023}
($\kappa=0.1$, $u_0=\sin(2x)\sin(2y)$, $N_p=256$) with \emph{identical} ESP
parameters ($S=1$, $T_e=1$, $N_s=10$, $N_{\rm Iter}=20$) extended from
Wang's stopping time $t=20$ to $t=120$. Top: Wang's standard FFT ESP breaks
odd--odd symmetry after $t\approx80$ and collapses into diagonal stripes
($D=0.997$). Middle: $A_{\mathrm R}$ isotypic ESP preserves odd--odd symmetry
($D_{\mathrm R}\equiv0$) while allowing $\pi$-periodicity-breaking
$4\times4\to2\times2$ coarsening. Bottom: full $A_{\mathrm R}\times\mathcal T_2$
isotypic ESP keeps all 15 transverse channels identically zero ($D\equiv0$),
relaxing from $E(0)=6.7175$ to the symmetric $4\times4$ checkerboard plateau
$E=4.2351$ for $t\in[20,120]$
(right panels: defect and energy histories extending~\cite[Fig.~7, left]{Wang2023}).}\label{fig:wang41}
\end{figure}

To make a strictly controlled comparison in two dimensions---holding the PDE,
initial datum, spatial resolution, time integrator and \emph{every} algorithmic
parameter identical and changing \emph{only} the spectral representation---we
revisit Wang's 2D Allen--Cahn benchmark~\cite[Ex.~4.1, Fig.~5]{Wang2023}:
$\partial_tu=\kappa^2\Delta u-(u^3-u)$ on $[-\pi,\pi)^2$ with $\kappa=0.1$,
$u_0(x,y)=\sin(2x)\sin(2y)$, $N_p=256\times256$, and ESP with $S=1$, $T_e=1$,
$N_s=10$, $N_{\rm Iter}=20$. In~\cite{Wang2023} the simulation was stopped at
$t=20$, where the $4\times4$ checkerboard still appears intact. Because
$N_p=256$ is a non-equivariant 2D FFT length (Table~\ref{tab:fftN}), however,
Wang's standard FFT ESP creates an odd--odd defect of $3.9\times10^{-16}$
already at $t=1$, which grows exponentially to $9.8\times10^{-14}$ at $t=20$
and $5.4\times10^{-3}$ at $t=80$; by $t=90$--$100$ it completely destroys the
odd--odd symmetry ($D=0.997$ at $t=120$), reconnecting the checkerboard into
asymmetric diagonal stripes with an energy drop triggered by departure from
the symmetric branch, from the plateau $E=4.2351$ down to $3.3648$
(Figure~\ref{fig:wang41}, top row). Adding a~posteriori zero-mode suppression
makes the breakdown even faster ($D=0.746$ at $t=90$, $E=3.3589$ at $t=120$).

Remarkably, $u_0(x,y)=\sin(2x)\sin(2y)$ carries \emph{two} commuting abelian
symmetries: odd--odd reflection ($A_{\mathrm R}$) and $\pi$-periodicity
($\mathcal T_2$), generating the order-16 group
$G=A_{\mathrm R}\times\mathcal T_2\cong\Z_2^4$ (16 channels on a
$(N/4+1)\times(N/4+1)$ fundamental grid, whose non-symmetric 2D equivalence
with standard Fourier is $7.8\times10^{-16}$ after 20 ETDRK2 steps). Running
the \emph{exact same} ESP integrator ($S=1$, $T_e=1$, $N_s=10$,
$N_{\rm Iter}=20$) on both isotypic backends reveals the exact group-theoretic
structure of Theorem~\ref{thm:exact} (Figure~\ref{fig:wang41}):
\begin{itemize}
\item the reflection-only $A_{\mathrm R}$ isotypic ESP preserves the odd--odd
  reflection symmetry \emph{exactly} ($D_{\mathrm R}\equiv0$ throughout, so
  the zero lines $x=0$ and $y=0$ never break), while rounding inside the
  odd--odd subspace breaks the unprotected $\pi$-periodicity around
  $t\approx100$--$120$ (coarsening the $4\times4$ checkerboard toward a
  $2\times2$ odd--odd checkerboard, $E\to3.9452$);
\item the full $A_{\mathrm R}\times\mathcal T_2$ isotypic ESP preserves
  \emph{all 15 transverse channels identically zero} ($D\equiv0$) without any
  filtering or intervention: the energy relaxes smoothly from $E(0)=6.7175$ to
  $E(5)=4.2352$ and holds the symmetric $4\times4$ checkerboard plateau
  $E=4.2351$ throughout $t\in[20,120]$.
\end{itemize}

\subsection{Two dimensions: a point-odd (\texorpdfstring{$C_2$}{C2}) pattern over long times}
\begin{figure}[t]
\centering
\includegraphics[width=\textwidth]{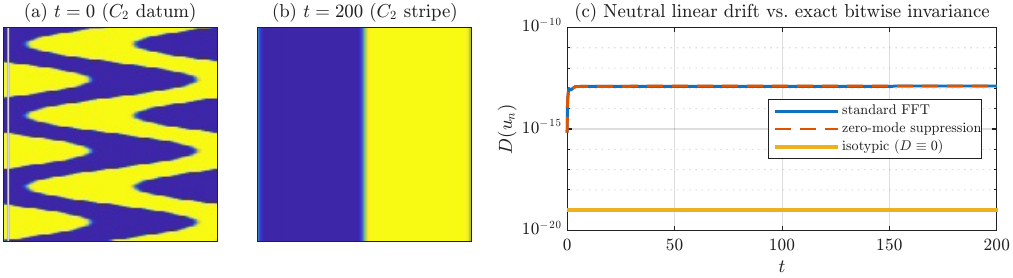}
\caption{Point-odd ($C_2$) pattern, $\eps=0.05$, $N=192$, ETDRK2,
$\tau=2\times10^{-3}$: (a)~initial datum, (b)~coarsened $C_2$-odd stripe
at $t=200$, and (c)~defect histories showing early growth to $7.4\times10^{-14}$
followed by slow drift to a plateau of $1.31\times10^{-13}$ for both standard
FFT and zero-mode suppression versus exact invariance ($D\equiv0$) for the
isotypic scheme.}\label{fig:C2d}
\end{figure}

The datum $u_0=\tanh\big(\sin(x+2.5\sin3y)/(\sqrt2\eps)\big)$ of our companion study
is odd under the point reflection $(x,y)\mapsto(-x,-y)$ but has no reflection symmetry
in a single variable; hence $\Sigma=\{(1,0),(0,1)\}$ and the channels $(0,0)$ and $(1,1)$
must vanish. We use $N=192$, $\tau=2\times10^{-3}$, $T=200$ ($10^5$ steps).
The wavy interfaces straighten and coarsen into a point-odd ($C_2$) stripe
pattern (Figure~\ref{fig:C2d}). The defect of the standard scheme grows
rapidly during the initial interface relaxation (reaching $7.4\times10^{-14}$
by $t=0.6$) and then drifts slowly to a plateau of $1.31\times10^{-13}$ around
$t\approx175$ as the stripe pattern becomes stationary; zero-mode suppression
behaves identically, because the transverse defect lives in all modes of the
$(0,0)$ and $(1,1)$ channels, not only in the mean. The isotypic scheme keeps
both channels identically zero. (For $N=128$, which gives zero defect in
Table~\ref{tab:fftN}, all three implementations give $D\equiv0$ up to $T=400$.)

\subsection{Three dimensions: an octant saddle}
\begin{figure}[t]
\centering
\includegraphics[width=\textwidth]{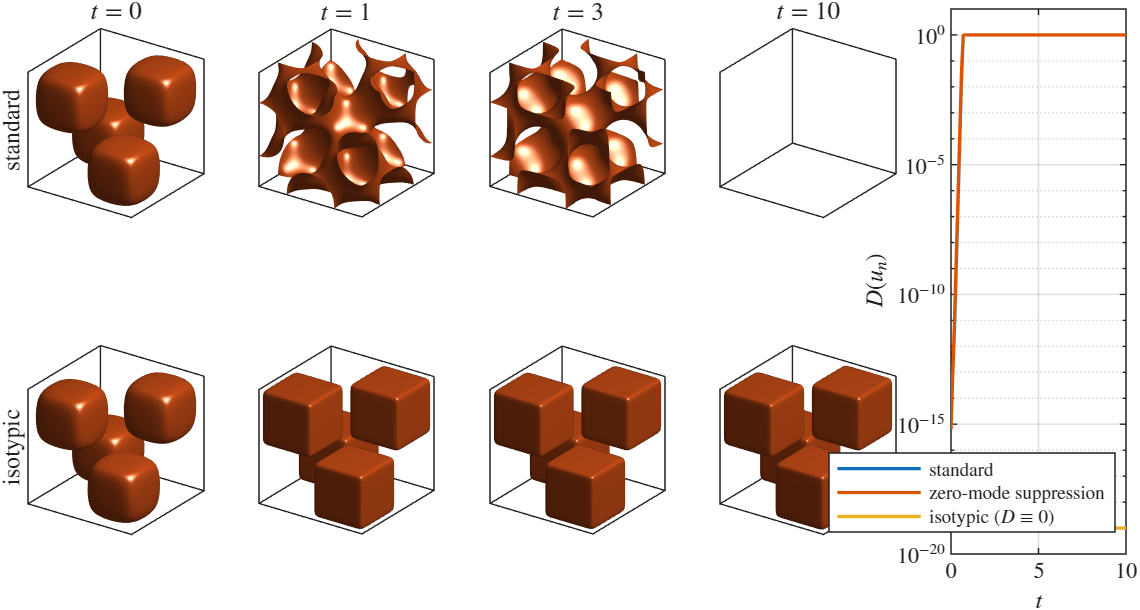}
\caption{Octant saddle in 3D, $\eps=0.1$, $N=96$, ETDRK2, $\tau=2\times10^{-3}$: level set
$u=0.5$ of the standard (top) and isotypic (bottom) solutions, and defect
histories (right).}\label{fig:O3d}
\end{figure}

The datum $u_0=\tanh\big(\sin x\sin y\sin z/(\sqrt2\eps)\big)$, $\eps=0.1$, is odd in all
three variables ($\Sigma=\{(1,1,1)\}$, seven transverse channels); the symmetric
evolution is the three-dimensional analogue of the X-junction, a
checkerboard of eight octants separated by the planes $x_i\in\{0,\pi\}$. With
$N=96$ (an FFT length for which the standard implementation is not equivariant,
Table~\ref{tab:fftN}), $\tau=2\times10^{-3}$, the defect of the standard scheme grows
from $6.4\times10^{-16}$ at $t=0$ by roughly two orders of magnitude per $0.1$ time
units ($\lambda_\perp\approx55$): $5.3\times10^{-12}$ at $t=0.2$, $3.3\times10^{-7}$ at $t=0.4$,
$2.1\times10^{-2}$ at $t=0.6$; at $t\approx0.7$ the octant structure has collapsed
(Figure~\ref{fig:O3d}). Again the zero-mode suppression changes nothing. The
isotypic scheme preserves all seven vanishing channels exactly up to $T=10$.

\subsection{A non-polynomial nonlinearity: the Flory--Huggins potential}
To test Proposition~\ref{prop:orbit} we replace $f$ by the logarithmic
nonlinearity $F(u)=\tfrac\theta2\big(\mathrm{log1p}(u)-\mathrm{log1p}(-u)\big)-\theta_cu$
on $|u|<1$, with $\theta=0.8$, $\theta_c=1.6$ (binodal values $\pm0.9575$),
evaluated in the odd-exact form shown (no clamping needed as the stabilized
ETDRK2 iterates remain strictly inside $(-1,1)$), and use the orbit evaluation
in the isotypic scheme. For non-symmetric data the isotypic and standard ETDRK2
solutions agree to $4.4\times10^{-16}$ (1D, $N=64$) and $1.0\times10^{-15}$
(2D, $N=32$) after 40 steps. For $\kappa=1$, $\gamma=1$, $S=0$,
$u_0=0.9\sin x$, $N=128$, $\tau=10^{-2}$ the zero state is unstable in the
constant direction with rate $\theta_c-\theta=0.8$; the standard scheme breaks
the symmetry at $T_\ast=50.7$ and converges to the wrong state
$u\equiv-0.9575$, while the isotypic scheme keeps $D\equiv0$ and decays to the
correct odd limit $u\equiv0$ (at $T=400$, $\max|u|=1.15\times10^{-35}$). In
2D ($\eps=0.05$ scaling, $u_0=0.9\sin x\sin y$, $N=96$, $\tau=10^{-3}$,
$S=8\eps^{-2}$, horizon $T=5$) the standard scheme breaks the odd--odd
symmetry at $T_\ast=0.45$, whereas the isotypic scheme keeps $D\equiv0$
throughout the simulated horizon.

\subsection{Cost}
Table~\ref{tab:cost} reports \emph{paired} median wall-clock times per ETDRK4
step over 5 repeated trials (four nonlinear evaluations, 4 threads, warm-up
excluded) starting from the same analytic initial functions across standard
Fourier, reflection isotypic ($A_{\mathrm R}$), and half-period translation
isotypic ($\mathcal T_2$) solvers (\texttt{exp\_cost.m}; standard Fourier and
$A_{\mathrm R}$ assemble bitwise-identical full-grid initial arrays, while
$\mathcal T_2$ initializes directly on its half-period fundamental orbit,
agreeing with standard Fourier to $2.0\times 10^{-15}$ in 2D $N=256$).
``Full'' evolves all $2^d$ channels (with no zero skipping, representing
worst-case cost for arbitrary non-symmetric data); ``reduced'' evolves only
the pre-selected active channel in $\Sigma$ via an a~priori mask
(\texttt{P.active}; automatic runtime detection of exact-zero channels would
incur a small additional scan overhead not timed here). By Theorem~\ref{thm:exact},
reduced and full computations agree identically
(\texttt{isequal(\allowbreak U\_full,\allowbreak U\_reduced)} verified
\texttt{true} in every case, checking signed exact equality up to IEEE
$\pm 0.0$). In our MATLAB prototype the full $A_{\mathrm R}$ scheme costs
$1.78\text{--}4.86\times$ standard FFT and full $\mathcal T_2$ costs
$3.03\text{--}5.65\times$ (reflecting the $2\cdot 4^d$ channel products,
twiddle overheads, and lack of native real-to-real transforms in MATLAB; with
native FFTW REDFT00/RODFT00~\cite{FrigoJohnson2005}, the transform work of all
$2^d$ channels equals one real $N^d$ FFT). In favourable multi-dimensional configurations, reduced mode is faster than
standard FFT: comparing medians against the standard FFT median (and reporting
the paired trial-by-trial ratio range $\min_r/\max_r(t_{r,j}/t_{r,\mathrm{std}})$
across the 5 trials), $A_{\mathrm R}$ reduced achieves median ratio
$0.69\times$ (paired range $[0.59\times,0.76\times]$) in 2D ($N=256$), while
$\mathcal T_2$ reduced achieves median ratio $0.92\times$ (paired range
$[0.89\times,0.95\times]$) in 3D ($N=128$) and $1.01\times$ (paired range
$[0.99\times,1.04\times]$) in 3D ($N=64$). Measuring run-to-run variability by
the maximum absolute relative deviation from each solver's 5-trial column
median, $\max_{r,j}|t_{r,j}-\mathrm{med}_j|/\mathrm{med}_j$, gives
$4.5\%$ for 3D ($N=128$), $8.6\%$ for 3D ($N=64$), $16.9\%$ for 2D ($N=256$),
and $40.5\%$ for the sub-millisecond 1D ($N=4096$) case (note that the five
solvers execute sequentially in fixed order on a workstation, so small ratios
near $1\times$ remain sensitive to CPU load, thermal state, and FFTW
planning).

\begin{table}[t]
\centering
\caption{Paired median wall-clock time (seconds) per ETDRK4 step across 5
repeated trials (4 threads, same analytic initial functions, pre-selected
active mask \texttt{P.active} for reduced runs, fixed sequential order,
\texttt{exp\_cost.m}) for standard Fourier, reflection isotypic
($A_{\mathrm R}$), and half-period translation isotypic ($\mathcal T_2$).
Maximum absolute relative deviations from column medians across all five
solvers are $40.5\%$ (1D $N=4096$), $16.9\%$ (2D $N=256$), $8.6\%$ (3D
$N=64$), and $4.5\%$ (3D $N=128$). Exact signed equality (\texttt{isequal})
between full and reduced runs holds for both $A_{\mathrm R}$ and
$\mathcal T_2$ in all cases.}\label{tab:cost}
\small
\begin{tabular}{rrcccccc}
\toprule
$d$ & $N$ & standard & $A_{\mathrm R}$ (full) & $A_{\mathrm R}$ (reduced) & $\mathcal T_2$ (full) & $\mathcal T_2$ (reduced) & exact eq.\\
\midrule
1 & 4096 & 2.32e-4 & 1.13e-3 & 6.30e-4 & 9.68e-4 & 7.20e-4 & yes\\
2 & 256  & 5.91e-3 & 1.05e-2 & 4.05e-3 & 1.79e-2 & 6.70e-3 & yes\\
3 & 64   & 2.53e-2 & 6.71e-2 & 3.23e-2 & 9.02e-2 & 2.57e-2 & yes\\
3 & 128  & 1.96e-1 & 7.56e-1 & 2.19e-1 & 1.11e+0 & 1.80e-1 & yes\\
\bottomrule
\end{tabular}
\end{table}

%======================================================================
\section{Arbitrarily shifted symmetries in 1D, 2D and 3D}\label{sec:numshift}

A central limitation of origin-centred parity filters~\cite{LQTY2022}---and
of fixed-origin reflection isotypic schemes---is that they break down as soon
as the odd symmetry centre is shifted away from the origin. Wang~\cite[Ex.~5.2]{Wang2023}
highlighted shifted initial data $\tilde x=x+\theta$ ($\theta\in\{0,\pi/8,\pi/4\}$)
as a key test where filters cannot be applied directly. Here we test the
half-period translation scheme ($\mathcal T_2$, Theorem~\ref{thm:t2}) and the
canonical-frame isotypic scheme (CF-ISM, Algorithm~\ref{alg:cfism} and
Theorem~\ref{thm:cfism}) on both grid-aligned and \emph{off-grid} shifts
($\theta=0.37123456789$, not aligned with $h\Z/2$; the theory applies equally
to true irrational shifts).

\subsection{Verification of half-grid Type~II and \texorpdfstring{$\mathcal T_2$}{T2} representations}
For general non-symmetric data in $d=1,2,3$ (\texttt{test\_shifted.m}), the
half-grid DCT-II/DST-II representation (Section~\ref{sec:iso_refl}) agrees
with the standard Fourier multiplier to $2.8\times10^{-16}$ ($d=1$, $N=64$),
$4.4\times10^{-16}$ ($d=2$, $N=32$) and $4.4\times10^{-16}$ ($d=3$, $N=16$),
and the $\mathcal T_2$ representation (Section~\ref{sec:iso_t2}) agrees to
$1.9\times10^{-16}$ ($d=1$, $N=128$), $4.4\times10^{-16}$ ($d=2$, $N=64$) and
$4.4\times10^{-16}$ ($d=3$, $N=32$). For CF-ISM applied to completely
non-symmetric 1D data ($u_0=0.5\sin x+0.3\cos 2x+0.2\sin(3x+0.4)$, $N=128$,
40 ETDRK4 steps), shifting the canonical-frame solution back by $-\mathbf c$
matches the standard Fourier solution to $7.8\times10^{-16}$. This confirms
Theorem~\ref{thm:identical} for all three shifted architectures.

\subsection{Refutation of ESP on~\texorpdfstring{\cite[Ex.~5.2]{Wang2023}}{Wang2023 Ex. 5.2} and exact preservation by \texorpdfstring{$\mathcal T_2$}{T2} and CF-ISM}
Table~\ref{tab:wang52} tests all four cases of~\cite[Ex.~5.2]{Wang2023}
($N=512$) across four shifts: $\theta=0$, $\pi/8$, $\pi/4$ and the off-grid
shift $\theta=0.37123456789$.
\begin{itemize}
\item \textbf{Standard ESP ($S=1$, $N_s=24$, $N_{\rm Iter}=20$, $T_e=1,2,2,10$
  for (a)--(d)):} Even at the parameter setting for Ex.~5.2(a) with
  $\theta=\pi/4$ (corresponding to~\cite[Fig.~10(a), right]{Wang2023}), ESP
  crosses $|\mathrm{mean}|>10^{-2}$ at $T_\ast=41$ and collapses to $u\equiv 1$.
  For the off-grid shift $\theta=0.3712$, ESP breaks at $T_\ast=38$ in
  Ex.~5.2(a) and at $T_\ast=56$ in Ex.~5.2(b). In the thin-interface cases
  Ex.~5.2(c),(d), ESP exhibits mean drift up to $9.9\times10^{-11}$ and
  $6.0\times10^{-9}$. Standard ETDRK4 ($S=0$, $\tau=0.02,0.02,0.02,0.05$)
  breaks consistently at $T_\ast\approx 37$--$43$ (Ex.~5.2(a)) and
  $T_\ast\approx 54$--$65$ (Ex.~5.2(b)) for all four shifts.
\item \textbf{$\mathcal T_2$ isotypic ($S=0$ in Table~\ref{tab:wang52}; pure
  double precision, parity-consistent channel initialization, zero shift
  detection):} For all odd-harmonic cases (Ex.~5.2(a),(b)), $\mathcal T_2$
  maintains exact zero even-wavenumber channels ($D_{\mathrm{anti}}\equiv0$)
  across every shift including $\theta=0.3712$. (By contrast, naive full-orbit
  double-precision sampling at $N=512$ injects $3.9\times 10^{-16}$ into the
  even channel before $t=0$, consistent with Proposition~\ref{prop:obstruction}.)
\item \textbf{CF-ISM ($d_0=320$ at $t=0$ only, pure double precision for
  $t>0$):} Across \emph{all 16 combinations}---including the mixed even-and-odd
  harmonic cases Ex.~5.2(c),(d) and the off-grid shift---CF-ISM recovers the
  unknown shift centre automatically and maintains $D\equiv0$ for the
  entire simulation. The same holds for the shifted Crank--Nicolson benchmark
  of Section~\ref{sec:num1d} ($\kappa=\gamma=1$, $N=512$, $T=100$): standard
  CN breaks at $T_\ast\approx 36.1$--$36.6$ for $\theta\in\{\pi/8,\pi/4,0.3712\}$,
  while CF-ISM initializes in $0.23\text{--}1.14\,\mathrm s$ and keeps
  $D\equiv0$ up to $T=100$.
\end{itemize}

\begin{table}[t]
\centering
\caption{Shifted benchmarks of~\cite[Ex.~5.2]{Wang2023} ($N=512$) with
$\theta\in\{0,\pi/8,\pi/4,0.3712\}$. Reported are first recorded $T_\ast$
($|\mathrm{mean}|>10^{-2}$, or ``---'' with $\max|\mathrm{mean}|$ if below
$10^{-2}$) for ESP ($S=1$, $N_s=24$, $N_{\rm Iter}=20$, $T_e=1,2,2,10$ for
(a)--(d)) and standard ETDRK4 ($S=0$, $\tau=0.02,0.02,0.02,0.05$), and maximum
internal transverse-channel defect $D$ for $\mathcal T_2$ isotypic ($S=0$,
ETDRK4) and CF-ISM ($S=0$, ETDRK4, $d_0=320$ at $t=0$, pure double for
$t>0$).}\label{tab:wang52}
\small
\begin{tabular}{llcccc}
\toprule
case & shift $\theta$ & ESP ($S=1$) & Std ETDRK4 ($S=0$) & $\mathcal T_2$ iso & CF-ISM\\
\midrule
5.2(a) $\kappa=1$, $T=200$ & $0$      & --- ($1.0\times10^{-18}$) & $T_\ast=42.6$ & $0$ & $0$\\
$\sin\tilde x$             & $\pi/8$  & --- ($6.3\times10^{-18}$) & $T_\ast=39.0$ & $0$ & $0$\\
                           & $\pi/4$  & $\mathbf{T_\ast=41}$      & $T_\ast=37.4$ & $0$ & $0$\\
                           & $0.3712$ & $\mathbf{T_\ast=38}$      & $T_\ast=37.4$ & $0$ & $0$\\
\midrule
5.2(b) $\kappa=0.9$, $T=200$ & $0$      & --- ($6.5\times10^{-19}$) & $T_\ast=65.0$ & $0$ & $0$\\
$\sin\tilde x$               & $\pi/8$  & --- ($2.3\times10^{-18}$) & $T_\ast=58.0$ & $0$ & $0$\\
                             & $\pi/4$  & --- ($2.3\times10^{-17}$) & $T_\ast=54.4$ & $0$ & $0$\\
                             & $0.3712$ & $\mathbf{T_\ast=56}$      & $T_\ast=54.4$ & $0$ & $0$\\
\midrule
5.2(c) $\kappa=0.1$, $T=200$ & $0$      & --- ($6.7\times10^{-11}$) & --- & n/a & $0$\\
$\tfrac12(\sin\tilde x+\sin4\tilde x)$ & $\pi/8$  & --- ($3.5\times10^{-11}$) & --- & n/a & $0$\\
                             & $\pi/4$  & --- ($1.8\times10^{-11}$) & --- & n/a & $0$\\
                             & $0.3712$ & --- ($9.9\times10^{-11}$) & --- & n/a & $0$\\
\midrule
5.2(d) $\kappa=\sqrt{0.001}$, $T=600$ & $0$      & --- ($2.3\times10^{-13}$) & --- & n/a & $0$\\
$\tfrac12(\sin\tilde x+\sin8\tilde x)$ & $\pi/8$  & --- ($8.7\times10^{-14}$) & --- & n/a & $0$\\
                             & $\pi/4$  & --- ($9.4\times10^{-14}$) & --- & n/a & $0$\\
                             & $0.3712$ & --- ($6.0\times10^{-9}$)  & --- & n/a & $0$\\
\bottomrule
\end{tabular}
\end{table}

\begin{figure}[t]
\centering
\includegraphics[width=\textwidth]{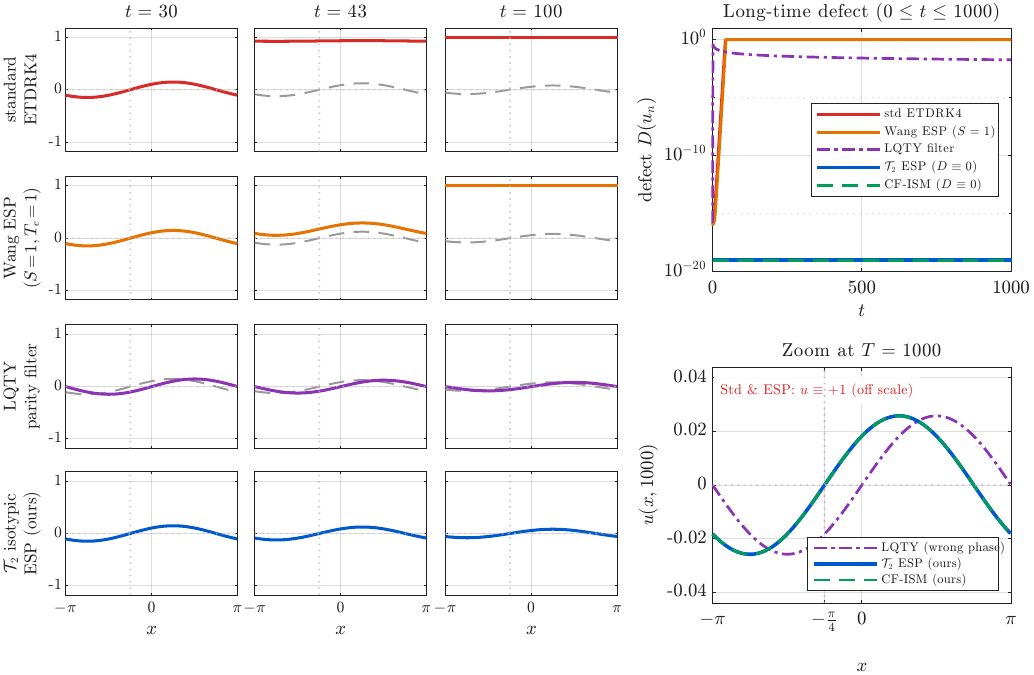}
\caption{Long-time shifted evolution on~\cite[Ex.~5.2(a)]{Wang2023}
($u_0(x)=\sin(x+\pi/4)$, $\kappa=\gamma=1$, $N=512$, $0\le t\le 1000$).
Left: solution profiles at $t=30,43,100$ (dashed grey: exact shifted odd
profile; vertical dotted line: true symmetry centre $x_0=-\pi/4$). Standard
ETDRK4 ($S=0$) and Wang's ESP ($S=1$, $T_e=1$, $N_s=24$, $N_{\rm Iter}=20$)
both break symmetry and collapse to $u\equiv +1$; the origin-centred parity
filter~\cite{LQTY2022} forces $u(0,t)=0$ and introduces an $O(1)$ phase error
from the first step; the proposed $\mathcal T_2$ isotypic implementation of
the \emph{exact same} $S=1$ ESP integrator (pure double precision, no shift
detection) preserves anti-periodicity identically ($D_{\mathrm{anti}}\equiv 0$,
with shifted reflection defect $\max D_{\mathrm{shift}}\approx 6.4\times10^{-16}$).
Top right: defect histories up to $T=1000$ ($D_{\mathrm{shift}}$ for standard
schemes/filter, $D_{\mathrm{anti}}$ for $\mathcal T_2$ ESP, and canonical-grid
$D$ for CF-ISM ETDRK4).
Bottom right: zoom at $T=1000$, where $\mathcal T_2$ ESP and CF-ISM ETDRK4
agree to $3.4\times10^{-12}$ (reflecting temporal discretization difference
between ESP and ETDRK4) while standard ETDRK4 and Wang's ESP are off-scale at
$u\equiv+1$.}\label{fig:wang52_evol}
\end{figure}

Figure~\ref{fig:wang52_evol} illustrates this failure mode visually on
Wang's~\cite[Ex.~5.2(a)]{Wang2023} benchmark with $\theta=\pi/4$ up to
$T=1000$, comparing all major existing approaches side by side. Standard
ETDRK4 ($S=0$) has largely departed the shifted odd branch by $t=43$
($D=0.93$) and collapses to $u\equiv+1$ well before $t=100$; Wang's ESP
($S=1$, $T_e=1$, $N_s=24$, $N_{\rm Iter}=20$) visibly lifts off the shifted
odd wave at $t=43$ ($\mathrm{mean}=0.174$, $D=0.175$) and collapses to
$u\equiv+1$ before $t=100$; and the origin-centred parity filter
of~\cite{LQTY2022} corrupts the shifted phase at the very first step by
forcing $u(0,t)=0$ instead of $u(-\pi/4,t)=0$. By contrast, running the
\emph{exact same} $S=1$ ESP integrator in the $\mathcal T_2$ isotypic
representation (pure IEEE double precision, zero shift detection, no
filtering) keeps $D_{\mathrm{anti}}\equiv0$ up to $T=1000$ (with
$\max D_{\mathrm{shift}}\approx 6.4\times 10^{-16}$), agreeing with CF-ISM
ETDRK4 to $3.4\times10^{-12}$ at $T=1000$.

\subsection{CF-ISM horizon law and faithfulness under off-grid shifts}
\begin{figure}[t]
\centering
\includegraphics[width=\textwidth]{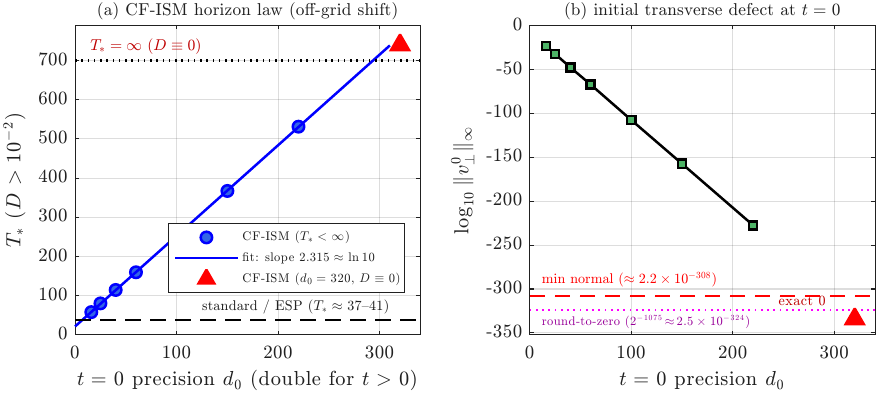}
\caption{CF-ISM horizon law for the mixed-harmonic off-grid shifted datum
$\tfrac12(\sin(x+\theta)+\sin(4(x+\theta)))$, $\theta=0.37123456789$,
$\kappa=1$, $S=0$, $N=128$, with pure IEEE double-precision ETDRK4 for all
$t>0$. (a)~Breakdown time $T_\ast$ versus $t=0$ precision $d_0$; (b)~initial
transverse defect $\|v^0_\perp\|_\infty$ crossing the minimum positive normal
($2.2\times10^{-308}$) and round-to-zero threshold
($2^{-1075}\approx 2.5\times10^{-324}$) at $d_0=320$.}\label{fig:cfism_horizon}
\end{figure}

To verify Theorem~\ref{thm:cfism} quantitatively, we take the mixed-harmonic
off-grid shifted datum $u_0(x)=\tfrac12(\sin(x+\theta)+\sin(4(x+\theta)))$
with $\theta=0.37123456789$, $\kappa=\gamma=1$, $S=0$, $N=128$, and vary the
$t=0$ initialization precision $d_0\in\{16,25,40,60,100,150,220,320\}$ while
always running ETDRK4 ($\tau=0.05$) in \emph{pure hardware IEEE double
precision} for $t>0$. Figure~\ref{fig:cfism_horizon} confirms
\eqref{eq:cfhorizon}: for $d_0\le 220$ the breakdown time grows linearly with
$d_0$ from $T_\ast=57.6$ ($d_0=16$) to $T_\ast=531.2$ ($d_0=220$), with fitted
slope $2.315\approx\ln 10=2.303$; once $d_0=320$, $\|v^0_\perp\|_\infty$ drops
below the IEEE-754 double round-to-zero threshold
$2^{-1075}\approx 2.47\times 10^{-324}$ and rounds naturally to exact
\texttt{0.0}, yielding $D\equiv0$ ($T_\ast=\infty$). Moreover, adding a
genuine constant asymmetry $\delta$ to $\sin(x+0.37123456789)$ gives initial
channel value $\|v^0_\perp\|_\infty=\delta$ to double precision and yields
$T_\ast=21.85$ ($\delta=10^{-10}$) and $T_\ast=69.00$ ($\delta=10^{-30}$)
under $S=0$, $D>10^{-2}$ (confirming the logarithmic scaling
$\Delta T_\ast\approx 20\ln 10\approx 46.1$ of Corollary~\ref{cor:faithful});
likewise, on the mixed-harmonic datum
$\tfrac12(\sin(x+\theta)+\sin(4(x+\theta)))+\delta$ ($N=128$, $d_0=100$,
\texttt{test\_shifted.m}), the mean-centred residual recovers $\delta$ to
exact double precision ($\max|v^0_{\mathrm{even}}-\delta|/\delta=0$) for both
$\delta=10^{-8}$ and $\delta=10^{-30}$. Furthermore, because
\texttt{cf\_gauge\_init.m} normalizes residuals by $\|u_0\|_\infty$ and uses
relative rank checks, the same mixed-harmonic waveform scaled to tiny
amplitude $A=10^{-18}$ still achieves exact double underflow
($\max|v^0_{\mathrm{even}}|=0$ at $d_0=320$), double-cancelled inputs such as
$(1+10^{-20}\sin(x+0.2)+3\times 10^{-21})-1$ and constant-dominated inputs
such as $10^{-32}\sin(x+0.02)+10^{-8}$ automatically fall back to $d_0$-digit
coarse sampling (recovering both the centre and the even channel without
discarding nonzero VPA data), and degenerate 2D inputs such as
$0.2+\sin x+0.1\cos 2x$ (independent of $y$) freeze the unidentifiable
coordinate cleanly without \texttt{NaN}s.

\subsection{Off-grid shifted 2D and 3D patterns}
\begin{figure}[t]
\centering
\includegraphics[width=\textwidth]{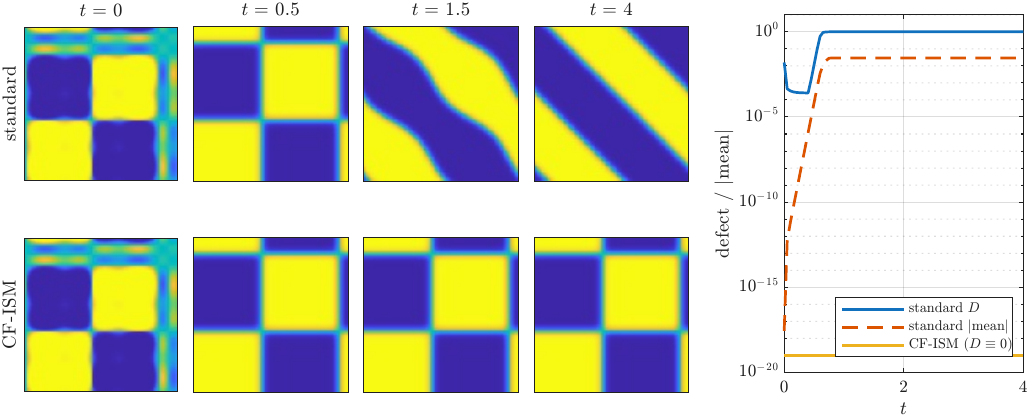}
\caption{Off-grid shifted 2D mixed-harmonic X-junction ($\eps=0.1$,
$N=128$, ETDRK2, $\tau=10^{-3}$): unshifted-grid standard Fourier (top;
initial off-grid spatial-aliasing defect $D(0)=1.59\times10^{-2}$ growing to
$D=1.00$) versus CF-ISM (bottom, canonical grid eliminating the initial
aliasing defect to $D\equiv0$, plotted back on the original grid) and defect
histories (right). In a same-grid control started from the canonical-grid
datum at $N=128$, standard Fourier also keeps $D=0$ through $T=4$ (cf.\
Table~\ref{tab:fftN}).}\label{fig:cfism_2d}
\end{figure}

\begin{figure}[t]
\centering
\includegraphics[width=\textwidth]{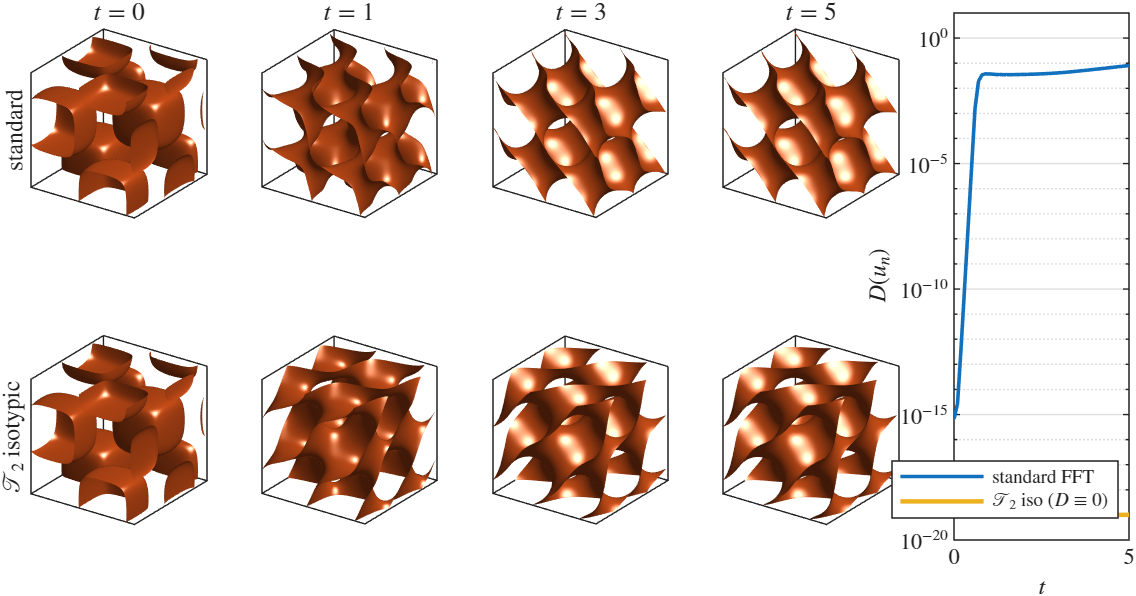}
\caption{Off-grid shifted 3D octant saddle ($\eps=0.1$, $N=96$, ETDRK2,
$\tau=2\times10^{-3}$): level set $u=0.5$ of standard Fourier (top) versus
$\mathcal T_2$ isotypic (bottom, pure double precision with no shift
detection), and anti-periodicity defect histories (right).}\label{fig:O3d_shifted}
\end{figure}

Finally we demonstrate both shifted mechanisms in two and three dimensions
with off-grid shifts $(\theta_1,\theta_2,\theta_3)=(0.37123456,0.71828182,1.13141592)$:
\begin{enumerate}
\item \textbf{2D off-grid shifted mixed-harmonic X-junction (CF-ISM).}
  Take
  \[
    u_0(x,y)=\tanh\!\Big(\frac{(\sin\tilde x+0.4\sin 4\tilde x)(\sin\tilde y+0.4\sin 4\tilde y)}{\sqrt2\eps}\Big),
    \qquad \tilde x=x+\theta_1,\quad \tilde y=y+\theta_2,
  \]
  with $\eps=0.1$, $N=128$, ETDRK2 ($\tau=10^{-3}$, $T=4$). Because of the $k=4$ harmonic, $\mathcal T_2$
  anti-periodicity does not hold; only shifted reflection symmetry about
  $(-\theta_1,-\theta_2)$ holds. On the unshifted collocation grid
  $\mathcal G_{128}$, spatial aliasing of the sharp $\tanh$ interface already
  injects an initial symmetry defect $D_{\mathrm{std}}(0)=1.59\times10^{-2}$
  (measured via Fourier back-shift onto the aligned frame; this aliasing
  defect is independent of double-precision rounding and decreases to
  $3.6\times10^{-4}$ at $N=256$ and $1.1\times10^{-7}$ at $N=512$). Starting
  from that unaligned datum, standard Fourier grows from $D(0)=1.59\times10^{-2}$
  to complete collapse ($D=1.00$, $|\mathrm{mean}|=2.9\times10^{-2}$), merging
  the four domains into a distorted diagonal stripe. By contrast, CF-ISM
  ($d_0=320$ at $t=0$, costing $6.7\,\mathrm s$ in VPA; pure double for $t>0$)
  aligns the collocation grid to the reflection centre, eliminating the
  initial spatial-aliasing defect ($D(0)=0$) and preserving the shifted
  checkerboard with canonical-grid $D\equiv0$ throughout
  (Figure~\ref{fig:cfism_2d}). (In a same-grid control where standard Fourier
  is started from the \emph{same} canonical-grid datum at $N=128$, standard
  Fourier also stays at $D=0$ through $T=4$ because $N=128$ is a zero-defect
  FFTW length in Table~\ref{tab:fftN}; Figure~\ref{fig:cfism_2d} therefore
  illustrates elimination of off-grid initial discretization/aliasing defects
  by canonical-frame alignment, whereas pure time-stepping FFT-rounding
  amplification on aligned data is isolated in Table~\ref{tab:law},
  Section~\ref{sec:num23d}, and Figure~\ref{fig:O3d_shifted}.) Similarly, for
  the off-grid shifted 2D point-odd ($C_2$) datum
  $\tanh(\sin(\tilde x+2.5\sin 3\tilde y)/(\sqrt2\eps))$
  ($N=64$, $\eps=0.05$, $\tau=2\times10^{-3}$, 200 ETDRK2 steps in
  \texttt{test\_shifted.m}), CF-ISM with $d_0=320$ locates a valid $C_2$
  centre automatically (to $0.0\pmod{\pi}$ in $x$ and $0.0\pmod{\pi/3}$ in
  $y$) and maintains $D\equiv0$.
\item \textbf{3D off-grid shifted octant saddle ($\mathcal T_2$).}
  Take
  \[
    u_0(x,y,z)=\tanh\!\Big(\frac{\sin\tilde x\sin\tilde y\sin\tilde z}{\sqrt2\eps}\Big),
    \qquad \tilde x=x+\theta_1,\quad \tilde y=y+\theta_2,\quad \tilde z=z+\theta_3,
  \]
  with $\eps=0.1$, $N=96$, ETDRK2 ($\tau=2\times10^{-3}$, $T=5$). Because $N=96$
  uses a mixed-radix FFT, standard Fourier first crosses $D>10^{-2}$ at the
  recorded time $T_\ast=0.70$ ($D(0.70)=1.90\times 10^{-2}$, reaching a peak
  defect $\max_{t\le 5}D=8.09\times10^{-2}$), whereas the $\mathcal T_2$
  isotypic scheme---running in 100\% pure IEEE double precision with \emph{no}
  shift detection and \emph{no} high-precision initialization---preserves all
  seven vanishing wavenumber-parity channels identically ($D\equiv0$, spatial
  mean $\equiv 0$) throughout, preserving the $\mathcal T_2$-anti-periodic
  triply periodic interface geometry while strictly forbidding any symmetry
  breaking into even-wavenumber or constant modes
  (Figure~\ref{fig:O3d_shifted}). On the off-grid shifted 2D X-junction
  ($N=256$, $\eps=0.05$, $T=10$), the fixed-origin reflection isotypic scheme
  develops a small sub-threshold anti-periodicity defect
  $\max D=6.7\times10^{-4}$ (below the $10^{-2}$ macroscopic threshold),
  whereas $\mathcal T_2$ again keeps $D\equiv0$. (At $N=256=2^8$, the first
  radix-2 butterfly of the standard FFT happens to separate even and odd
  wavenumbers, another instance of the FFT-length fragility noted in
  Table~\ref{tab:fftN}.)
\end{enumerate}

\subsection{Drop-in extensibility: shifted high-order schemes, CAC and Cahn--Hilliard}\label{sec:ext}
A practical strength of the isotypic spectral representation is that it is
\emph{not} tied to a particular time-stepping formula or PDE: it acts as a
modular drop-in backend that can be placed underneath existing or future
high-order, energy-stable or adaptive algorithms so that accurate local time
integration is not corrupted by long-time symmetry drift. Upgrading a
Fourier pseudo-spectral solver requires three representation-level
substitutions, leaving the time-stepping formula untouched:
\begin{enumerate}
\item replace the full-grid \texttt{fft}/\texttt{ifft} by the channel-wise
  transforms of $A_{\mathrm R}$ or $\mathcal T_2$;
\item replace pointwise evaluations of $u^3$ (or general odd $F(u)$) by the
  $\Z_2^d$-graded channel product (or Walsh--Hadamard orbit evaluation);
\item evaluate any scalar inner products $\langle a,b\rangle$ or scalar
  auxiliary variables using the exact orthogonal channel weights (so all
  scalars multiply channels diagonally), and use symmetry-preserving stopping
  norms for iterative or adaptive controllers.
\end{enumerate}
Beyond the twelve unshifted 1D integrators tested in Section~\ref{sec:num1d}
(the five schemes of Table~\ref{tab:law} plus the seven high-order schemes of
Section~\ref{sec:highorder}; adding SAV1 below gives thirteen integrators in
total), we verify this drop-in extensibility on shifted high-order integrators
as well as two further phase-field PDEs (\texttt{exp\_extensibility.m}):
\begin{itemize}
\item \textbf{Shifted high-order integrators (pure double $\mathcal T_2$):}
  On the off-grid shifted datum $u_0(x)=\sin(x+0.37123456789)$
  ($\kappa=\gamma=1$, $S=0$, $N=512$, $T=200$, anti-periodicity defect
  monitored every $0.1$ time units), standard FFT implementations of
  BDF2-SAV ($\tau=0.01$), SBDF4 ($\tau=0.01$), IFRK4 ($\tau=0.02$) and
  SDC4 ($\tau=0.05$) all first cross $D>10^{-2}$ at the recorded time
  $T_\ast=37.30$ and collapse to $u\equiv\pm1$ (matching ETDRK4 in
  Table~\ref{tab:wang52}); standard Strang2 creates a nonzero defect
  $5.6\times10^{-17}$. Plugging in the $\mathcal T_2$ isotypic backend (pure
  IEEE double precision, parity-consistent initialization) keeps $D\equiv0$
  for all five schemes. (First-order SAV1 likewise has non-symmetric
  equivalence $2.2\times10^{-16}$ in $u$ and $4.4\times10^{-16}$ in $r$,
  breaks at $T_\ast=44.00$ on $u_0=\sin x$, and keeps $D\equiv0$ with the
  isotypic backend.)
\item \textbf{Conservative Allen--Cahn (CAC):} The nonlocal Lagrange
  multiplier $\partial_tu=\Delta u-\eps^{-2}(u^3-u)+\eps^{-2}\overline{u^3-u}$
  enforces spectral zero-mode conservation ($\widehat u(\mathbf 0)\equiv0$,
  with physical-grid sample mean $|\mathrm{mean}(u)|\le1.5\times10^{-17}$ due
  to summation rounding). On the 2D odd--odd X-junction ($N=256$, $\eps=0.05$,
  ETDRK2, $\tau=10^{-3}$, monitored every $0.02$ time units), standard FFT CAC
  still breaks the odd--odd symmetry completely at $T_\ast=0.220$
  ($\max D=1.00$), demonstrating that exact mass conservation cannot prevent
  multi-dimensional symmetry collapse. Plugging in the isotypic backend keeps
  $D\equiv0$ throughout.
\item \textbf{Fourth-order Cahn--Hilliard equation:} For
  $\partial_tu=-\Delta(\eps^2\Delta u+u-u^3)$ with stabilized ETDRK2
  ($L_k=\eps^2|k|^4+S|k|^2$, $N_k(u)=|k|^2\mathcal F[(1+S)u-u^3]$, $S=2$), the
  isotypic and standard implementations agree to $5.3\times10^{-16}$ on
  non-symmetric 2D data ($N=32$, $\tau=5\times10^{-4}$, 20 steps). On the 2D
  odd--odd X-junction ($N=192$, $\eps=0.1$, $\tau=5\times10^{-4}$, monitored
  every 25 steps), standard FFT Cahn--Hilliard amplifies the symmetry defect
  by more than two orders of magnitude from $5.0\times10^{-16}$ at step~1 to
  $1.86\times10^{-13}$ by $T=1.0$, whereas the drop-in isotypic
  Cahn--Hilliard scheme maintains $D\equiv0$.
\end{itemize}
Ready-to-use drop-in backend modules, worked examples (native MATLAB
\texttt{ode15s}, Python \texttt{scipy.\allowbreak integrate.\allowbreak solve\_ivp}, standalone
Chebfun \texttt{spin}/\texttt{spin2}-style and Shenfun-style ETDRK drivers,
and a Dedalus~v3 channel-equation recipe), and side-by-side accuracy and
wall-clock timing comparisons are provided at
\url{https://github.com/minghuili123456789/isotypic-phase-field/tree/main/dropin-library-backends}.

\subsection{Summary comparison with existing approaches}
Table~\ref{tab:compare} summarizes the capabilities of existing methods and
the proposed isotypic framework across all tests in
Sections~\ref{sec:num1d}--\ref{sec:numshift}. Standard schemes and ESP break
at $T_\ast=\lambda_\perp^{-1}\log(1/\umach)+O(1)$ (or rely on fragile
sub-resolution absorption); origin-centred parity filters~\cite{LQTY2022}
fail on shifted data and destroy genuine small asymmetries; a~posteriori
zero-mode suppression fails in 2D and 3D whenever the unstable transverse
mode has zero mean. Only the isotypic framework ($A_{\mathrm R}$,
$\mathcal T_2$, CF-ISM) satisfies every requirement simultaneously without
artificial intervention.

\begin{table}[t]
\centering
\caption{Comparison of existing approaches and the proposed non-interventional
isotypic framework across all 1D/2D/3D benchmarks
(Sections~\ref{sec:num1d}--\ref{sec:numshift}). All methods run all $t>0$
time steps in pure hardware IEEE double precision. For CF-ISM, exactness
($D\equiv0$) requires the pre-rounding initial transverse bound to satisfy
$C_0 10^{-d_0}<2^{-1075}$ (Theorem~\ref{thm:cfism}(ii), attained at $d_0=320$
for the tested data); \texttt{cf\_gauge\_init.m} is implemented and tested in
1D and 2D, while 3D shifted benchmarks use $\mathcal T_2$ (3D CF-ISM follows
the same theory under $D^2\widetilde{\mathcal J}_{\mathrm c}\succ0$).}\label{tab:compare}
\scriptsize
\setlength{\tabcolsep}{3.5pt}
\resizebox{\textwidth}{!}{%
\begin{tabular}{lccccc}
\toprule
method & non-interv. & 1D origin & 1D shifted & 2D/3D multi-mode & faithful ($\delta>0$)\\
\midrule
Standard Fourier (IMEX/ETD/CN)   & yes & fails ($T_\ast$) & fails ($T_\ast$) & fails (FFT-$N$ dep.) & $\delta\gtrsim\umach$ only\\
Wang~\cite{Wang2023} ESP ($S=1$) & yes & fails ($T_e=20$) & fails ($\pi/4$, off-grid) & fails & $\delta\gtrsim\umach$ only\\
LQTY~\cite{LQTY2022} parity filter & no & exact & fails ($O(1)$ phase) & origin-only & fails (erases $\delta$)\\
Zero-mode suppression ($\widehat u(\mathbf 0)\!\gets\!0$) & no & $D\sim\umach$ & $D\sim\umach$ ($\kappa=1$) & fails ($t=0.3$--$0.6$) & fails (erases mean)\\
\midrule
Isotypic $A_{\mathrm R}$ (ours)  & yes & exact ($D\equiv0$) & grid/half-grid & exact ($D\equiv0$) & rel.\ acc.\ (poly.)\\
Isotypic $\mathcal T_2$ (ours)   & yes & exact (odd harm.)  & exact (any shift) & exact (2D/3D anti-per.) & rel.\ acc.\ (poly.)\\
CF-ISM (ours, $d_0=320$ at $t=0$)& yes & exact ($D\equiv0$) & exact (if $v^0_\perp\!=\!0$) & exact in 2D (if $v^0_\perp\!=\!0$) & rel.\ acc.\ (poly.)\\
\bottomrule
\end{tabular}}
\end{table}

%======================================================================
\clearpage
\section{Conclusions}\label{sec:concl}

We have shown that the long-time loss of symmetry in phase-field
simulations is not caused by temporal or spatial truncation errors---every
scheme in $\mathfrak I$ is exactly equivariant in exact arithmetic---but by
non-equivariant rounding errors amplified by transverse instability of the
symmetric state. Under the discounted-seed and instability hypotheses of
Theorem~\ref{thm:breakdown}, breakdown occurs on the logarithmic time scale
$T_\ast=\lambda_\perp^{-1}\log(1/\umach)+O(1)$, which we observed across
thirteen integrators: IMEX1, ETDRK2, ETDRK4, Crank--Nicolson, ESP, SAV1,
BDF2-SAV, ARS(4,4,3), SBDF3, SBDF4, IFRK4, SDC4, and Strang2 splitting; cases
where standard schemes appear to preserve symmetry are explained by
sub-resolution absorption or failure of the seed condition and are fragile
with respect to stabilization, step size, shift, and FFT length (as shown
directly on both 1D and 2D ESP benchmarks of~\cite{Wang2023}).

The isotypic spectral framework removes the cause rather than the symptom.
Both the reflection ($A_{\mathrm R}$, Type~I and Type~II) and half-period
translation ($\mathcal T_2$) realizations are algebraically identical to the
Fourier pseudo-spectral method on the same collocation grid and preserve every
target subgroup subspace identically ($D\equiv 0$, up to IEEE signed zeros,
assuming zero-transverse initialization and finite intermediates) without
projections, filters, or zero-mode resets. For the polynomial channel
algebra~\eqref{eq:chalg}, transverse rounding errors scale relative to the
transverse component itself, propagating genuine asymmetries faithfully as
verified against a 120-digit symbolic reference down to $\delta=10^{-60}$ and
observed down to $\delta=10^{-200}$ (while Walsh--Hadamard orbit evaluation
guarantees exact zero-channel invariance for general odd potentials). For
unknown off-grid shifts, $\mathcal T_2$ preserves anti-periodicity exactly in
pure double precision from parity-consistent initialization, while CF-ISM
confines $d_0$-digit gauge initialization to $t=0$ only and runs 100\% of the
$t>0$ time stepping in hardware double precision, extending the horizon as
$\lambda_\perp^{-1}d_0\ln 10+O(1)$ and attaining exact zero transverse
channels once $C_0 10^{-d_0}<2^{-1075}$. When the data lie in a symmetric
subspace (or via exact-zero-channel short-circuiting), reduced evaluation
matches or beats standard FFT speed while giving identically equal results.

From a software and algorithmic perspective, the isotypic framework serves as
a \emph{modular drop-in spectral layer} for PDE solvers and
scientific-computing workflows (demonstrated for MATLAB \texttt{ode15s},
Python/SciPy, and standalone Chebfun/Shenfun-style drivers, alongside a
Dedalus~v3 channel recipe) rather than a standalone time stepper. Researchers
can retain their preferred high-order, energy-stable, maximum-bound-preserving
or adaptive integrators by replacing the transform, nonlinear evaluation, and
(where needed) orthogonal inner products and symmetry-preserving stopping
norms. The same half-domain channel splitting extends directly to homogeneous
Neumann and Dirichlet boundary conditions on $[0,\pi]^d$. Natural further
extensions left for future work include full algorithmic verification of
non-Abelian dihedral/cubic groups ($D_4/B_3$), symmetric finite-difference
and finite-element discretizations, coupled multi-phase gradient flows, and
reliable numerical bifurcation tracking enabled by sub-roundoff faithfulness.

\section*{Data and code availability}
{\sloppy
All code, reproducible experiment scripts, supplementary boundary-condition
evidence, and drop-in library adapters are openly available on GitHub under the
MIT license:
\begin{itemize}
\item \textbf{Main paper repository (core library and all manuscript
  experiments, organized by section):}\\
  \url{https://github.com/minghuili123456789/isotypic-phase-field}
\item \textbf{Supplementary boundary-condition evidence (1D/2D homogeneous
  Neumann and 1D homogeneous Dirichlet scripts and comparison figures):}\\
  \url{https://github.com/minghuili123456789/isotypic-phase-field/tree/main/supplementary-bc-evidence}
\item \textbf{Drop-in library backends (minimal adapters and worked examples
  for MATLAB \texttt{ode15s}, Chebfun \texttt{spin}/\texttt{spin2}/\texttt{spin3},
  Python NumPy/SciPy, Dedalus~v3, Shenfun, and custom solvers):}\\
  \url{https://github.com/minghuili123456789/isotypic-phase-field/tree/main/dropin-library-backends}
\end{itemize}
}

\appendix
\section{Proof of Theorem~\ref{thm:breakdown}}\label{app:proofs}
Write $\Psi_{n,m}=J^\perp_{n-1}\cdots J^\perp_m$ ($\Psi_{m,m}=I$) and
$\mu:=\ee^{\lambda_\perp\tau}>1$. Iterating \eqref{eq:transverse} from $e_0=0$
gives the exact discrete variation-of-constants representation
\begin{equation}\label{eq:voc}
 e_n=\sum_{m=0}^{n-1}\Psi_{n,m+1}\,\big(r_{\perp,m}+q_m(w_m,e_m)\big).
\end{equation}
Recall that for $A\cong\Z_2^d$ the orthogonal projections $\Pi_{\bsig}$ satisfy
$\|\Pi_{\bsig}\|_{\infty\to\infty}\le 1$ and $e_n=\sum_{\bsig\notin\Sigma}(u_n)_{\bsig}$,
so the defect $D(u_n)=\max_{\bsig\notin\Sigma}\|(u_n)_{\bsig}\|_\infty$ and
$\|e_n\|_\infty$ are equivalent via
\begin{equation}\label{eq:Dnorm}
 2^{-d}\|e_n\|_\infty \le D(u_n) \le \|e_n\|_\infty .
\end{equation}

\paragraph{Part (i) (exponentially weighted bootstrap)}
Set $\gamma_\tau:=1-\mu^{-1}=1-\ee^{-\lambda_\perp\tau}$,
$K_1:=2KC/\gamma_\tau$, and define the weighted maximum
$\beta_n:=\max_{0\le m\le n}\mu^{-m}\|e_m\|$. Let
$T_1:=\lambda_\perp^{-1}\log(\eta/(K_1\umach))$ so that for $k\tau\le T_1$
we have $K_1\umach\mu^k\le\eta$. Assume $\eta$ is small enough that
$KC_q\eta/\gamma_\tau\le\tfrac14$, and $\umach$ is small enough that
$KC_q C_w (T_1/\tau)^2\umach/(2\mu)\le\tfrac14$. We prove by induction on $n$
that $\beta_n\le K_1\umach$ (and hence $\|e_k\|\le K_1\umach\mu^k\le\eta$) for
all $n\tau\le T_1$. The base case $n=0$ holds since $e_0=0$. Suppose
$\beta_{n-1}\le K_1\umach$ and $n\tau\le T_1$. For every $1\le k\le n$,
multiplying \eqref{eq:voc} by $\mu^{-k}$ and using
$\|\Psi_{k,m+1}\|\le K\mu^{k-m-1}$ and
$\|q_m(w_m,e_m)\|\le C_q(C_w m\umach+\|e_m\|)\|e_m\|$ gives
\begin{align*}
 \mu^{-k}\|e_k\|
 &\le K\sum_{m=0}^{k-1}\mu^{-m-1}\Big(C\umach
      + C_q C_w m\umach\,\beta_{n-1}\mu^m
      + C_q\,\beta_{n-1}^2\mu^{2m}\Big) \\
 &\le \frac{KC\umach}{\mu-1}
      + \frac{KC_q C_w\umach k^2}{2\mu}\,\beta_{n-1}
      + \frac{KC_q(\beta_{n-1}\mu^{k-1})}{\gamma_\tau}\,\beta_{n-1}.
\end{align*}
Since $(\mu-1)^{-1}=\mu^{-1}/\gamma_\tau<\gamma_\tau^{-1}$,
$\beta_{n-1}\mu^{k-1}\le K_1\umach\mu^k\le\eta$, and our smallness hypotheses
bound the last two coefficients by $\tfrac14$ each, we obtain
\[
 \mu^{-k}\|e_k\| \le \frac{KC\umach}{\gamma_\tau} + \Big(\frac14+\frac14\Big)\beta_{n-1}
 \le \frac12 K_1\umach + \frac12 K_1\umach = K_1\umach .
\]
Hence $\beta_n\le K_1\umach$, completing the induction. For $n\tau\le T_1$,
\eqref{eq:Dnorm} gives $D(u_n)\le\|e_n\|_\infty\le\eta$.

\paragraph{Part (ii) (discounted-seed lower bound)}
In the steady case $J_n^\perp\equiv J^\perp$, $P_\mu J^\perp=\mu P_\mu$.
Applying $P_\mu$ to \eqref{eq:voc} gives
\[
 P_\mu e_n = \mu^n\sum_{m=0}^{n-1}\mu^{-m-1}P_\mu r_{\perp,m}
           + \sum_{m=0}^{n-1}\mu^{n-m-1}P_\mu q_m(w_m,e_m).
\]
Let $n_\ast$ be the smallest integer $n\ge n_0$ such that
$c_0\umach\mu^n\ge 2^{d+1}\|P_\mu\|_\infty\eta$, so that
\[
 n_\ast\tau=\lambda_\perp^{-1}\log\!\big(2^{d+1}\|P_\mu\|_\infty\eta/(c_0\umach)\big)+O(1)
 \quad\text{and}\quad
 c_0\umach\mu^{n_\ast}<2^{d+1}\|P_\mu\|_\infty\mu\eta.
\]
Suppose for contradiction that $D(u_m)<\eta$ (hence $\|e_m\|_\infty<2^d\eta$
by \eqref{eq:Dnorm}) for all $m\le n_\ast$. Using $\|e_m\|_\infty<2^d\eta$
directly in the quadratic term gives
$\|q_m(w_m,e_m)\|_\infty\le C_q(C_w m\umach+2^d\eta)\|e_m\|_\infty$ for every
$m\le n_\ast$. Repeating the weighted induction of Part~(i) on
$\beta_k=\max_{0\le m\le k}\mu^{-m}\|e_m\|_\infty$ over $k\le n_\ast$
(without restricting $n_\ast\tau$ to $T_1$) yields, for $1\le k\le n_\ast$,
\[
 \mu^{-k}\|e_k\|_\infty
 \le \frac{KC\umach}{\gamma_\tau}
   + \Big(\frac{KC_q C_w n_\ast^2\umach}{2\mu}
        + \frac{2^d K C_q\eta}{\gamma_\tau}\Big)\beta_{k-1}
 \le \frac12 K_1\umach + \frac12\beta_{k-1},
\]
provided $2^d K C_q\eta/\gamma_\tau\le\tfrac14$ and
$K C_q C_w n_\ast^2\umach/(2\mu)\le\tfrac14$. Hence $\beta_{n_\ast}\le K_1\umach$,
so $\|e_m\|_\infty\le K_1\umach\mu^m$ for all $m\le n_\ast$, and in particular
$K_1\umach\mu^{n_\ast}\le M_0\eta$ with the $(\eta,\umach)$-independent
constant $M_0:=2^{d+1}K_1\|P_\mu\|_\infty\mu/c_0$. Using
$\|q_m(w_m,e_m)\|_\infty\le C_q(C_w m\umach+\|e_m\|_\infty)\|e_m\|_\infty$, the
projected remainder at $n=n_\ast$ is bounded by
\begin{align*}
 \Big\|\sum_{m=0}^{n_\ast-1}\mu^{n_\ast-m-1}P_\mu q_m(w_m,e_m)\Big\|_\infty
 &\le \|P_\mu\|_\infty C_q\Big(\frac{C_w K_1 n_\ast^2\umach}{2\mu}
      + \frac{K_1 M_0\eta}{\mu-1}\Big)\umach\mu^{n_\ast}.
\end{align*}
Choosing $\eta$ and $\umach$ small enough that
\[
 \frac{\|P_\mu\|_\infty C_q K_1 M_0\eta}{\mu-1}\le\frac{c_0}{4}
 \qquad\text{and}\qquad
 \frac{\|P_\mu\|_\infty C_q C_w K_1 n_\ast^2\umach}{2\mu}\le\frac{c_0}{4},
\]
the discounted-seed hypothesis \eqref{eq:seedlb} yields
\[
 \|P_\mu e_{n_\ast}\|_\infty \ge c_0\umach\mu^{n_\ast} - \tfrac12 c_0\umach\mu^{n_\ast}
 = \tfrac12 c_0\umach\mu^{n_\ast} \ge 2^d\|P_\mu\|_\infty\eta .
\]
Hence $\|e_{n_\ast}\|_\infty\ge 2^d\eta$, and by \eqref{eq:Dnorm} we obtain
$D(u_{n_\ast})\ge\eta$, a contradiction. Therefore $D(u_m)\ge\eta$ for some
$m\le n_\ast$. Combining Parts~(i) and~(ii) establishes
$T_\ast=\lambda_\perp^{-1}\log(1/\umach)+O(1)$. \qed

\section{Proof of Proposition~\ref{prop:orbit}}\label{app:orbit}
Let $A\cong\Z_2^d$ with independent generators $\tau_1,\dots,\tau_d$, and let
$H\le\Gamma_A=A\times\{\pm1\}$. Consider the synthesis Walsh--Hadamard
transform $W_{\mathrm s}=W_{\mathrm s}^{(d)}\circ\cdots\circ W_{\mathrm s}^{(1)}:\R^{2^d}\to\R^{2^d}$
executed in the fixed stage order $i=1,\dots,d$, where stage $i$ maps each
pair $(v_g,v_{g\oplus e_i})$ ($g_i=0$) to $(v_g+v_{g\oplus e_i},v_g-v_{g\oplus e_i})$,
and the analysis transform $W_{\mathrm a}=W_{\mathrm a}^{(d)}\circ\cdots\circ W_{\mathrm a}^{(1)}$
with factors $1/2$ at each stage. Throughout, equalities between floating-point
zeros ignore the IEEE sign bit ($\pm 0.0$ both represent numerical zero), and
all intermediate values are assumed finite and within the domain of $F$.

\emph{Claim:} For every $h\in A$ and every finite vector $v\in\mathbb F^{2^d}$
of IEEE-754 floating-point numbers, the computed transforms satisfy the
\emph{exact floating-point intertwining identities} (up to signed zeros)
\begin{equation}\label{eq:wh_intertwine}
 \fl(W_{\mathrm s})(\chi(h)\odot v) = \Pi_h\,\fl(W_{\mathrm s})(v),\qquad
 \fl(W_{\mathrm a})(\Pi_h w) = \chi(h)\odot\fl(W_{\mathrm a})(w),
\end{equation}
where $(\chi(h)\odot v)_{\bsig}=\chi_{\bsig}(h)v_{\bsig}$ and
$(\Pi_h w)_g=w_{gh}$.

\emph{Proof of Claim.} Write $h=\tau_1^{h_1}\cdots\tau_d^{h_d}$ ($h_i\in\{0,1\}$)
and track the signed permutation stage by stage in the fixed execution order
$i=1,\dots,d$ (without assuming that floating-point stages along different axes
commute). Define $S_i$ on intermediate arrays after stage $i$ as the
composition of the index swaps $\Pi_{\tau_r}^{h_r}$ for $r\le i$ with the sign
multipliers $(-1)^{\sigma_r h_r}$ for $r>i$ (so $S_0=\chi(h)\odot$ and
$S_d=\Pi_h$). At stage $i$, we claim that
$\fl(W_{\mathrm s}^{(i)})\circ S_{i-1} = S_i\circ \fl(W_{\mathrm s}^{(i)})$:
\begin{itemize}
\item The index permutations $\Pi_{\tau_r}^{h_r}$ ($r<i$) act only on indices
  distinct from bit $i$ and therefore relabel the butterflies of stage $i$
  without altering any pair $(p,q)$ or arithmetic operation; similarly, the
  signs $(-1)^{\sigma_r h_r}$ ($r>i$) multiply both entries of each butterfly
  pair $(p,q)$ by the \emph{same} sign $\pm1$, and in IEEE round-to-nearest
  $\fl((-p)\pm(-q))=-\fl(p\pm q)$ numerically (up to signed zero when $p=\pm q$).
\item For bit $i$ itself: if $h_i=0$, nothing changes; if $h_i=1$, the sign
  $(-1)^{\sigma_i}$ maps $(p,q)\mapsto(p,-q)$, so the stage-$i$ butterfly
  computes $(\fl(p-q),\fl(p+q))$ instead of $(\fl(p+q),\fl(p-q))$, which is
  precisely the swap $\Pi_{\tau_i}$ applied after stage $i$.
\end{itemize}
Composing over $i=1,\dots,d$ gives
$\fl(W_{\mathrm s})\circ S_0 = S_d\circ\fl(W_{\mathrm s})$, proving the first
identity of \eqref{eq:wh_intertwine}. The identical stage-by-stage tracking
proves the second identity (noting that in IEEE round-to-nearest
$\fl((-a)/2)=-\fl(a/2)$ holds even in the subnormal range, and $\fl(0/2)=0$).

Now let $u\in V_{\Sigma_H}$ have exact floating-point zeros in all channels
$\bsig\notin\Sigma_H$. Then for every $(h,s)\in H$, $\chi_{\bsig}(h)=s$ on all
nonzero channels, so $\chi(h)\odot u = s\,u$. Applying \eqref{eq:wh_intertwine}
and $\fl(W_{\mathrm s})(s\,u)=s\,\fl(W_{\mathrm s})(u)$ gives
$\Pi_h\fl(W_{\mathrm s})(u)=s\,\fl(W_{\mathrm s})(u)$. Since $F$ is odd-exact
($\fl(F(-x))=-\fl(F(x))$), $w:=\fl(F(\fl(W_{\mathrm s})(u)))$ satisfies
$\Pi_h w = s\,w$. Applying the second identity of \eqref{eq:wh_intertwine}
yields $\chi(h)\odot\fl(W_{\mathrm a})(w)=s\,\fl(W_{\mathrm a})(w)$ for every
$(h,s)\in H$. For any $\bsig\notin\Sigma_H$, there exists $(h,s)\in H$ with
$\chi_{\bsig}(h)=-s$, whence $-\fl(W_{\mathrm a})(w)_{\bsig}=\fl(W_{\mathrm a})(w)_{\bsig}$,
forcing $\fl(W_{\mathrm a})(w)_{\bsig}=0$ (exact IEEE zero). \qed

\section{Proof of Theorem~\ref{thm:relative} and Corollary~\ref{cor:faithful}}
\label{app:relative}
We work with the polynomial channel algebra \eqref{eq:chalg}, assuming all
intermediate quantities remain finite and above the underflow threshold,
using the IEEE-754 model $\fl(a\circ b)=(a\circ b)(1+\xi)$, $|\xi|\le\umach$,
for $\circ\in\{+,-,\times\}$, and the standard normwise forward error bound
$\|\fl(\mathcal T v)-\mathcal Tv\|_\infty\le c_1\umach\log N\,\|\mathcal T\|_\infty\|v\|_\infty$
for an FFT-based channel transform $\mathcal T$ of length $O(N)$~\cite{Higham2002},
together with uniformly bounded stage multipliers on the active spectrum.

Fix a transverse channel $\bsig\notin\Sigma$. Consider a single evaluation of
the cubic nonlinearity \eqref{eq:chalg} on $u=u_\Sigma+u_\perp$ with
$\|u_\Sigma\|_\infty,\|u_\perp\|_\infty\le B$. Partition $\Z_2^d$ into
$\Sigma$ and $\Sigma^{\mathrm c}$. For each $\boldsymbol\rho\in\Z_2^d$:
\begin{itemize}
\item If $\bsig\oplus\boldsymbol\rho\notin\Sigma$, then
  $\|u_{\bsig\oplus\boldsymbol\rho}\|_\infty\le\|u_\perp\|_\infty$, while
  $\|(u^2)_{\boldsymbol\rho}\|_\infty\le 2^d B^2$ and its computed value
  satisfies
  $\|\fl((u^2)_{\boldsymbol\rho})-(u^2)_{\boldsymbol\rho}\|_\infty\le c_2\umach B^2$.
  Multiplying by $u_{\bsig\oplus\boldsymbol\rho}$ gives error at most
  $(c_2+2^d)\umach B^2\|u_\perp\|_\infty$.
\item If $\bsig\oplus\boldsymbol\rho\in\Sigma$, then by \eqref{eq:closed} every
  summand $u_{\mathbf a}u_{\mathbf a\oplus\boldsymbol\rho}$ in
  $(u^2)_{\boldsymbol\rho}$ has at least one factor outside $\Sigma$, so each
  product has modulus at most $B\|u_\perp\|_\infty$ and rounding error at most
  $\umach B\|u_\perp\|_\infty$. Summing the $2^d$ products yields
  $\|\fl((u^2)_{\boldsymbol\rho})-(u^2)_{\boldsymbol\rho}\|_\infty\le c_3\umach B\|u_\perp\|_\infty$,
  and multiplying by $u_{\bsig\oplus\boldsymbol\rho}$ (bounded by $B$) gives
  error at most $c_4\umach B^2\|u_\perp\|_\infty$.
\end{itemize}
Summing the $2^d$ outer products introduces an additional error bounded by
$2^d\umach\sum_{\boldsymbol\rho}|(u^2)_{\boldsymbol\rho}u_{\bsig\oplus\boldsymbol\rho}|\le 4^d\umach B^2\|u_\perp\|_\infty$.
Thus $\|\fl((u^3)_{\bsig})-(u^3)_{\bsig}\|_\infty\le C_d\umach B^2\|u_\perp\|_\infty$.
Next, the forward and inverse transforms and radial multipliers act
channel-wise, so on channel $\bsig\notin\Sigma$ their rounding errors are
bounded by $c_1\umach\log N$ times the norm of channel $\bsig$ itself. Now let
$\delta_{\perp}:=\max_{0\le r<k}\|u^{n-r}_\perp\|_\infty$ be the maximum
transverse norm across the $k$ history inputs $\mathbf U=(u^n,\dots,u^{n-k+1})$
(with $k=1$ for one-step methods), and let
$v^{(j)}=v^{(j)}_\Sigma+v^{(j)}_\perp$ and
$\tilde v^{(j)}=\tilde v^{(j)}_\Sigma+\tilde v^{(j)}_\perp$ be the exact and
computed inputs to nonlinear stage $j\in\{1,\dots,q\}$, all bounded in
$\|\cdot\|_\infty$ by the common bound $B$ (along with all history inputs
$u^{n-r}$). For a $k$-step scheme whose linear combination includes cached
nonlinearities $N(u^{n-r})$ ($r\ge 1$), each cached evaluation is formed from
$u^{n-r}$ by the same channel cubic \eqref{eq:chalg}, so its exact transverse
norm is $O((1+B)^2\delta_\perp)$ and its rounding error is
$O(\umach(1+B)^2\delta_\perp)$ (counting cached evaluations among the
bounded history inputs, or equivalently taking $q\ge 1$ to cover the nonlinear
evaluations entering the update). We track by induction on $j$ the exact
transverse norm $\|v^{(j)}_\perp\|_\infty\le M_j\delta_\perp$, the tangential
stage error $\|\tilde v^{(j)}_\Sigma-v^{(j)}_\Sigma\|_\infty\le A_j\umach$,
and the transverse stage error
$\|\tilde v^{(j)}_\perp-v^{(j)}_\perp\|_\infty\le E_j\umach\delta_\perp$,
starting from $M_1=O(1)$, $A_1=0$, $E_1=0$ (for linear combinations of the $k$
exact history inputs, adding at most $c_{\mathrm h}\umach\delta_\perp$ in
combining the history states). At stage $j$:
\begin{itemize}
\item On any transverse channel $\bsig\notin\Sigma$, split
  $\fl((\tilde v^{(j)})^3)_{\bsig}-((v^{(j)})^3)_{\bsig}
  =\big[\fl((\tilde v^{(j)})^3)_{\bsig}-((\tilde v^{(j)})^3)_{\bsig}\big]
   +\big[((\tilde v^{(j)})^3)_{\bsig}-((v^{(j)})^3)_{\bsig}\big]$.
  The first bracket is bounded by
  $C_d\umach B^2\|\tilde v^{(j)}_\perp\|_\infty\le C_d\umach B^2(M_j+E_j\umach)\delta_\perp$.
  In the second bracket, by \eqref{eq:closed} every channel product in
  $(\tilde v^{(j)})^3-(v^{(j)})^3$ for $\bsig\notin\Sigma$ contains at least
  one transverse factor, so it is bounded by
  $C_d'\big(B^2\|\tilde v^{(j)}_\perp-v^{(j)}_\perp\|_\infty + B\|\tilde v^{(j)}_\Sigma-v^{(j)}_\Sigma\|_\infty(\|\tilde v^{(j)}_\perp\|_\infty+\|v^{(j)}_\perp\|_\infty)\big)
  \le C_d''\big(B^2 E_j + B A_j M_j\big)\umach\delta_\perp$.
\item Applying the bounded channel transforms, stage multipliers, and linear
  combinations with history states propagates
  $M_{j+1}\le c_0(1+B)^2 M_j$,
  $A_{j+1}\le c_0(1+B)^2 A_j + c_0'(1+B)^3$, and
  $E_{j+1}\le c_0(1+B)^2 E_j + c_0''\big((1+B)^2 M_j + B A_j M_j\big)$.
\end{itemize}
By induction, $M_j=O((1+B)^{2j-2})$, $A_j=O((1+B)^{2j-1})$ for $j\ge 2$
($A_1=0$), and hence $B A_j M_j=O((1+B)^{4j-2})$ dominates the recurrence,
yielding $E_{q+1}\le C(1+B)^{4q-2}$ (reducing to $C(1+B)^2$ for $q=1$), which
proves \eqref{eq:relbound}. (By contrast, Walsh--Hadamard orbit evaluation mixes
$u_\Sigma$ and $u_\perp$ into $u_\Sigma\pm u_\perp$, so sub-roundoff
$\|u_\perp\|\ll\umach\|u_\Sigma\|$ is absorbed during synthesis; hence
Theorem~\ref{thm:relative} applies specifically to the polynomial channel
algebra \eqref{eq:chalg}.)

For Corollary~\ref{cor:faithful}, write the exact perturbed discrete trajectory
as $u_m^{\mathrm{ex}}=s_m^{\mathrm{ex}}+e_m^{\mathrm{ex}}$ and the computed
trajectory as $\tilde u_m=\tilde s_m+\tilde e_m$, with
$\tilde w_m:=\tilde s_m-s_m^{\mathrm{ex}}\in V_\Sigma$ satisfying
$\|\tilde w_m\|_\infty\le C_w m\umach$, and set
$J_m^\perp:=\partial_e\Phi_\tau^\perp(s_m^{\mathrm{ex}},0)$ along the exact
perturbed tangential path $s_m^{\mathrm{ex}}$. Set
$d_m:=\tilde e_m-e_m^{\mathrm{ex}}$ ($d_0=0$) and define the running
maximum relative error
\[
 r_k:=\max_{0\le m\le k}\frac{\|d_m\|_\infty}{\|e_m^{\mathrm{ex}}\|_\infty}.
\]
We prove by induction that $r_k\le 1$ (in fact $r_k\le M_1 k^2\umach\le 1$) on
$0\le k\le n$, which automatically guarantees
$\|\tilde e_m\|_\infty\le 2\|e_m^{\mathrm{ex}}\|_\infty\le 2\eta_0$ on the
entire window. By exact-arithmetic equivariance,
$\Phi_\tau^\perp(s_m^{\mathrm{ex}},0)=\Phi_\tau^\perp(\tilde s_m,0)=0$, so
expanding $\Phi_\tau^\perp(s_m^{\mathrm{ex}},e_m^{\mathrm{ex}})$ around
$(s_m^{\mathrm{ex}},0)$ and $\Phi_\tau^\perp(\tilde s_m,\tilde e_m)$ around
$(\tilde s_m,0)$ is \emph{exact} (with no omitted tangential feedback): on any
step $m<k$ with $r_{k-1}\le 1$, Theorem~\ref{thm:relative} and the smoothness
of $\Phi_\tau$ yield
\begin{align*}
 e_{m+1}^{\mathrm{ex}} &= J_m^\perp e_m^{\mathrm{ex}} + q_m^{\mathrm{ex}}(e_m^{\mathrm{ex}}), \\
 \tilde e_{m+1} &= J_m^\perp\tilde e_m + \Delta J_m^\perp\tilde e_m + \rho_m + \tilde q_m(\tilde e_m),
\end{align*}
where $q_m^{\mathrm{ex}}(e):=\Phi_\tau^\perp(s_m^{\mathrm{ex}},e)-J_m^\perp e$,
$\tilde q_m(e):=\Phi_\tau^\perp(\tilde s_m,e)-\partial_e\Phi_\tau^\perp(\tilde s_m,0)e$,
$\|\Delta J_m^\perp\|_\infty:=\|\partial_e\Phi_\tau^\perp(\tilde s_m,0)-\partial_e\Phi_\tau^\perp(s_m^{\mathrm{ex}},0)\|_\infty\le C_J\|\tilde w_m\|_\infty\le C_J C_w m\umach$,
$\|\rho_m\|_\infty\le C\umach\|\tilde e_m\|_\infty$, and by smooth dependence
on $(s,e)$ the remainder difference satisfies
\[
 \|\tilde q_m(\tilde e_m)-q_m^{\mathrm{ex}}(e_m^{\mathrm{ex}})\|_\infty
 \le C_q\big(\|\tilde e_m\|_\infty+\|e_m^{\mathrm{ex}}\|_\infty\big)\|d_m\|_\infty
   + C_q\|\tilde w_m\|_\infty\|\tilde e_m\|_\infty^2 .
\]
Absorbing $C_q\|\tilde w_m\|_\infty\|\tilde e_m\|_\infty^2\le 2C_q\eta_0 C_w m\umach\|\tilde e_m\|_\infty$
into $C_1(m+1)\umach\|\tilde e_m\|_\infty$ (with $C_1:=C+(C_J+2C_q\eta_0)C_w$)
and using $\|\tilde e_m\|_\infty+\|e_m^{\mathrm{ex}}\|_\infty\le(2+r_{k-1})\|e_m^{\mathrm{ex}}\|_\infty\le 3\|e_m^{\mathrm{ex}}\|_\infty$,
subtracting gives $d_{m+1}=J_m^\perp d_m+\xi_m$ with
\[
 \|\xi_m\|_\infty \le C_1(m+1)\umach(1+r_{k-1})\|e_m^{\mathrm{ex}}\|_\infty
                    + 3 C_q\|e_m^{\mathrm{ex}}\|_\infty^2 r_{k-1}.
\] Applying discrete
variation of constants to $d_k=\sum_{m=0}^{k-1}\Psi_{k,m+1}\xi_m$ and dividing
by $\|e_k^{\mathrm{ex}}\|_\infty$:

\emph{Part~(i) (general compatibility):} Using \eqref{eq:propseed} and
$\|e_m^{\mathrm{ex}}\|_\infty\le\eta_0$,
\begin{align*}
 \frac{\|d_k\|_\infty}{\|e_k^{\mathrm{ex}}\|_\infty}
 &\le \sum_{m=0}^{k-1}\frac{\|\Psi_{k,m+1}\|_\infty\|e_m^{\mathrm{ex}}\|_\infty}{\|e_k^{\mathrm{ex}}\|_\infty}
      \Big(C_1(m+1)\umach(1+r_{k-1}) + 3C_q\eta_0 r_{k-1}\Big) \\
 &\le \Lambda\Big(C_1\frac{k(k+1)}{2}\umach(1+r_{k-1}) + 3C_q k\eta_0 r_{k-1}\Big).
\end{align*}
Whenever $\Lambda\big(C_1 n(n+1)\umach + 6C_q n\eta_0\big)\le\tfrac12$ (and
$\Lambda C_1 n(n+1)\umach\le\tfrac12$), induction gives
$r_k\le \Lambda C_1 k(k+1)\umach\le 1$ for all $k\le n$.

\emph{Part~(ii) (exponentially growing regime):} Now suppose
$\|\Psi_{k,m+1}\|_\infty\le K\mu^{k-m-1}$ and
$c_w\delta\mu^m\le\|e_m^{\mathrm{ex}}\|_\infty\le C_w\delta\mu^m$ on
$0\le m\le n$ with $\mu=\ee^{\lambda_\perp\tau}>1$. Retaining the factor
$\|e_m^{\mathrm{ex}}\|_\infty^2\le C_w^2\delta^2\mu^{2m}$ in the quadratic term
yields a genuine geometric sum:
\begin{align*}
 \sum_{m=0}^{k-1}\frac{\|\Psi_{k,m+1}\|_\infty\|e_m^{\mathrm{ex}}\|_\infty^2}{\|e_k^{\mathrm{ex}}\|_\infty}
 &\le \sum_{m=0}^{k-1}\frac{K\mu^{k-m-1}(C_w\delta\mu^m)^2}{c_w\delta\mu^k}
  = \frac{K C_w^2}{c_w}\delta\sum_{m=0}^{k-1}\mu^{m-1} \\
 &\le \frac{K C_w^2}{c_w\mu(\mu-1)}\,\delta\mu^k
  \le \frac{K C_w^2}{c_w^2\mu(\mu-1)}\,\|e_k^{\mathrm{ex}}\|_\infty
  \le \frac{K C_w^2}{c_w^2\mu(\mu-1)}\,\eta_0 .
\end{align*}
Meanwhile, the rounding/tangential sum satisfies
\[
 \sum_{m=0}^{k-1}\frac{\|\Psi_{k,m+1}\|_\infty\|e_m^{\mathrm{ex}}\|_\infty}{\|e_k^{\mathrm{ex}}\|_\infty}\,C_1(m+1)\umach(1+r_{k-1})
 \le \frac{K C_w C_1}{c_w\mu}\,\frac{k(k+1)}{2}\umach(1+r_{k-1}).
\]
Combining both bounds, provided $\eta_0$ is small enough (independently of $n$)
and $n^2\umach$ is small enough that
\[
 \frac{3K C_q C_w^2\eta_0}{c_w^2\mu(\mu-1)}\le\frac14
 \qquad\text{and}\qquad
 \frac{K C_w C_1 n(n+1)\umach}{2c_w\mu}\le\frac14,
\]
we obtain
\[
 \frac{\|d_k\|_\infty}{\|e_k^{\mathrm{ex}}\|_\infty}
 \le \frac{K C_w C_1 k(k+1)}{2c_w\mu}\umach + \Big(\frac14+\frac14\Big)r_{k-1},
\]
whence by induction $r_k\le M_1 k(k+1)\umach\le\tfrac12$ with
$M_1:=K C_w C_1/(c_w\mu)$ for all $k\le n$, proving \eqref{eq:rel_exp}.

\emph{Part~(iii) (first-crossing times):} Around a steady symmetric state
$u^\ast$ with tangential bound $\|\Psi_{k,m}^\Sigma\|_\infty\le K_\Sigma\mu_\Sigma^{k-m}$
($\mu_\Sigma<\mu^2$), the joint bootstrap of Proposition~\ref{prop:obstruction}
gives both $\|s_m^{\mathrm{ex}}-u^\ast\|_\infty\le C_{\mathrm{fb}}\|e_m^{\mathrm{ex}}\|_\infty^2$
and $c_w\delta\mu^m\le\|e_m^{\mathrm{ex}}\|_\infty\le C_w\delta\mu^m$ (so
$\|J_m^\perp-J_\ast^\perp\|_\infty=O(\|e_m^{\mathrm{ex}}\|_\infty^2)$ and the
exact perturbed propagator $\Psi_{k,m}$ inherits the bound
$\|\Psi_{k,m}\|_\infty\le K'\mu^{k-m}$ on the window
$\|e_m^{\mathrm{ex}}\|_\infty\le\eta_0$ for small $\eta_0$). Furthermore,
$e_m^{\mathrm{ex}}=\delta\mu^m\big(P_\mu w + O(\rho^m + \eta_0)\big)$ (where
$\rho<1$ is the subdominant spectral ratio of $J_\ast^\perp$). Fix $\alpha\in(0,\tfrac12(1-\mu^{-1}))$ and set $q_-:=\mu(1-\alpha)>1$,
$q_+:=\mu(1+\alpha)$. Because $e_m^{\mathrm{ex}}=\delta\mu^m(P_\mu w+O(\rho^m+\eta_0))$
with $\rho<1$ and $D(P_\mu w)>0$, there exists an index $m_0=m_0(\alpha)$
(independent of $\delta$ and $\umach$) such that for $\eta_0$ small enough and
every $m\ge m_0$ with $\|e_m^{\mathrm{ex}}\|_\infty\le\eta_0$, the exact
defect obeys the dominant-mode two-sided step-growth bounds
\[
 q_- D(u_m^{\mathrm{ex}}) \le D(u_{m+1}^{\mathrm{ex}}) \le q_+ D(u_m^{\mathrm{ex}}).
\]
Set $\eta_0:=4\cdot 2^d q_+\eta$ and define the $\umach$-independent cutoff
index
\[
 N_+:=\min\{m\ge 1:\ D(u_m^{\mathrm{ex}})\ge 2\eta\}.
\]
Take $\delta$ small enough that $N_+-1\ge m_0$ (and $n_--1\ge m_0$, where
$n_-$ is defined below). For every $m<N_+$, $D(u_m^{\mathrm{ex}})<2\eta$, so
$\|e_m^{\mathrm{ex}}\|_\infty\le 2^d D(u_m^{\mathrm{ex}})<2^{d+1}\eta<\eta_0$.
At $m=N_+$, because $N_+-1\ge m_0$ and $\|e_{N_+-1}^{\mathrm{ex}}\|_\infty<\eta_0$,
the dominant-mode upper growth bound applies directly to $N_+-1$ (even if
$D(u_{N_+-1}^{\mathrm{ex}})<\tfrac12\eta$), yielding
$D(u_{N_+}^{\mathrm{ex}})\le q_+ D(u_{N_+-1}^{\mathrm{ex}})<2q_+\eta$ and
hence $\|e_m^{\mathrm{ex}}\|_\infty\le 2^d D(u_m^{\mathrm{ex}})<2\cdot 2^d q_+\eta<\eta_0$
for \emph{all} $0\le m\le N_+$, including $N_+-1$ and all predecessors of
crossing. Now set $\epsilon:=2^d C' N_+^2\umach$ (non-circularly, since $N_+$
does not depend on $\umach$) and take $\umach$ small enough that
$\epsilon\le\tfrac14$. By Part~(ii) and \eqref{eq:Dnorm},
$|D(\tilde u_m)-D(u_m^{\mathrm{ex}})|\le\|d_m\|_\infty\le C' m^2\umach\|e_m^{\mathrm{ex}}\|_\infty\le\epsilon\,D(u_m^{\mathrm{ex}})$,
giving $(1-\epsilon)D(u_m^{\mathrm{ex}})\le D(\tilde u_m)\le(1+\epsilon)D(u_m^{\mathrm{ex}})$
for all $0\le m\le N_+$. Let $n_-$ and $n_+$ be the first indices at which
$D(u_m^{\mathrm{ex}})\ge\eta/(1+\epsilon)$ and
$D(u_m^{\mathrm{ex}})\ge\eta/(1-\epsilon)$, respectively. Because
$\eta/(1-\epsilon)\le\tfrac43\eta<2\eta$, we have
$m_0\le n_-\le n_{\mathrm{ex}}\le n_+\le N_+$, and the computed crossing index
$n_{\mathrm{comp}}$ of $D(\tilde u_m)\ge\eta$ is sandwiched between
$n_-\le n_{\mathrm{comp}}\le n_+$. By the dominant-mode lower growth bound
for $m\ge n_-\ge m_0$, $D(u_{n_-+j}^{\mathrm{ex}})\ge q_-^j\eta/(1+\epsilon)$,
so
\[
 0 \le n_+ - n_- \le \left\lceil\frac{\log\big((1+\epsilon)/(1-\epsilon)\big)}{\log q_-}\right\rceil
   = O\!\Big(1 + \frac{\epsilon}{\lambda_\perp\tau}\Big),
\]
whence $|T_\ast(\delta)-T_\ast^{\mathrm{ex}}(\delta)|=O\big(\tau+N_+^2\umach/\lambda_\perp\big)$,
valid throughout the IEEE double normal range ($\gtrsim 2.2\times 10^{-308}$). \qed

\section{Proof of Theorem~\ref{thm:cfism}}\label{app:cfism}
(i) We treat the two reflection modes of Algorithm~\ref{alg:cfism} in turn.
In the point-odd case (\texttt{mode=\allowbreak'point'}, covering $d=1$ and 2D
$C_2$), $v_0(-\mathbf x)=-v_0(\mathbf x)$ implies that every antipodal pair
sum $r_j(\mathbf x_0):=(u_0(\mathbf x_0+\tilde{\mathbf x}_j)+u_0(\mathbf x_0-\tilde{\mathbf x}_j))/s_{\mathrm{var}}$
vanishes identically at $\mathbf c=\mathbf x_0$. In the \emph{coordinate-wise
odd mode} (\texttt{mode='refl'} in $d\ge 2$, e.g.\ the 2D odd--odd X-junction),
$v_0(R_i\mathbf x)=-v_0(\mathbf x)$ for each axis $i=1,\dots,d$ (where $R_i$
flips only coordinate $i$; note that flipping two axes simultaneously leaves
$v_0$ unchanged, so antipodal sums do \emph{not} vanish), and therefore every
per-axis residual
$r_{i,j}(\mathbf x_0):=(u_0(\mathbf x_0+\tilde{\mathbf x}_j)+u_0(\mathbf x_0+R_i\tilde{\mathbf x}_j))/s_{\mathrm{var}}$
vanishes identically at $\mathbf c=\mathbf x_0$. In both modes,
$\bar r(\mathbf x_0)=0$, $\widetilde{\mathcal J}_{\mathrm c}(\mathbf x_0)=0$,
and $D\widetilde{\mathcal J}_{\mathrm c}(\mathbf x_0)=\mathbf 0$. Under the
assumption that $D^2\widetilde{\mathcal J}_{\mathrm c}(\mathbf x_0)\succ0$ and
that the coarse FFT minimizer $\mathbf c^{(0)}$ lies in the local basin of
attraction of Gauss--Newton for $\widetilde{\mathcal J}_{\mathrm c}$ with
approximate Jacobian $\widetilde J(\mathbf c)$ satisfying
$\|\widetilde J(\mathbf c)-J(\mathbf c)\|\le\theta_J$ sufficiently small, the
zero-residual property $r_{\mathrm c}(\mathbf x_0)=\mathbf 0$ ensures that
inexact Gauss--Newton contracts at each step as
\[
 \|\mathbf c_{k+1}-\mathbf x_0\|
 \le C\big(\|\widetilde J-J\|\,\|\mathbf c_k-\mathbf x_0\|+\|\mathbf c_k-\mathbf x_0\|^2\big)
   + O(10^{-d_0}).
\]
With central finite differences ($h_q=10^{-\lfloor d_0/2\rfloor}$ in
$d_0$-digit arithmetic, where cancellation yields
$\|\widetilde J-J\|=O(10^{-d_0/2})$), starting inside the quadratic-contraction
region gives quadratic doubling until $O(10^{-d_0/2})$ followed by an
$O(10^{-d_0/2})$ contraction, reaching $\|\mathbf c_q-\mathbf x_0\|\le C_1 10^{-d_0}$
in $O(\log d_0)$ iterations (with the actual step count depending on the
starting accuracy and Jacobian condition number: under stopping criterion
$\|\Delta\mathbf c\|_\infty<10^{-(d_0-6)}$ and a 14-iteration cap in
\texttt{cf\_gauge\_init.m}, Figure~\ref{fig:cfism_horizon}'s 1D mixed-harmonic
datum starting from double-precision polish converges in 1--2 VPA iterations
at $d_0=320$, while the 2D point-odd branch starting directly from the coarse
FFT grid takes 7--8 iterations; if the user callback $u_0$ itself suffers
internal cancellation limiting attainable accuracy below $10^{-(d_0-6)}$, as
in $(1+10^{-20}\sin(x+0.2)+3\times 10^{-21})-1$, the loop exits at the
14-iteration cap with \texttt{converged=false} while still recovering the
centre to full attainable precision $\approx 10^{-309}$ and exact
double-precision channels; a fixed tolerance $\theta_J$ independent of $d_0$
gives linear contraction requiring $O(d_0/\log(1/\theta_J))$ iterations).
Evaluating the canonical-frame channels in $d_0$-digit arithmetic via
$2^d$-point reflection butterflies about $\mathbf c_q$---which form linear
combinations of opposite pairs vanishing at $\mathbf c=\mathbf x_0$ for every
forbidden channel (even channels in \texttt{mode='point'}, and all channels
$\bsig\neq(1,\dots,1)$ in \texttt{mode='refl'})---gives
$\|v^0_\perp\|_\infty\le C\|Dv_0\|_\infty\|\mathbf c_q-\mathbf x_0\|+C_2 10^{-d_0}\le C_0 10^{-d_0}$.
Rounding $v^0_\perp$ to IEEE double precision preserves this bound (up to
relative error $\umach$) in the IEEE normal range.

To bound the subsequent pure double-precision evolution from the nonzero
initial seed $e_0=v^0_\perp$ (with $\|e_0\|_\infty\le C_0 10^{-d_0}$), we
include the initial propagator term $\Psi_{n,0}e_0$ in discrete variation of
constants along the exact perturbed tangential path (with tangential feedback
$O(\|e_m\|_\infty^2)$ controlled by $\|\Psi_{k,m}^\Sigma\|_\infty\le K_\Sigma\mu_\Sigma^{k-m}$,
$\mu_\Sigma<\mu^2$, as in Proposition~\ref{prop:obstruction}, and tangential
rounding error $\|\tilde w_m\|_\infty\le C_w m\umach$): by
Theorem~\ref{thm:relative} (polynomial channel algebra, normal range),
\[
 e_n = \Psi_{n,0}e_0 + \sum_{m=0}^{n-1}\Psi_{n,m+1}\big(\Delta J_m^\perp e_m + \rho_m + q_m(0,e_m)\big),
\]
with $\|\Delta J_m^\perp\|_\infty\le C_J C_w m\umach$,
$\|\rho_m\|_\infty\le C\umach\|e_m\|_\infty$, and
$\|q_m(0,e_m)\|_\infty\le C_q\|e_m\|_\infty^2$. Setting
$\beta_n:=\max_{0\le m\le n}\mu^{-m}\|e_m\|_\infty$ ($\mu=\ee^{\lambda_\perp\tau}>1$)
and $T_1:=\lambda_\perp^{-1}\log\big(\eta/(2K C_0 10^{-d_0})\big)$, the base
case $\beta_0=\|e_0\|_\infty\le 2K C_0 10^{-d_0}$ holds. For $1\le k\le n$ with
$n\tau\le T_1$ and induction hypothesis $\beta_{n-1}\le 2K C_0 10^{-d_0}$,
multiplying by $\mu^{-k}$ and noting that $\mu^{-m-1}\mu^m=\mu^{-1}$ on the
linear terms while $\sum_{m=0}^{k-1}\mu^{m-1}<\mu^{k-1}/(\mu-1)$ on the
quadratic term yields
\begin{align*}
 \mu^{-k}\|e_k\|_\infty
 &\le K\|e_0\|_\infty
    + K\sum_{m=0}^{k-1}\mu^{-m-1}\Big(C\umach\beta_{n-1}\mu^m
      + C_J C_w m\umach\beta_{n-1}\mu^m + C_q\beta_{n-1}^2\mu^{2m}\Big) \\
 &\le K\|e_0\|_\infty
    + \Big[\frac{K}{\mu}\Big(C k\umach + \frac{C_J C_w k(k-1)\umach}{2}\Big)
         + \frac{K C_q\eta}{\mu(\mu-1)}\Big]\beta_{n-1} \\
 &\le K C_0 10^{-d_0} + \frac12\beta_{n-1} \le 2K C_0 10^{-d_0}
\end{align*}
provided
\[
 \frac{K}{\mu}\Big(C n\umach + \frac{C_J C_w n(n-1)\umach}{2}\Big)
 + \frac{K C_q\eta}{\mu(\mu-1)} \le \frac12
\]
for $n\le T_1/\tau$. Hence $\|e_n\|_\infty\le 2K C_0 10^{-d_0}\mu^n\le\eta$ and
$D(u^n)\le\|e_n\|_\infty\le\eta$ for all $n\tau\le T_1$, proving
$T_\ast(d_0)\ge\lambda_\perp^{-1}d_0\ln 10+O(1)$. If in addition
$J_n^\perp\equiv J^\perp$ has a simple dominant unstable eigenvalue
$\mu=\ee^{\lambda_\perp\tau}>1$ with non-degenerate seed ratio
$\|P_\mu v^0_\perp\|_\infty\ge c_0\|v^0_\perp\|_\infty$ and
$\|v^0_\perp\|_\infty\asymp 10^{-d_0}$, Corollary~\ref{cor:faithful}(iii)
yields the matching upper bound $T_\ast(d_0)=\lambda_\perp^{-1}d_0\ln 10+O(1)$.

(ii) Since the smallest positive IEEE-754 double subnormal is
$2^{-1074}\approx 4.94\times 10^{-324}$, the round-to-nearest threshold to
zero is $2^{-1075}\approx 2.47\times 10^{-324}$. Whenever
$C_0 10^{-d_0}<2^{-1075}$, every entry of $v^0_\perp$ rounds to exact IEEE
zero \texttt{0.0} upon conversion to double precision. Theorem~\ref{thm:exact}
then guarantees $D(u^n)\equiv 0$ for all $n\ge 0$ (assuming finite
intermediates).

(iii) If $u_0(\mathbf x)=v_0(\mathbf x-\mathbf x_0)+\delta$ with $\delta$
constant, then for every $\mathbf c$ and every point $\tilde{\mathbf x}_j$ (in
either \texttt{mode='point'} or \texttt{mode='refl'}), each pair sum shifts by
the same constant $2\delta$, which cancels identically upon subtracting the
sample mean $\bar r(\mathbf c)$ (and likewise cancels from the
finite-difference Jacobian). Under raw amplitude scaling
$s_\delta=\|u_0\|_\infty$, the normalized objective and its Hessian scale by
the positive factor $(s_0/s_\delta)^2$ (preserving $\mathbf x_0$ and
positive-definiteness); under the mean-removed variation scaling
$s_{\mathrm{var}}=\|u_0-\bar u_0\|_\infty$ used in \texttt{cf\_gauge\_init.m}
(with $d_0$-digit VPA coarse fallback when double-precision variation vanishes
under a dominant constant $\delta$), $s_{\mathrm{var}}$ itself is invariant to
$\delta$, leaving $\widetilde{\mathcal J}_{\mathrm c}$ and
$D^2\widetilde{\mathcal J}_{\mathrm c}(\mathbf x_0)\succ0$ unchanged. Provided
the smallest singular value of the normalized Jacobian exceeds the rank
tolerance, Gauss--Newton recovers $\mathbf c_q=\mathbf x_0+O(10^{-d_0})$, and
evaluating the unnormalized channels in $d_0$-digit arithmetic yields
$v^0_\perp=\delta+O(10^{-d_0})$, which recovers $\delta$ to full IEEE double
precision once $\delta$ lies in the normal range and
$C_0 10^{-d_0}\ll\umach|\delta|$ (otherwise, if the coarse/Jacobian rank check
triggers degenerate fallback, the routine retains a finite origin and full
high-precision channels without discarding $\delta$). Corollary~\ref{cor:faithful}
then applies directly under its propagator--seed hypotheses. \qed

\end{document}